%
%
%
\def\PHHai{Ph\`ung H\^o Hai }
\def\emph#1{{\it#1}}
%
%
\hsize=5in
\baselineskip=12pt
\vsize=20cm
\parindent=10pt
\pretolerance=40
\predisplaypenalty=0
\displaywidowpenalty=0
\finalhyphendemerits=0
\hfuzz=2pt
\frenchspacing
\footline={\ifnum\pageno=1\else\hfil\tenrm\number\pageno\hfil\fi}
%
%
\input amssym.def
\font\tenbbold=bbold10
\font\sevenbbold=bbold7
\newfam\bbfam
\textfont\bbfam=\tenbbold
\scriptfont\bbfam=\sevenbbold
\def\titlefonts{\baselineskip=1.44\baselineskip
	\font\titlef=cmbx12
	\titlef
	}
\font\tenib=cmmib10 
\font\sevenib=cmmib7
\skewchar\tenib='177
\skewchar\sevenib='177
\font\tenbsy=cmbsy10
\skewchar\tenbsy='60
\def\boldfonts{\bf
	\textfont1=\tenib
	\scriptfont1=\sevenib
	\textfont2=\tenbsy
	}
\font\ninerm=cmr9
\font\ninebf=cmbx9
\font\ninei=cmmi9
\skewchar\ninei='177
\font\ninesy=cmsy9
\skewchar\ninesy='60
\font\nineit=cmti9
\def\reffonts{\baselineskip=0.9\baselineskip
	\textfont0=\ninerm
	\def\rm{\fam0\ninerm}%
	\textfont1=\ninei
	\textfont2=\ninesy
	\textfont\bffam=\ninebf
	\def\bf{\fam\bffam\ninebf}%
	\def\it{\nineit}%
	}
%
%
\def\frontmatter{\vbox{}\vskip1cm\bgroup
	\leftskip=0pt plus1fil\rightskip=0pt plus1fil
	\parindent=0pt
	\parfillskip=0pt
	\pretolerance=10000
	}
\def\endfrontmatter{\egroup\bigskip}
\def\title#1{{\titlefonts#1\par}}
\def\author#1{\bigskip#1\par}
\def\address#1{\medskip{\reffonts\it#1}}
\def\email#1{\smallskip{\reffonts{\it E-mail: }\rm#1}}
\def\thanks#1{\footnote{}{\reffonts\rm\noindent#1\hfil}}
\def\section#1\par{\ifdim\lastskip<\bigskipamount\removelastskip\fi
	\penalty-250\bigskip
	\vbox{\leftskip=0pt plus1fil\rightskip=0pt plus1fil
	\parindent=0pt
	\parfillskip=0pt
	\pretolerance=10000{\boldfonts#1}}\nobreak\medskip
	}
\def\proclaim#1. {\medbreak\bgroup{\noindent\bf#1.}\ \it}
\def\endproclaim{\egroup
	\ifdim\lastskip<\medskipamount\removelastskip\medskip\fi}
\newdimen\itemsize
\def\setitemsize#1 {{\setbox0\hbox{#1\ }
	\global\itemsize=\wd0}}
\def\item#1 #2\par{\ifdim\lastskip<\smallskipamount\removelastskip\smallskip\fi
	{\leftskip=\itemsize
	\noindent\hskip-\leftskip
	\hbox to\leftskip{\hfil\rm#1\ }#2\par}\smallskip}
\def\Proof#1. {\ifdim\lastskip<\medskipamount\removelastskip\medskip\fi
	{\noindent\it Proof\if\space#1\space\else\ \fi#1.}\ }
\def\endproof{\hfill\hbox{}\quad\hbox{}\hfill\llap{$\square$}\medskip}
%
%
\newcount\citation
\newtoks\citetoks
\def\citedef#1\endcitedef{\citetoks={#1\endcitedef}}
\def\endcitedef#1\endcitedef{}
\def\citenum#1{\citation=0\def\curcite{#1}%
	\expandafter\checkendcite\the\citetoks}
\def\checkendcite#1{\ifx\endcitedef#1?\else
	\expandafter\lookcite\expandafter#1\fi}
\def\lookcite#1 {\advance\citation by1\def\auxcite{#1}%
	\ifx\auxcite\curcite\the\citation\expandafter\endcitedef\else
	\expandafter\checkendcite\fi}
\def\cite#1{\makecite#1,\cite}
\def\makecite#1,#2{[\citenum{#1}\ifx\cite#2]\else\expandafter\clearcite\expandafter#2\fi}
\def\clearcite#1,\cite{, #1]}
%
%
\def\references{\section References\par
	\bgroup
	\parindent=0pt
	\reffonts
	\rm
	\frenchspacing
	\setbox0\hbox{99. }\leftskip=\wd0
	}
\def\endreferences{\egroup}
\newtoks\authtoks
\newif\iffirstauth
\def\checkendauth#1{\ifx\auth#1%
    \iffirstauth\the\authtoks
    \else{} and \the\authtoks\fi,%
  \else\iffirstauth\the\authtoks\firstauthfalse
    \else, \the\authtoks\fi
    \expandafter\nextauth\expandafter#1\fi
}
\def\nextauth#1,#2;{\authtoks={#1 #2}\checkendauth}
\def\auth#1{\nextauth#1;\auth}
\newif\ifbookref
\def\nextref#1 {\par\hskip-\leftskip
	\hbox to\leftskip{\hfil\citenum{#1}.\ }%
	\initnextref}
\def\initnextref{\bookreffalse\firstauthtrue\ignorespaces}
\def\paper#1{{\it#1},}
\def\book#1{\bookreftrue{\it#1},}
\def\journal#1{#1}
\def\Vol#1{{\bf#1}}
\def\publisher#1{#1,}
\def\Year#1{\ifbookref #1.\else(#1)\fi}
\def\Pages#1{\makepages#1.}
\long\def\makepages#1-#2.#3{#1--#2\ifx\par#3.\fi#3}
\def\inRus{{ \rm(in Russian)}}
\def\etransl#1{English translation in \journal{#1}}
%
%
\newsymbol\smallsetminus 2272
\let\setm\smallsetminus
\newsymbol\square 1203
\newsymbol\triangledown 104F
\let\tri\triangledown
\def\bbk{{\fam\bbfam k}}
\let\Rar\Rightarrow
\let\hrar\hookrightarrow
\let\ot\otimes
\let\sbs\subset
\let\sps\supset
\let\<\langle
\let\>\rangle
\def\chr{\mathop{\rm char}\nolimits}
\def\End{\mathop{\rm End}\nolimits}
\def\Ext{\mathop{\rm Ext}\nolimits}
\def\gr{\mathop{\rm gr}\nolimits}
\def\Hom{\mathop{\rm Hom}\nolimits}
\def\id{{\rm id}}
\def\Id{\mathop{\rm Id}\nolimits}
\def\Im{\mathop{\rm Im}}
\def\Ker{\mathop{\rm Ker}}
\def\op{^{\rm op}}
\def\St{\mathop{\rm St}}
\def\triv{_{\rm triv}}
\def\mapr#1{{}\mathrel{\smash{\mathop{\longrightarrow}\limits^{#1}}}{}}
\def\lmapr#1#2{{}\mathrel{\smash{\mathop{\count0=#1
  \loop
    \ifnum\count0>0
    \advance\count0 by-1\smash{\mathord-}\mkern-4mu
  \repeat
  \mathord\rightarrow}\limits^{#2}}}{}}
\def\lmapd#1#2#3{\llap{$\vcenter{\hbox{$\scriptstyle{#2}$}}$}
  \left\downarrow\vcenter to#1{}\right.
  \rlap{$\vcenter{\hbox{$\scriptstyle{#3}$}}$}}
\def\diagram#1{\vbox{\halign{&\hfil$##$\hfil\cr #1}}}
\let\al\alpha
\let\be\beta
\let\la\lambda
\let\ph\varphi
\let\si\sigma
\let\ze\zeta
\let\De\Delta
\let\La\Lambda
\let\Si\Sigma
\let\Up\Upsilon
\def\frB{{\frak B}}
\def\frS{{\frak S}}
\def\calD{{\cal D}}
\def\calH{{\cal H}}
\def\calM{{\cal M}}
\def\calR{{\cal R}}
\def\calT{{\cal T}}
\def\bbS{{\Bbb S}}
\def\bbT{{\Bbb T}}
\def\bbZ{{\Bbb Z}}
\def\alt{_{\rm alt}}
\def\ARM{\vphantom{\calM}^{A(R)}\calM}
\def\ARpM{\vphantom{\calM}^{A(R')}\calM}
\def\chit{\widetilde\chi}

\def\grFM{\gr^F\!\calM}
\def\Homk{\Hom_{\,\bbk}}
\def\Hsptri{{\calH'}^{\mskip-2mu\tri}_{\mskip-5mu i}}
\def\laDmu{\hbox{}_\la^{\vphantom{e}}\calD_\mu}
\def\MAR{{\cal M}^{A(R)}}
\def\MARp{{\cal M}^{A(R')}}
\def\Md{{M'}^*}
\def\Mt{\widetilde M}
\def\Mtt{\widetilde M\vphantom{M}'}
\def\muD{\hbox{}_\mu\calD}
\def\muDla{\hbox{}_\mu\calD_\la}
\def\muDlass{\hbox{}_\mu\mskip-1mu\calD_{\mskip-1mu\la}}
\def\muspDla{\muspD_\la^e}
\def\muspD{\hbox{}_{\mu^1}\calD}
\def\oT{\mathaccent"017 T}
\def\Rd{{R'}^*}
\def\Rt{\widetilde R}
\def\Sd{\bbS^{\mskip1mu!}}
\def\Sit{\widetilde\Si}
\def\Tt{\widetilde T}
\def\Ut{\widetilde U}
\def\Vd{{V'}^*}
\def\vxy{v_{x,\mskip1mu y}}

\citedef
Back81
Bich03
Bj-B
Dip-J86
Du92
Geck-P
Gur90
Hai97
Hai02
Hai05
Hay92
Hie93
Hla87
Lar-T91
Lyu87
Man
Mro14
Pol-P
Pr70
Resh-TF89
Sch96
\endcitedef

\frontmatter

\title{On the graded algebras associated with Hecke symmetries}
\author{Serge Skryabin}
\address{Institute of Mathematics and Mechanics,
Kazan Federal University,\break
Kremlevskaya St.~18, 420008 Kazan, Russia}
\email{Serge.Skryabin@kpfu.ru}

\endfrontmatter

\section
Introduction

Let $V$ be a finite dimensional vector space over a field $\bbk$. 
A \emph{Hecke symmetry}\/ with parameter $0\ne q\in\bbk$ is any linear 
operator $R:V^{\ot2}\to V^{\ot2}$ which satisfies the Hecke relation 
$(R+\id)(R-q\cdot\id)=0$ and the braid relation $R_1R_2R_1=R_2R_1R_2$ where 
$R_1=R\ot\id_V$ and $R_2=\id_V\ot R$ are linear operators on $V^{\ot3}$; we 
will be saying that $R$ is a Hecke symmetry \emph{on $V$}.
The symmetries with parameter $q=1$ were considered by Lyubashenko 
\cite{Lyu87}. Many notions and results originated in his work were later 
generalized to the case $q\ne1$.

The $R$-symmetric algebra $\bbS(V,R)$ and the $R$-skewsymmetric algebra $\La(V,R)$ 
are two factor algebras of the tensor algebra $\bbT(V)$. They are regarded as 
analogs of the symmetric and the exterior algebras of $V$. Since the braid 
equation is just a slightly different form of the quantum Yang-Baxter 
equation, there is also a bialgebra $A(R)$ given by the 
Faddeev-Reshetikhin-Takhtajan construction \cite{Resh-TF89}. This bialgebra 
coacts on $V$ universally with respect to the property that the induced 
coaction on $V^{\ot2}$ commutes with $R$ \cite{Lar-T91}. In particular, 
$\bbS(V,R)$ and $\La(V,R)$ are $A(R)$-comodule algebras. The graded algebras 
$\bbS(V,R)$, $\La(V,R)$, $A(R)$ are \emph{quadratic}\/ in the sense that they 
are generated by homogeneous elements of degree 1 and their defining relations 
are of degree 2.

Gurevich's work on Hecke symmetries \cite{Gur90} was motivated by the 
construction of quantum groups not necessarily arising as deformations 
of the classical objects. Even disregarding the Hopf algebraic aspect, 
Hecke symmetries provide a large class of graded algebras with predictable 
properties meaningful from the viewpoint of noncommutative algebraic geometry. 
However, general results on these algebras have been known under the assumption 
that the $q$-integers
$$
[n]_q=1+q+\ldots+q^{n-1}\in\bbk
$$
are nonzero for all integers $n>0$ (this means that $q$ is not a root of 1 
with the exception that $q=1$ is allowed when $\chr\bbk=0$). The assumption 
$\chr\bbk=0$ was also used, but it is less relevant. The present paper makes 
an attempt to investigate several questions without the aforementioned 
restriction on $q$. Particularly, we are interested in Koszulness and 
Gorensteinness of those graded algebras.

A Hecke symmetry $R$ on $V$ gives rise, for each $n>0$, to a representation of 
the Hecke algebra $\calH_n$ of type $A_{n-1}$ in the vector space 
$\bbT_n(V)=V^{\ot n}$. If $[n]_q\ne0$ for all $n$, then these Hecke algebras 
are semisimple, and we will refer to this case as the \emph{semisimple case}. 
Semisimplicity was the main driving force in the earlier results on the graded 
algebras associated with Hecke symmetries.

For $q$ a root of 1 we cannot be too optimistic, as an example at the end of 
section 3 shows. This example suggests that the properties of the graded 
algebras depend on the kind of the Hecke algebra representations in the tensor 
powers of $V$. We will say that an indecomposable $\calH_n$-module has a 
\emph{$1$-dimensional} (respectively, \emph{trivial}) \emph{source} if it is a 
direct summand of an $\calH_n$-module induced from a $1$-dimensional 
(respectively, the trivial $1$-dimensional) representation of a parabolic 
subalgebra. This terminology is explained by the Hecke algebra version of the 
Green correspondence in the modular representation theory of finite groups 
(see Du \cite{Du92}). The following two conditions on a Hecke symmetry $R$ 
will emerge in the statements:

\proclaim
The 1-dimensional source condition.
For each $n>0$ all indecomposable direct summands of $\,V^{\ot n}$ regarded as 
an $\calH_n$-module with respect to the representation arising from $R$ have 
$1$-dimensional sources.
\endproclaim

\proclaim
The trivial source condition.
For each $n>0$ all indecomposable direct summands of $\,V^{\ot n}$ regarded as 
an $\calH_n$-module with respect to the representation arising from $R$ have 
trivial sources.
\endproclaim

In the semisimple case the trivial source condition is obviously satisfied 
since indecomposable modules are simple, and each simple $\calH_n$-module is a 
direct summand of the cyclic free module. We have to consider the weaker 
1-dimensional source condition in order to include the supersymmetry on a 
$\bbZ/2\bbZ$-graded vector space (in that example $q=1$, so that $\calH_n$ is 
just the group algebra of the symmetric group $\frS_n$, and $V^{\ot n}$ is a 
direct sum of $\calH_n$-modules induced from not necessarily trivial 
$1$-dimensional representations of parabolic subalgebras, but $\calH_n$ is not 
semisimple when $n\ge\chr\bbk>0$). This condition is also satisfied for the 
Hecke symmetries related to the standard quantum supergroups. There is one 
Hecke symmetry on a 2-dimensional vector space for which $q$ is a 4th root of 
1 and the 1-dimensional source condition fails to hold (see section 3). 
However, this Hecke symmetry is not closed. This raises the following

\proclaim
Question.
Does the $1$-dimensional source condition hold for every closed Hecke symmetry\/{\rm?}
\endproclaim

A closed symmetry extends to a braiding on a monoidal subcategory of the 
category of finite dimensional vector spaces containing $V$ and its dual 
objects (see \cite{Gur90} for the precise definition). The results we are 
going to present do not depend on the closedness of $R$.

\proclaim
Theorem 3.1.
Suppose that $R$ satisfies the $1$-dimensional source condition. Then the 
$R$-symmetric algebra $\bbS(V,R)$ and the $R$-skewsymmetric algebra 
$\La(V,R)$ are Koszul. Their Hilbert series satisfy the relation 
\ $h_{\bbS(V,R)}(t)\,h_{\La(V,R)}(-t)=1$.
\endproclaim

In the semisimple case exactness of certain complexes proved by Gurevich 
amounts to the conclusion of the theorem stated above, although 
Koszulness of graded algebras was not mentioned in \cite{Gur90} explicitly. By 
a systematic use of various projectors all considerations in \cite{Gur90} were 
done in terms of subspaces rather than factor spaces of the tensor powers of 
the initial vector space. The realization of Koszul complexes based on 
projectors is not appropriate for arbitrary $q$, however. Koszulness in the 
semisimple case was also considered by \PHHai \cite{Hai97}. We will discuss 
more general results in a moment.

\proclaim
Theorem 4.5.
Suppose that $R$ satisfies the trivial source condition. Suppose also that\/ 
$\dim\La_n(V,R)=1$ and $\La_{n+1}(V,R)=0$ for some $n>0$. Then\/ $\La(V,R)$ is 
a Frobenius algebra{\rm,} while\/ $\bbS(V,R)$ is a Gorenstein algebra of 
global dimension $n$.
\endproclaim

The subscripts here and elsewhere in the paper indicate the homogeneous 
components of graded spaces. The first conclusion in Theorem 4.5 is 
equivalent to nondegeneracy of the bilinear pairings 
$\La_i(V,R)\times\La_{n-i}(V,R)\to\La_n(V,R)$ arising from the multiplication. 
This was proved by Gurevich in the semisimple case. 

By a \emph{Gorenstein algebra}\/ we mean any positively graded algebra 
$A=\bigoplus_{i=0}^\infty A_i$ with $A_0=\bbk$ for which 
$\Ext_A^\bullet(\bbk,A)$ is 1-dimensional. No Noetherian conditions are 
requested. The two conclusions in Theorem 4.5 are closely related (see 
\cite{Pol-P}), and the essential job will be to prove the first one.

Let now $R'$ be a second Hecke symmetry on another finite dimensional vector 
space $V'$. Suppose that $R$ and $R'$ have the same parameter $q$ of the Hecke 
relation. Then there is a graded algebra $A(R',R)$ whose construction 
generalizes that of $A(R)$. In fact $A(R',R)=A(R)$ when $R'=R$. The algebra 
$A(R',R)$ in a different notation was introduced by \PHHai \cite{Hai02} under 
the name ``quantum hom-space". We will consider yet another graded algebra 
$E(R',R)$ whose relationship with $A(R',R)$ is similar to that between 
$\La(V,R)$ and $\bbS(V,R)$. In a different notation this algebra was also 
introduced in \cite{Hai02} under the name ``exterior algebra on the quantum 
hom-space".

\proclaim
Theorem 6.2.
Suppose that both $R$ and $R'$ satisfy the $1$-dimensional source condition. 
Then the graded algebras $A(R',R)$ and $E(R',R)$ are Koszul. Their Hilbert 
series satisfy the relation \ $h_{A(R',R)}(t)\,h_{E(R',R)}(-t)=1$.
\endproclaim

In the semisimple case this was proved in \cite{Hai02}. The argument used by 
\PHHai is based on the observation that Koszulness of $A(R',R)$ is equivalent 
to a certain property of $\calH_n\op\ot\calH_n$-modules concerned with 
distributivity of collections of subspaces in these modules. If $R_{k,q}$ is 
the Hecke symmetry corresponding to the standard quantum $GL_k$ with parameter 
$q$, then $A(R_{k,q})$ is known to possess a PBW basis, and its Koszulness 
follows from Priddy's theorem \cite{Pr70}. If $[n]_q\ne0$ for all $n$, the 
semisimplicity of $\calH_n\op\ot\calH_n$ ensures then the desired property for 
all simple modules, and therefore for arbitrary $\calH_n\op\ot\calH_n$-modules. 
This in turn yields Koszulness of $A(R',R)$ for arbitrary Hecke symmetries 
with the same parameter $q$.

In the present paper we prove directly exactness of certain complexes defined 
with respect to representations of $\calH_n\op\ot\calH_n$, and thus we derive 
Theorem 6.2 solely from Hecke algebra considerations, avoiding the use of 
Priddy's theorem.

\proclaim
Theorem 6.6.
Suppose that both $R$ and $R'$ satisfy the trivial source condition. 
If $\dim E_n(R',R)=1$ and $E_{n+1}(R',R)=0$ for some integer 
$n>0,$ then $E(R',R)$ is a Frobenius algebra{\rm,} while $A(R',R)$ 
is a Gorenstein algebra of global dimension $n$.
\endproclaim

If $\dim V'=1$ and $R'$ is the multiplication by $q$, then $A(R',R)=\bbS(V,R)$ 
and $E(R',R)=\La(V,R)$. Thus the first two results discussed in this 
introduction are a special case of the two subsequent ones. We nevertheless 
provide separate proofs in this special case. It serves as a model for the 
more complicated Theorems 6.2 and 6.6 where we need some lengthy verifications 
done in section 5 of the paper.

The quantum hom-spaces were considered in \cite{Hai02} in connection with the 
quantum version of the classical invariant theory. There is an even more 
obvious role played by the algebra $A(R',R)$. Being equipped with an 
$A(R')\,$-$A(R)$ bicomodule structure, it gives rise to functors between the 
corepresentation categories of the two bialgebras $A(R)$ and $A(R')$. 
For each coalgebra $C$ let $^C\!\calM$ and $\calM^C$ stand for the categories 
of left and right comodules. The sign $\,\square_C\,$ denotes the cotensor 
product of comodules.

\proclaim
Theorem 7.2.
Suppose that for each $n>1$ the indecomposable $\calH_n$-modules isomorphic to 
direct summands of\/ $\bbT_n(V')$ are the same as those isomorphic to direct 
summands of\/ $\bbT_n(V)$. Then the functors
$$
A(R',R)\square_{A(R)}\hbox{}
?\qquad{\rm and}\qquad?\,\hbox{}\square_{A(R)}A(R,R')
$$
are braided monoidal equivalences $\,\ARM\mapr{}\ARpM\,$ and $\,\MAR\mapr{}\MARp$.
\endproclaim

This should be compared with the monoidal equivalences between the 
corepresentation categories of two Hopf algebras. According to a result of 
Schauenburg \cite{Sch96} such equivalences are given by the cotensor product 
functors determined by the so-called bi-Galois algebras. In this way Bichon 
\cite{Bich03} and Mrozinski \cite{Mro14} showed that the categories of 
comodules over certain Hopf algebras associated with bilinear forms are 
monoidally equivalent to the respective categories defined for the standard 
quantizations of $SL_2$ and $GL_2$. Those Hopf algebras correspond to a 
special class of Hecke symmetries. \PHHai dealt with the Hopf envelopes of the 
bialgebras $A(R)$ for arbitrary Hecke symmetries under the previously 
mentioned restriction on $q$ and the characteristic of $\bbk$. By the main 
result of \cite{Hai05} the category of comodules over such a Hopf algebra is 
determined, up to monoidal equivalence, by the parameter $q$ and the birank of 
$R$.

Theorem 7.2 is a similar, to some extent, result for bialgebras, though the 
bicomodule algebra $A(R,R')$ is definitely not bi-Galois (Galois algebras 
exist only for Hopf algebras). In fact, everything what is needed for the 
equivalence here is present already in the construction of the algebras 
involved. Note that there are no restrictions on $R$ in Theorem 7.2.

\section
1. Hecke algebra preliminaries

We denote by $\frS_n$ the symmetric group of permutations of the set 
$\{1,\ldots,n\}$. Let $\frB_n=\{\tau_1,\ldots,\tau_{n-1}\}$ be the set of 
basic transpositions $\tau_i=(i,i+1)$. The \emph{length} $\ell(\si)$ of a 
permutation $\si\in\frS_n$ is the smallest number of factors in the 
expressions of $\si$ as product of basic transpositions. By the letter $e$ we 
denote the identity permutation and also the trivial subgroup of $\frS_n$.

The subgroups of $\frS_n$ generated by subsets of $\frB_n$ are called 
\emph{Young subgroups} and are parametrized traditionally by 
\emph{compositions} of $n$, i.e. by finite sequences of positive integers 
summing up to $n$. The Young subgroup $\frS_\la$ labelled by a composition 
$\la=(\la_1,\ldots,\la_k)$ is generated by the set
$$
\frB_\la=\{\tau_j\in\frB_n\mid j\ne\la_1+\ldots+\la_i\hbox{ for each }i\le k\}
$$
and is isomorphic to $\frS_{\la_1}\times\ldots\times\frS_{\la_k}$. In 
particular, the subgroup $\frS_{i,\,n-i}$ corresponding to the composition 
$(i,n-i)$ is generated by $\{\tau_j\in\frB_n\mid j\ne i\}$. For each $\la$ the 
pair $(\frS_\la,\frB_\la)$ is a Coxeter system. We will use standard facts 
concerning Coxeter groups and the respective Hecke algebras. For reference 
see, e.g., \cite{Bj-B}, \cite{Geck-P}.

Given a pair of Young subgroups $\frS_\la\sbs\frS_\nu$, each coset of 
$\frS_\la$ in $\frS_\nu$ contains a unique element of minimal length called the 
\emph{distinguished coset representative}. We denote by 
$\calD(\frS_\nu/\frS_\la)$ and $\calD(\frS_\la\backslash\frS_\nu)$ the sets of 
distinguished representatives of the respective cosets. Recall that
$$
\eqalign{
\calD(\frS_\nu/\frS_\la)&{}=\{\pi\in\frS_\nu\mid\pi\tau_i>\pi\hbox{ for all 
$\tau_i\in\frB_\la$}\}\cr
&{}=\{\pi\in\frS_\nu\mid\ell(\pi\si)=\ell(\pi)+\ell(\si)
\hbox{ for all $\si\in\frS_\la$}\}\cr
}
$$
where comparison $\pi\tau_i>\pi$ refers to the Bruhat order (recall that 
$|\ell(\pi\tau_i)-\ell(\pi)|=1$, and $\pi\tau_i>\pi$ if and only if 
$\ell(\pi\tau_i)>\ell(\pi)$). The set 
$\,\calD(\frS_\la\backslash\frS_\nu)=\calD(\frS_\nu/\frS_\la)^{-1}\,$ admits 
similar characterizations. For another Young subgroup $\frS_\mu$ of $\frS_\nu$ 
the set of \emph{distinguished $\frS_\mu$-$\,\frS_\la$ double coset 
representatives} 
is
$$
\calD(\frS_\mu\backslash\frS_\nu/\frS_\la)=
\calD(\frS_\nu/\frS_\la)\cap\calD(\frS_\mu\backslash\frS_\nu).
$$
We will also use shorter notation: $\calD_\la=\calD(\frS_n/\frS_\la)$, 
$\muD=\calD(\frS_\mu\backslash\frS_n)$, and
$\muDla=\calD(\frS_\mu\backslash\frS_n/\frS_\la)$.

Let $\bbk$ be the ground field. The Hecke algebra $\calH_n=\calH_n(q)$ 
of type $A_{n-1}$ with parameter $q\in\bbk$ is presented by generators 
$T_1,\ldots,T_{n-1}$ and relations
$$
\openup1\jot
\displaylines{
T_iT_jT_i=T_jT_iT_j\hbox{\ \ when $|i-j|=1$},\qquad
T_iT_j=T_jT_i\hbox{\ \ when $|i-j|>1$},\cr
(T_i-q)(T_i+1)=0\quad\hbox{for $i=1,\ldots,n-1$}.\cr
}
$$
It has a standard basis $\{T_\si\mid\si\in\frS_n\}$ characterized by the 
properties that $T_e=1$ and $T_{\tau_i\si}=T_iT_\si$ whenever $\tau_i\si>\si$. 
Recall that $T_\pi T_\si=T_{\pi\si}$ for each pair $\pi,\si\in\frS_n$ such that 
$\ell(\pi\si)=\ell(\pi)+\ell(\si)$.

The definition of $\calH_n(q)$ makes sense also when $q=0$. We will use this 
algebra called the \emph{$0$-Hecke algebra} on several occasions.

The elements $\{T_\si\mid\si\in\frS_\la\}$ form a basis for 
the \emph{parabolic subalgebra} $\calH_\la$ generated by 
$\{T_i\mid\tau_i\in\frB_\la\}$. For each $k<n$ we identify the symmetric 
group $\frS_k$ with the subgroup of $\frS_n$ generated by 
$\{\tau_i\mid0<i<k\}$ and the Hecke algebra $\calH_k$ with the subalgebra of 
$\calH_n$ generated by $\{T_i\mid0<i<k\}$. By a convention $\frS_0=\frS_1=e$ 
and $\calH_0=\calH_1=\bbk$.

For each 1-dimensional representation of a parabolic subalgebra $\calH_\la$ of 
$\calH_n$ given by an algebra homomorphism $\chi:\calH_\la\to\bbk$ we denote by 
$\bbk(\chi)$ the corresponding 1-dimensional $\calH_\la$-module. The induced 
$\calH_n$-module
$$
M=\calH_n\ot_{\calH_\la}\bbk(\chi)
$$
has a basis $\{T_\si c\mid\si\in\calD_\la\}$ where $c=1\ot1$ is the 
\emph{canonical generator} of $M$. Here and later by a \emph{basis}\/ we mean 
a basis over the ground field $\bbk$.

By Deodhar's lemma (see \cite{Geck-P, 2.1.2}) for each $\tau_i\in\frB_n$ and 
$\si\in\calD_\la$ we have either $\tau_i\si\in\calD_\la$ or 
$\tau_i\in\si\frB_\la\si^{-1}$. Hence
$$
\eqalign{
\calD_\la&{}=A_i(\la)\cup\tau_i A_i(\la)\cup B_i(\la),\quad
\hbox{a disjoint union, where}\cr
\noalign{\smallskip}
A_i(\la)&{}=\{\si\in\calD_\la\mid\tau_i\si\in\calD_\la
{\rm\ and\ }\tau_i\si>\si\}\,,\cr
\tau_iA_i(\la)&{}=\{\si\in\calD_\la\mid\tau_i\si\in\calD_\la
{\rm\ and\ }\tau_i\si<\si\}\,,\cr
B_i(\la)&{}=\{\si\in\calD_\la\mid\tau_i\si\notin\calD_\la\}
=\{\si\in\calD_\la\mid\si^{-1}\tau_i\mskip1mu\si\in\frB_\la\}.
}
$$
If $\si\in B_i(\la)$, then $\tau_i\si=\si\tau_j>\si$ for some 
$\tau_j\in\frB_\la$, and it follows that
$T_iT_\si=T_{\tau_i\si}=T_\si T_j$. In this case $T_\si c$ is an eigenvector 
for the linear operator $(T_i)_M$ by which $T_i$ acts on $M$. Let us denote by 
$\chi_\si(T_i)$ the corresponding eigenvalue.

Note that $\si^{-1}(i)<\si^{-1}(i+1)$ when $\tau_i\si>\si$. Since 
$\si^{-1}\tau_i\si$ is the transposition of $\si^{-1}(i)$ and $\si^{-1}(i+1)$, 
the equality $\tau_i\si=\si\tau_j$ implies then that $j=\si^{-1}(i)$. 
Thus the action of $T_i$ on the basis elements of $M$ is given by the formulas
$$
T_iT_\si c=\cases{
T_{\tau_i\si}c & if $\si\in A_i(\la)$, \cr
\noalign{\smallskip}
(q-1)\,T_\si c+q\,T_{\tau_i\si}c & if $\si\in\tau_iA_i(\la)$, \cr
\noalign{\smallskip}
\chi_\si(T_i)\,T_\si c\quad
\hbox{with $\chi_\si(T_i)=\chi(T_{\si^{-1}\tau_i\si})=\chi(T_{\si^{-1}(i)})$} 
& if $\si\in B_i(\la)$ \cr
}
$$

The restriction of the induced $\calH_n$-module $M$ to a second parabolic 
subalgebra $\calH_\mu$ of $\calH_n$ is given by the Mackey formula 
\cite{Dip-J86, Th. 2.7}:
$$
M=\bigoplus_{\pi\in\,\muDlass}M(\pi),\qquad M(\pi)\cong
\calH_\mu\ot_{\calH_{\nu(\pi)}}\bbk(\chi_\pi),
$$
where $\nu(\pi)$ is the composition of $n$ such that 
$\frS_{\nu(\pi)}=\frS_\mu\cap\pi\frS_\la\pi^{-1}$ and 
$\chi_\pi$ is the 1-dimensional representation of the corresponding parabolic 
subalgebra $\calH_{\nu(\pi)}$ with the values $\chi_\pi(T_i)$ on the 
generators $T_i$ of $\calH_{\nu(\pi)}$ defined in the preceding paragraph. 
Here $M(\pi)$ is the $\calH_\mu$-submodule of $M$ generated by $T_\pi c$. Its 
basis over $\bbk$ is formed by the elements $T_\si c$ with $\si\in\calD_\la 
\cap\frS_\mu\pi\frS_\la$.

The assignments $T_i\mapsto q$ (respectively, $T_i\mapsto -1$) for each $i$ 
such that $\tau_i\in\frB_\la$ define the \emph{trivial} (respectively, 
\emph{alternating}) representation $\calH_\la\to\bbk$. We denote by $\bbk\triv$ 
and $\bbk\alt$ the corresponding $\calH_\la$-modules. They make sense for each 
parabolic subalgebra of $\calH_n$, in particular, for $\calH_n$ itself. If 
$q\ne0$, then every 1-dimensional $\calH_n$-module is isomorphic to either 
$\bbk\triv$ or $\bbk\alt$. If $q=-1$, then $\bbk\triv=\bbk\alt$, and in fact 
all 1-dimensional representations of any parabolic subalgebra coincide.

If $q\ne0$ then the induction functor from any parabolic subalgebra is 
isomorphic to the coinduction functor. In other words,
$$
\Hom_{\calH_n}(N,\,\calH_n\ot_{\calH_\la}U)\cong\Hom_{\calH_\la}(N,U)
$$
for any left $\calH_n$-module $N$ and left $\calH_\la$-module $U$. This is 
a general property of the Hecke algebras of Coxeter groups which we view as 
part of the Frobenius reciprocity (see \cite{Dip-J86, Th. 2.5, 2.6} and 
\cite{Geck-P, 9.1.7}). In particular, $\calH_n$ can be replaced in the 
isomorphism above by any parabolic subalgebra larger than the given $\calH_\la$.

\proclaim
Lemma 1.1.
Let $M=\calH_n\ot_{\calH_\la}\bbk(\chi)$ and $N=\calH_n\ot_{\calH_\mu}\bbk(\ze)$ 
be $\calH_n$-modules induced from $1$-dimensional representations of two 
parabolic subalgebras of $\calH_n$. If $q\ne0$ then
$$ 
\dim\Hom_{\calH_n}(N,M)=f(\ze,\chi)
$$
where $f(\ze,\chi)=\#\{\pi\in\muDla\mid\ze(T_i)=\chi_\pi(T_i)
{\rm\ for\ all\ }i{\rm\ with\ }\tau_i\in\frS_\mu\cap\pi\frS_\la\pi^{-1}\}$.
\endproclaim

\Proof.
Consider the Mackey decomposition $M=\bigoplus M(\pi)$ with 
respect to $\calH_\mu$. Then
$$
\Hom_{\calH_n}(N,M)\cong\Hom_{\calH_\mu}\bigl(\bbk(\ze),M\bigr)
\cong\bigoplus_{\pi\in\,\muDlass}\Hom_{\calH_\mu}\bigl(\bbk(\ze),M(\pi)\bigr).
$$
Recall that for each $\pi$ the $\calH_\mu$-module $M(\pi)$ is induced from the 
$1$-dimensional representation $\chi_\pi$ of the parabolic subalgebra 
$\calH_{\nu(\pi)}$ corresponding to the Young subgroup 
$\frS_\mu\cap\pi\frS_\la\pi^{-1}$. By the Frobenius reciprocity
$$
\Hom_{\calH_\mu}\bigl(\bbk(\ze),M(\pi)\bigr)\cong
\Hom_{\calH_{\nu(\pi)}}\bigl(\bbk(\ze),\bbk(\chi_\pi)\bigr),
$$
and this space is nonzero if and only if $\ze$ agrees with $\chi_\pi$ on 
$\calH_{\nu(\pi)}$.
\endproof

The next lemma will be crucial for establishing the relation between the 
Hilbert series of the pairs of graded algebras in Theorems 3.1 and 6.2.

\proclaim
Lemma 1.2.
Suppose that $q\ne0$. 
Let $M$ and $N$ be finite dimensional $\calH_n$-modules whose indecomposable 
direct summands all have $1$-dimensional sources. Then
$$ 
\dim\Hom_{\calH_n}(N,M)=\dim\Hom_{\calH_n}(M,N).
$$
\endproclaim

\Proof.
Since $\Hom_{\calH_n}$ is an additive functor of both arguments, it suffices to 
check the desired equality when $M$ and $N$ are indecomposable modules, so that 
they are direct summands of $\calH_n$-modules induced from $1$-dimensional 
representations of parabolic subalgebras. If $\calH_n$ is semisimple, then the 
indecomposable modules are simple, and the conclusion is obviously true since 
for two simple modules $\Hom_{\calH_n}(N,M)=0$ unless $N\cong M$.

If $\calH_n$ is not semisimple, we will apply the specialization argument. 
Let $M$ and $N$ be as in Lemma 1.1. In this case Lemma 1.1 gives the exact 
values for the dimensions of $\Hom_{\calH_n}(N,M)$ and $\Hom_{\calH_n}(M,N)$, 
but we have to show that
$$ 
\dim\Hom_{\calH_n}(N',M')=\dim\Hom_{\calH_n}(M',N')
$$
whenever $N'$ is a direct summand of $N$ and $M'$ is a direct summand of $M$.

Let $O$ be the completion of the polynomial ring $\bbk[t]$ in one 
indeterminate $t$ at its maximal ideal generated by $t-q$. Then $O$ is a 
complete discrete valuation ring with residue field isomorphic to $\bbk$. The 
specialization homomorphism $O\to\bbk$ sends $t$ to $q$. Denote by $K$ the 
field of fractions of $O$. Let $\calH_n(t)$ be the Hecke algebra of type 
$A_{n-1}$ with parameter $t$ over the ring $O$. Then 
$\calH_n(t)\ot_O\bbk\cong\calH_n$, while
$$
\calH_n(t)_K=\calH_n(t)\ot_OK
$$
is a semisimple Hecke algebra of type $A_{n-1}$ over the field $K$.

Let $\calH_\la(t)$ and $\calH_\mu(t)$ be the parabolic subalgebras of 
$\calH_n(t)$ corresponding to the two compositions $\la$ and $\mu$ of $n$. 
Define ring homomorphisms $\chi^t:\calH_\la(t)\to O$ and 
$\ze^t:\calH_\mu(t)\to O$ by the formulas
$$
\chi^t(T_i)=\cases{t & if $\chi(T_i)=q$,\cr \noalign{\smallskip}
-1 & otherwise,}\qquad\quad
\ze^t(T_i)=\cases{t & if $\ze(T_i)=q$,\cr \noalign{\smallskip}
-1 & otherwise.}
$$
Let $O(\chi^t)=O$ with the $\calH_\la(t)$-module structure given by $\chi^t$, 
and $O(\ze^t)=O$ the similar $\calH_\mu(t)$-module with respect to $\ze^t$. 
Put
$$
M(t)=\calH_n(t)\ot_{\calH_\la(t)}O(\chi^t),\qquad\quad
N(t)=\calH_n(t)\ot_{\calH_\mu(t)}O(\ze^t).
$$
These are $O$-free $\calH_n(t)$-modules such that $M(t)\ot_O\bbk\cong M$ and 
$N(t)\ot_O\bbk\cong N$. The $\calH_n(t)_K$-modules $M(t)\ot_OK$ and $N(t)\ot_OK$ 
are induced from $1$-dimensional representations of parabolic subalgebras of 
$\calH_n(t)_K$. So Lemma 1.1 yields
$$
\dim_K\Hom_{\calH_n(t)_K}\bigl(N(t)\ot_OK,\,M(t)\ot_OK\bigr)=f(\ze^t,\chi^t).
$$
By exactness of the functor $?\ot_OK$, we have
$$
\Hom_{\calH_n(t)_K}\bigl(N(t)\ot_OK,\,M(t)\ot_OK\bigr)\cong
\Hom_{\calH_n(t)}\bigl(N(t),\,M(t)\bigr)\ot_OK.
$$
Since the $O$-module $\Hom_{\calH_n(t)}\bigl(N(t),M(t)\bigr)$ is torsionfree, 
it has to be free of rank equal to $f(\ze^t,\chi^t)$. For each pair of indices 
$i,j$ such that $\tau_i\in\frB_\mu$ and $\tau_j\in\frB_\la$ it is seen from 
the definition of $\chi^t,\,\ze^t$ that $\ze^t(T_i)=\chi^t(T_j)$ if and only 
if $\ze(T_i)=\chi(T_j)$ since $\ze(T_i)$ and $\chi(T_j)$ can be equal to only 
$q$ or $-1$. Therefore $\ze^t(T_i)=\chi^t_\pi(T_i)$ for some $\pi\in\muDla$ 
and $i$ such that $\tau_i\in\frS_\mu\cap\pi\frS_\la\pi^{-1}$ if and only if 
$\ze(T_i)=\chi_\pi(T_i)$. It follows that $f(\ze^t,\chi^t)=f(\ze,\chi)$, as 
defined in Lemma 1.1. 

A homomorphism $\ph:N(t)\to M(t)$ induces zero map $N\to M$ if and only if 
$\,\Im\ph\sbs(t-q)M(t)$, so that 
$(t-q)^{-1}\ph\in\Hom_{\calH_n(t)}\bigl(N(t),M(t)\bigr)$ for such a $\ph$. 
This shows that $\,\Hom_{\calH_n(t)}\bigl(N(t),M(t)\bigr)\ot_O\bbk\,$ embeds in 
$\,\Hom_{\calH_n}(N,M)$, but then
$$
\Hom_{\calH_n(t)}\bigl(N(t),M(t)\bigr)\ot_O\bbk\cong\Hom_{\calH_n}(N,M)
$$
by comparison of dimensions.

As a special case we get $\bigl(\End_{\calH_n(t)}M(t)\bigr)\ot_O\bbk\cong
\End_{\calH_n}M$. As is well-known, in this situation each idempotent of 
$\End_{\calH_n}M$ can be lifted to an idempotent of $\End_{\calH_n(t)}M(t)$. 
Direct summands of an arbitrary module are determined by idempotents in its 
endomorphism ring. Thus we can find a direct summand $M'(t)$ of the 
$\calH_n(t)$-module $M(t)$ such that $M'(t)\ot_O\bbk\cong M'$. Similarly, 
there is a direct summand $N'(t)$ of $N(t)$ satisfying $N'(t)\ot_O\bbk\cong 
N'$. Being submodules of $O$-free modules, both $M'(t)$ and $N'(t)$ are 
themselves $O$-free.

Since every $\calH_n$-module homomorphism $N\to M$ lifts to an 
$\calH_n(t)$-module homomorphism $N(t)\to M(t)$, it follows that every 
$\calH_n$-module homomorphism $N'\to M'$ lifts to an $\calH_n(t)$-module 
homomorphism $N'(t)\to M'(t)$. This entails
$$
\Hom_{\calH_n(t)}\bigl(N'(t),M'(t)\bigr)\ot_O\bbk\cong\Hom_{\calH_n}(N',M').
$$
Since $\Hom_{\calH_n(t)}\bigl(N'(t),M'(t)\bigr)$ is a free $O$-module, we 
deduce that
$$
\eqalign{
\dim_{\mskip1mu\bbk}\Hom_{\calH_n}(N',M')
&{}=\dim_K\Hom_{\calH_n(t)}\bigl(N'(t),M'(t)\bigr)\ot_OK\cr
&{}=\dim_K\Hom_{\calH_n(t)_K}\bigl(N'(t)\ot_OK,\,M'(t)\ot_OK\bigr).\cr
}
$$
By symmetry
$$
\dim_{\mskip1mu\bbk}\Hom_{\calH_n}(M',N')
=\dim_K\Hom_{\calH_n(t)_K}\bigl(M'(t)\ot_OK,\,N'(t)\ot_OK\bigr),
$$
and the equality $\dim\Hom_{\calH_n}(N',M')=\dim\Hom_{\calH_n}(M',N')$ follows 
from the already discussed semisimple case.
\endproof

A suitable version of Lemma 1.1 is valid also when $q=0$. Later we will need 
only a special case of that fact:

\proclaim
Lemma 1.3.
Let $M=\calH_n\ot_{\calH_\la}\bbk(\chi)$. Then
$$
\dim\Hom_{\calH_n}(\bbk\alt,M)=\cases{
1 & \rm if $\chi$ is the alternating representation of $\calH_\la$,\cr
\noalign{\smallskip}
0 & \rm otherwise.
}
$$
\endproclaim

\Proof.
This conclusion is a consequence of the Frobenius reciprocity when $q\ne0$. 
For $q=0$ it is derived as follows. Since $T_i^2=-T_i$, it is seen from the 
formulas for the action of $T_i$ that $T_iM$ is contained in the linear span 
of the basis elements $T_\si c$ with $\si\in\calD_\la\setm A_i(\la)$. If 
$v\in M$ is such that $T_iv=-v$ for all $i=1,\ldots,n-1$, then 
$v\in\bigcap_{i=1}^{n-1}T_iM$.

Let $w_n$ and $w_\la$ be the longest elements of $\frS_n$ and $\frS_\la$, 
respectively. By \cite{Geck-P, 2.2.1} $d_\la=w_nw_\la$ is the unique element of 
maximal length in $\calD_\la$, and $\calD_\la$ consists precisely of all 
suffixes of $d_\la$ which are elements $\si\in\frS_n$ satisfying 
$\ell(d_\la)=\ell(d_\la\si^{-1})+\ell(\si)$. If $\si\in\calD_\la$ and $\si\ne d_\la$, then 
there exists $\tau_i\in\frB_n$ such that $\tau_i\si>\si$ and $\tau_i\si$ is a 
suffix of $d_\la$, so that $\tau_i\si\in\calD_\la$, i.e. $\si\in A_i(\la)$. This 
shows that $T_\si c$ can be involved with nonzero coefficient in the 
expression for $v$ only when $\si=d_\la$, i.e. $v$ has to be a scalar multiple 
of $T_{d_\la}c$.

Now $T_{d_\la}c$ is an eigenvector for all operators $(T_i)_M$. Moreover, the 
1-dimensional $\calH_n$-submodule of $M$ spanned by $T_{d_\la}c$ is isomorphic 
to $\bbk\alt$ if and only if $\chi_{d_\la}(T_i)$ is equal to $-1$ for each $i$ 
such that $d_\la\in B_i(\la)$. Since the conjugation by $w_n$ (respectively, 
by $w_\la$) map $\frB_n$ (respectively, $\frB_\la$) onto itself, we have 
$d_\la\frB_\la d_\la^{-1}\sbs\frB_n$. This means that 
$$
\{\tau_i\mid d_\la\in B_i(\la)\}=d_\la\frB_\la d_\la^{-1}.
$$
Since $d_\la\tau_j\mskip1mu d_\la^{-1}=\tau_{d_\la(j)}$ and  
$\chi_{d_\la}(T_{d_\la(j)})=\chi(T_j)$ for each $\tau_j\in\frB_\la$, the 
previous condition on $\chi_{d_\la}$ is equivalent to the condition that 
$\chi$ is alternating.
\endproof

\proclaim
Lemma 1.4.
Suppose that $q\ne0$. Let $M=\bbk(\chi)\ot_{\calH_\la}\calH_n$ be the right 
$\calH_n$-module induced from a $1$-dimensional representation $\chi$ of the 
parabolic subalgebra $\calH_\la$. Considering the dual space $M^*$ as a left 
$\calH_n$-module with respect to the natural action of $\calH_n,$ we have 
$\,M^*\cong\calH_n\ot_{\calH_\la}\bbk(\chi)$.
\endproclaim

\Proof.
Since $M^*\cong\Hom_{\calH_\la}\bigl(\calH_n,\bbk(\chi)\bigr)$, the conclusion 
is a consequence of the Frobenius reciprocity.
\endproof

Put $\Tt_i=q-1-T_i$ for each $i=1,\ldots,n-1$. The assignment $T_i\mapsto\Tt_i$ 
extends to an involutive automorphism of $\calH_n$. We denote by $\Mt$ the 
$\calH_n$-module whose underlying vector space coincides with $M$, but the new 
action of $T_i$ is given by the original action of $\Tt_i$.

\proclaim
Lemma 1.5.
If $M=\calH_n\ot_{\calH_\la}\bbk(\chi),$ then $\Mt\cong\calH_n\ot_{\calH_\la}\bbk(\chit)$ 
where $\chit:\calH_\la\to\bbk$ is the $1$-dimensional representation such that 
$$
\chit(T_i)=\chi(\Tt_i)=q-1-\chi(T_i)\quad
\hbox{\rm for each $i$ with $\tau_i\in\frB_\la$}.
$$
\endproclaim

\Proof.
Let $c$ be the canonical generator of $M$. Then $\Tt_ic=\chit(T_i)c$ for each 
$i$ such that $\tau_i\in\frB_\la$. Clearly $c$ generates also $\Mt$. Hence 
there is a surjective homomorphism $\calH_n\ot_{\calH_\la}\bbk(\chit)\to\Mt$ 
which has to be bijective since the two modules here have equal dimensions.  
\endproof

\setitemsize(a)
\proclaim
Lemma 1.6.
Suppose that $q\ne-1$. Let $M$ and $N$ be two left $\calH_n$-modules. Given a 
$\bbk$-linear map $\ph:N\to M,$ the following conditions are equivalent\/{\rm:}

\item(a)
for each $i=1,\ldots,n-1$ there exists a $\bbk$-linear map $\psi_i:N\to M$ such 
that $\ph=\psi_iT_i-T_i\psi_i\,$ {\rm(}i.e. 
$\ph(x)=\psi_i(T_ix)-T_i\psi_i(x)$ for all $x\in N${\rm),}

\item(b)
$\ph$ is a homomorphism of $\calH_n$-modules $N\to\Mt$.

If $\,q=-1,$ then {\rm(a) $\Rar$ (b)}.
\endproclaim

\Proof.
Note that $T_i+\Tt_i=q-1$ and $T_i\Tt_i=-q$. If $\psi_i$ satisfies (a), then 
$$
\ph T_i-\Tt_i\ph=\psi_iT^2_i-(T_i+\Tt_i)\psi_iT_i+T_i\Tt_i\psi_i
=\psi_i\bigl(T^2_i-(q-1)T_i-q\bigr)=0.
$$
Hence (a) implies that $\ph T_i=\Tt_i\ph$\/ for all $i$, i.e. 
$\ph\in\Hom_{\calH_n}(N,\Mt)\strut$. Conversely, if $\ph T_i=\Tt_i\ph$, 
then taking $\psi_i=(q+1)^{-2}(\Tt_i-T_i)\ph$, we get
$$
\psi_iT_i-T_i\psi_i=(q+1)^{-2}(\Tt_i-T_i)^2\ph=\ph
$$
since $\,(\Tt_i-T_i)^2=(\Tt_i+T_i)^2-4T_i\Tt_i=(q-1)^2+4q=(q+1)^2$. 
\endproof

\proclaim
Lemma 1.7.
Suppose that $q=-1$. Let $M=\calH_n\ot_{\calH_\la}\bbk$ and 
$N=\calH_n\ot_{\calH_\mu}\bbk$ be $\calH_n$-modules induced from the 
$1$-dimensional representations of two parabolic subalgebras. Denote by $c$ 
and $c'$ their canonical generators. For a homomorphism $\ph:N\to M$ the 
following conditions are equivalent\/{\rm:}

\item(a)
$\ph$ factors through a free $\calH_n$-module{\rm,}

\item(b)
$\ph(c')\in x_\mu M$ where $x_\mu=\sum_{\si\in\frS_\mu}T_\si,$

\item(c)
$\ph(c')\in(T_i+1)M$ for each $i$ such that $\tau_i\in\frB_\mu\,,$

\item(d)
for each $i=1,\ldots,n-1$ there exists a $\bbk$-linear map $\psi_i:N\to M$ 
such that\hfil\break $\ph=T_i\psi_i+\psi_iT_i+2\psi_i,$ i.e.  
$\ph=\psi_iT_i-\Tt_i\psi_i$ where $\,\Tt_i=-2-T_i$.

\noindent
The space of all homomorphisms $N\to M$ satisfying {\rm(a) -- (d)} has a basis 
indexed by the set\/ 
$\{\pi\in\muDla\mid\frS_\mu\cap\pi\frS_\la\pi^{-1}=e\}$ 
of distinguished representatives of the double cosets with the 
trivial intersection property.
\endproclaim

\Proof.
(a)$\,\Rar\,$(b)
Since the algebra $\calH_\mu$ is Frobenius, its socle contains the left module 
$\bbk$ with multiplicity 1. The unique 1-dimensional left ideal of $\calH_\mu$ 
is spanned by $x_\mu$. If $F$ is any free $\calH_n$-module, then $F$ is free 
also as an $\calH_\mu$-module, and therefore every homomorphism $N\to F$ sends 
$c'$ into $x_\mu F$. From a factorization $N\to F\to M$ of $\ph$ we deduce that 
$\ph(c')\in x_\mu M$.

(b)$\,\Rar\,$(c)
This is clear since 
$x_\mu=(T_i+1)\sum_{\si\in\calD(\<\tau_i\>\backslash\frS_\mu)}T_\si$ where 
$\<\tau_i\>$ is the 2-element subgroup of $\frS_\mu$ generated by $\tau_i$.

(c)$\,\Rar\,$(d)
Let us fix $i$ and construct the desired map $\psi_i$ by specifying its 
values on the basis elements $\{T_\si c'\mid\si\in\calD_\mu\}$ of $N$.

If $\si\in B_i(\mu)$, then $\tau_i\si=\si\tau_j$ for some $\tau_j\in\frB_\mu$. 
Hence $T_iT_\si=T_{\tau_i\si}=T_\si T_j$. The equality 
$\ph(v)=T_i\psi_i(v)+\psi_i(T_iv)+2\psi_i(v)$ for the element $v=T_\si c'$ is 
equivalent to $\ph(v)=(T_i+1)\psi_i(v)$ since $T_iv=-v$. An element 
$\psi_i(v)\in M$ satisfying the required equality can be found. Indeed,
$\ph(v)=T_\si\ph(c')\in(T_i+1)M$ since $\ph(c')\in(T_j+1)M$.

Suppose now that $\si\in A_i(\mu)$. In this case we put
$$
\psi_i(T_\si c')=0\qquad{\rm and}\qquad\psi_i(T_{\tau_i\si}c')=\ph(T_\si c').
$$
Then the equality $\ph(v)=T_i\psi_i(v)+\psi_i(T_iv)+2\psi_i(v)$ holds for 
$v=T_\si c'$ and also for $v=T_{\tau_i\si}c'=T_iT_\si c'$. In this way the 
value of $\psi_i$ has been determined on all basis elements of $N$.

(d)$\,\Rar\,$(c)
If $\tau_i\in\frB_\mu$, then $T_ic'=-c'$, whence
$$
\ph(c')=(T_i\psi_i+\psi_iT_i+2\psi_i)(c')=(T_i+1)\psi_i(c')\in(T_i+1)M.
$$

(c)$\,\Rar\,$(b)
For $\pi\in\muDla$ let $\frS_{\nu(\pi)}=\frS_\mu\cap\pi\frS_\la\pi^{-1}$. 
We have $M=\bigoplus_{\pi\in\,\muDla}M(\pi)$ where $M(\pi)$ is the 
$\calH_\mu$-submodule of $M$ with a basis 
$$
\{T_{\si\pi}c\mid\si\in\calD(\frS_\mu/\frS_{\nu(\pi)})\}.
$$
Since $\Hom_{\calH_\mu}\bigl(\bbk,M(\pi)\bigr)\cong\bbk$ for each $\pi$, the 
vector space $\Hom_{\calH_n}(N,M)$ has a basis $\{\ph_\pi\mid\pi\in\muDla\}$ 
where the homomorphism $\ph_\pi:N\to M$ is defined by the rule
$$ 
\ph_\pi(c')=\sum_{\si\in\calD(\frS_\mu/\frS_{\nu(\pi)})}T_{\si\pi}c\in M(\pi).
$$
We can write $\ph=\sum\al_\pi\ph_\pi$ with $\al_\pi\in\bbk$. If 
$\tau_i\in\frB_\mu$, then $(T_i+1)M$ is a direct sum of its 
subspaces $(T_i+1)M(\pi)$, whence the condition $\ph(c')\in(T_i+1)M$ implies 
that $\ph_\pi(c')\in(T_i+1)M$ for each $\pi\in\muDla$ such that $\al_\pi\ne0$.
If $\tau_i\in\frB_{\nu(\pi)}$ for some $\pi$, then $T_iT_\pi c=-T_\pi c$, 
and since the elements $T_\si c$ with $\si\in\calD_\la$, $\si\ne\pi$ span a 
$T_i$-invariant subspace, we deduce that $\ph_\pi(c')\notin(T_i+1)M$. In this 
case $\al_\pi=0$. On the other hand, if $\frS_{\nu(\pi)}=e$, then the 
$\calH_\mu$-module $M(\pi)$ is freely generated by $T_\pi c$, which shows 
that $\ph_\pi(c')=x_\mu T_\pi c\in x_\mu M$.

(b)$\,\Rar\,$(a)
Let $p\in M$ be such that $\ph(c')=x_\mu p$. Then $\ph$ is the composite of 
the homomorphism $N\to\calH_n$ sending $c'$ to $x_\mu$ and the homomorphism 
$\calH_n\to M$ sending $1$ to $p$.

In the course of the proof we have seen that $\ph$ satisfies the equivalent 
conditions (b) and (c) if and only if $\ph$ is a linear combination of the 
homomorphisms $\ph_\pi$ with $\pi\in\muDla$ and $\frS_{\nu(\pi)}=e$. This 
establishes the final conclusion.
\endproof

\section
2. Complexes associated with representations of $\calH_n$

Given an ordered collection of subspaces $U_1,\ldots,U_{n-1}$ of a vector space 
$M$, we denote by $K_\bullet\bigl(M;(U_i)\bigr)$ the complex of vector spaces
$$
0\mapr{}K_n\mapr{}\ldots\mapr{}K_i\mapr{\partial_i}K_{i-1}\mapr{}
\ldots\mapr{}K_0\mapr{}0
$$
constructed as follows. For $0\le i\le n$ put
$$
K_i=\Up_i/(\Up_i\cap\Si_i)\quad\hbox{where\quad $\Up_i=\bigcap_{j<i}U_j$\quad 
and\quad $\Si_i=\sum_{j>i}U_j$}
$$
with the convention that $\Up_0=\Up_1=M$ and $\Si_{n-1}=\Si_n=0$. The 
differentials $\partial_i$ are induced by the inclusions $\,\Up_i\sbs\Up_{i-1}$, 
$\,\Si_i\sbs\Si_{i-1}$.

Exactness of this complex gives an inductive step for the verification that 
the lattice of subspaces of $M$ generated by $U_1,\ldots,U_{n-1}$ is 
distributive \cite{Pol-P, Ch. 1, Prop. 7.2}. Complexes of this kind are 
responsible for Koszulness of the graded algebras, as discussed in the next 
section.

We are interested in the case when $M$ is a left module over the Hecke algebra 
$\calH_n=\calH_n(q)$ and the subspaces $U_1,\ldots,U_{n-1}$ are defined by one 
of the two conditions below:

\item(a)
\quad$U_i=\Ker\,(T_i-q)_M$ for each $i$,

\item(b)
\quad$U_i=\Im\,(T_i+1)_M$ for each $i$.

\noindent
where we denote by $x_M$ the linear operator by which an element $x\in\calH_n$ 
acts on $M$. Since $(T_i-q)(T_i+1)=0$, we have
$$
(T_i+1)M\sbs U_i\sbs\Ker\,(T_i-q)_M
$$
both in (a) and (b). If $q\ne-1$ then (a) is equivalent to (b), and so there 
is a difference between the two conditions only when $q=-1$.

For each $i=0,\ldots,n$ we have identified $\frS_i$ with the subgroup of 
$\frS_n$ generated by the set of basic transpositions $\{\tau_j\mid 0<j<i\}$. 
Denote by $\frS^\tri_i$ the subgroup of $\frS_n$ generated by 
$\{\tau_j\mid i<j<n\}$. Thus $\frS^\tri_i\cong\frS_{n-i}$. In particular, 
$\frS^\tri_i$ is the trivial subgroup $e$ for $i=n$ and for $i=n-1$.

\setitemsize(iii)
\proclaim
Lemma 2.1.
For $0\le i<n$ consider the following elements of $\calH_n\,${\rm:}
$$
x_i=\sum_{\si\in\calD_i}T_\si\,,\qquad
y_i=\sum_{\si\in\calD^\tri_i}
(-1)^{\ell(\si)}q^{n-1-i-\ell(\si)}T_\si.
$$
where $\calD_i=\calD(\frS_{i+1}/\frS_i)$ and 
$\calD^\tri_i=\calD(\frS^\tri_{i+1}\backslash\frS^\tri_i)$ are the sets of distinguished 
coset representatives. Then

\item(i)
$x_i$ maps $\Up_i$ to $\Up_{i+1}$ and $y_i$ maps $\Si_i$ to $\Si_{i+1},$

\item(ii)
$x_iy_i$ induces a linear map $s_i:K_i\to K_{i+1},$

\item(iii)
$\partial_{i+1}s_i+s_{i-1}\partial_i=[n]_q\Id,\,$ assuming that $s_{-1}=0$.

If\/ $[n]_q\ne0$ then the complex $K_\bullet\bigl(M;(U_i)\bigr)$ is exact.
\endproclaim

\Proof.
The inclusions $x_i\Up_i\sbs\Up_{i+1}$ and $y_i\Si_i\sbs\Si_{i+1}$ are special 
cases of Lemma 2.2 (see below) applied, respectively, to the pairs of Young 
subgroups $\frS_i\sbs\frS_{i+1}$ and $\frS^\tri_{i+1}\sbs\frS^\tri_i$. 
The subspace $U_k$ is stable under the action of any $T_j$ with $|j-k|>1$ 
since $T_jT_k=T_kT_j$. It follows that $\Up_i$ is stable under any $T_\si$ 
with $\si\in\frS^\tri_i$, and therefore $y_i\Up_i\sbs\Up_i$. On the other 
hand, $\Si_{i+1}$ is stable under any $T_\si$ with $\si\in\frS_{i+1}$, which 
yields $x_i\Si_{i+1}\sbs\Si_{i+1}$. Hence $x_iy_i$ maps $\Up_i$ to $\Up_{i+1}$ 
and $\Si_i$ to $\Si_{i+1}$. Thus (i) and (ii) have been checked.
Next, note that
$$
\eqalign{
\calD_i&{}=\{e,\,\tau_i,\,\tau_{i-1}\tau_i,\ldots,\,
\tau_1\tau_2\cdots\tau_i\},\cr
\calD^\tri_i&{}=\{e,\,\tau_{i+1},\,\tau_{i+1}\tau_{i+2},\ldots,\,
\tau_{i+1}\tau_{i+2}\cdots\tau_{n-1}\}\cr
}
$$
Suppose that $0<i<n$. Then $\calD_i=\{e\}\cup\calD_{i-1}\tau_i$ and 
$\calD^\tri_{i-1}=\{e\}\cup\tau_i\calD^\tri_i\,$. Since $T_e=1$, 
$\,T_{\si\tau_i}=T_\si T_i$ for $\si\in\calD_{i-1}$ and 
$T_{\tau_i\si}=T_iT_\si$ for $\si\in\calD^\tri_i$, we get
$$
\openup1\jot
\displaylines{
x_i=1+x_{i-1}T_i\,,\qquad y_{i-1}=q^{n-i}-T_iy_i\,,\cr
x_iy_i+x_{i-1}y_{i-1}=(1+x_{i-1}T_i)y_i+x_{i-1}\bigl(q^{n-i}-T_iy_i\bigr)
=y_i+q^{n-i}x_{i-1}\,.\cr
}
$$
If $j>i$ then $T_jv+v\in(T_j+1)M\sbs U_j\sbs\Si_i$ for all 
$v\in M$. Hence $T_\si v\equiv(-1)^{\ell(\si)}v$ modulo $\Si_i$ for 
$\si\in\calD^\tri_i$, and therefore
$$
y_iv\equiv{\textstyle\sum\limits_{k=0}^{n-i-1}}q^kv\quad\hbox{modulo }\Si_i\,.
$$
If $j<i$, then $T_jv=qv$ for all $v\in\Up_i$. Hence $T_\si v=q^{\ell(\si)}v$ 
for $\si\in\calD_{i-1}$, and therefore 
$\,x_{i-1}v=\sum\limits_{k=0}^{i-1}q^kv\,$. 
It follows that
$$
(x_iy_i+x_{i-1}y_{i-1})v=y_iv+q^{n-i}x_{i-1}v
\equiv{\textstyle\sum\limits_{k=0}^{n-1}}q^kv\quad\hbox{modulo $\Up_i\cap\Si_i$}
$$
for all $v\in\Up_i$. Since the map $\partial_{i+1}s_i+s_{i-1}\partial_i:K_i\to 
K_i$ is induced by the action of the element $x_iy_i+x_{i-1}y_{i-1}$ on $M$, 
this proves (iii) when $0<i<n$.

Note also that $x_0=1$, while $y_0v\equiv[n]_qv$ modulo $\Si_0$ for 
all $v\in M$. At the upper boundary $y_{n-1}=1$, while $x_{n-1}v=[n]_qv$ for 
all $v\in\Up_n$. Hence $\partial_1s_0=[n]_q\Id$ on $K_0$ and 
$s_{n-1}\partial_n=[n]_q\Id$ on $K_n$, which yields (iii) for $i=0$ and $i=n$.
The final conclusion is immediate from (iii).
\endproof

\proclaim
Lemma 2.2.
For each Young subgroup $\frS_\la$ of $\frS_n$ put
$$
\Up(\la)=\bigcap_{\{j\,\mid\,\tau_j\in\frB_\la\}}U_j\qquad{\rm and}\qquad
\Si(\la)=\sum_{\{j\,\mid\,\tau_j\in\frB_\la\}}U_j.
$$
Suppose that $\frS_\la$ and $\frS_\mu$ are two Young subgroups such that 
$\frS_\la\sbs\frS_\mu$. Then
$$
x\Up(\la)\sbs\Up(\mu)\qquad{\rm and}\qquad y\,\Si(\mu)\sbs\Si(\la)
$$
where \ $x=\sum_{\si\in\calD(\frS_\mu/\frS_\la)}T_\si$ \ and \ 
$y=\sum_{\si\in\calD(\frS_\la\backslash\frS_\mu)}
(-1)^{\ell(\si)}q^{m-\ell(\si)}T_\si$ \ with\hfill\break 
$m=\max\{\ell(\si)\mid\si\in\calD(\frS_\la\backslash\frS_\mu)\}\strut$.
\endproclaim

\Proof.
Fix any basic transposition $\tau_i\in\frB_\mu$. As in section 1 we have
$$
\calD(\frS_\mu/\frS_\la)=A\cup\tau_iA\cup B,\qquad
\calD(\frS_\la\backslash\frS_\mu)=A^{-1}\cup A^{-1}\tau_i\cup B^{-1}
$$
(disjoint unions) where $A=\frS_\mu\cap A_i(\la)$ and 
$B=\frS_\mu\cap B_i(\la)$. Hence
$$
\openup1\jot
\eqalign{
x&{}=(T_i+1)\sum_{\si\in A}T_\si+\sum_{\si\in B}T_\si\,,\cr
y&{}=\sum_{\si\in A^{-1}}(-1)^{\ell(\si)}q^{m-1-\ell(\si)}T_\si(q-T_i)
+\sum_{\si\in B^{-1}}(-1)^{\ell(\si)}q^{m-\ell(\si)}T_\si\,.\cr
}
$$
Recall that $(T_i+1)M\sbs U_i$. Suppose that $\si\in B$. Then 
$\tau_i\si=\si\tau_j>\si$ for some $\tau_j\in\frB_\la$, and therefore 
$T_iT_\si=T_\si T_j$. For each $v\in\Up(\la)$ we have $T_iT_\si v=qT_\si v$ 
since $T_jv=qv$. Moreover, in case (b) $v=(T_j+1)w$ for some $w\in M$, and 
then $T_\si v=(T_i+1)T_\si w$. It follows that $T_\si v\in U_i$ both in (a) 
and (b). Since this inclusion holds for each $\si\in B$, we deduce that 
$x\Up(\la)\sbs U_i$.

Similarly, if $\si\in B^{-1}$, then $T_\si T_i=T_jT_\si$ for some $j$ such 
that $\tau_j\in\frB_\la$. For each $v\in U_i$ we have $T_iv=qv$, whence 
$T_jT_\si v=qT_\si v$. Moreover, in case (b) $v=(T_i+1)w$ for some $w\in M$, 
and then $T_\si v=(T_j+1)T_\si w$. Hence $T_\si v\in U_j\sbs\Si(\la)$ both in 
(a) and (b). Since $q-T_i$ annihilates $U_i$, it follows that 
$yU_i\sbs\Si(\la)$.
\endproof

\proclaim
Lemma 2.3.
Let $\chi:\calH_\la\to\bbk$ be a $1$-dimensional representation of a parabolic 
subalgebra $\calH_\la$ of $\calH_n$. Consider the standard basis 
$\{v_\si\mid\si\in\calD_\la\}$ of the induced $\calH_n$-module 
$M=\calH_n\ot_{\calH_\la}\bbk(\chi)$ where $v_\si=T_\si c$ in the notation of 
section\/ $1$. Then\/ $\Ker\,(T_i-q)_M$ is spanned by the elements
$$
\{v_{\tau_i\si}+v_\si\mid\si\in A_i(\la)\}\cup
\{v_\si\mid\si\in B_i(\la){\rm\ and\ }\chi_\si(T_i)=q\}
$$
and $\,\Im\,(T_i+1)_M$ is spanned by the elements
$$
\{v_{\tau_i\si}+v_\si\mid\si\in A_i(\la)\}\cup
\{v_\si\mid\si\in B_i(\la){\rm\ and\ }\chi_\si(T_i)\ne-1\}.
$$
\endproclaim

\Proof.
The module $M$ is a direct sum of $T_i$-invariant subspaces $M(\si)$ with 
$\si$ in $A_i(\la)\cup B_i(\la)$ where $M(\si)$ is spanned by two elements 
$v_\si$, $v_{\tau_i\si}$ for $\si\in A_i(\la)$ and by the single element 
$v_\si$ for $\si\in B_i(\la)$. Obviously $\Ker\,(T_i-q)_M$ and $\Im\,(T_i+1)_M$ 
are sums of their intersections with those subspaces. From the formulas for the 
action of $T_i$ on $M$ it is clear that
$$
\Ker\,(T_i-q)_M\cap M(\si)
=\Im\,(T_i+1)_M\cap M(\si)=\bbk(v_{\tau_i\si}+v_\si)
$$
when $\si\in A_i(\la)$. If $\si\in B_i(\la)$, then $v_\si$ is an eigenvector 
for the operator $(T_i)_M$ with the eigenvalue $\chi_\si(T_i)$. In this case 
$v_\si\in\Ker\,(T_i-q)_M$ if and only if $\chi_\si(T_i)=q$, and 
$v_\si\in\Im\,(T_i+1)_M$ if and only if $\chi_\si(T_i)\ne-1$.
\endproof

By Lemma 2.1 the complex $K_\bullet\bigl(M;(U_i)\bigr)$ is exact for any left 
$\calH_n$-module $M$ when $[n]_q\ne0$. If $q=0$, then this result does apply 
since $[n]_q=1$ in this case. For arbitrary $q$ we have to restrict the class 
of $\calH_n$-modules:

\proclaim
Proposition 2.4.
Let $M$ be a finite dimensional $\calH_n$-module whose indecomposable direct 
summands all have $1$-dimensional sources. With subspaces $U_1,\ldots,U_{n-1}$ 
defined by either\/ {\rm(a)} or\/ {\rm(b)} the complex 
$K_\bullet\bigl(M;(U_i)\bigr)$ is exact.
\endproclaim

\Proof.
Note that each $U_i$ depends on $M$ functorially, and therefore the 
construction of $K_\bullet\bigl(M;(U_i)\bigr)$ gives a functor from the 
category of $\calH_n$-modules to the category of complexes. Since this functor 
is additive, the conclusion of Proposition 2.4 holds for any given 
$\calH_n$-module $M$ if and only if the conclusion holds for each 
indecomposable direct summand of $M$. This shows that it suffices to give the 
proof assuming that $M=\calH_n\ot_{\calH_\la}\bbk(\chi)$ for some 
$1$-dimensional representation $\chi:\calH_\la\to\bbk$ of a parabolic 
subalgebra of $\calH_n$.

By Lemma 2.1 the conclusion is true when $q=0$. Suppose that $q\ne0$. Consider 
the parabolic subalgebra $\calH_\la(0)$ of the 0-Hecke algebra $\calH_n(0)$ 
corresponding to the same composition $\la$ of $n$. Let 
$\oT_1,\ldots,\oT_{n-1}$ stand for the canonical generators of $\calH_n(0)$. 
There is a $1$-dimensional representation $\xi:\calH_\la(0)\to\bbk$ defined 
on the generators $\{\oT_i\mid0<i<n,\ \tau_i\in\frS_\la\}$ of $\calH_\la(0)$ 
as follows. When $q\ne-1$ set
$$
\xi(\oT_i)=0\hbox{ \ if }\chi(T_i)=q\qquad\hbox{and}\qquad
\xi(\oT_i)=-1\hbox{ \ if }\chi(T_i)=-1.
$$
When $q=-1$ set
$$
\xi(\oT_i)=0\hbox{ \ in case (a)}\qquad\hbox{and}\qquad
\xi(\oT_i)=-1\hbox{ \ in case (b)}.
$$
Let $\{v_\si\mid\si\in\calD_\la\}$ be the standard basis of $M$ and
$\{v^0_\si\mid\si\in\calD_\la\}$ a similar basis of the induced 
$\calH_n(0)$-module $M^0=\calH_n(0)\ot_{\calH_\la(0)}\bbk(\xi)$. The 
assignments $v_\si\mapsto v^0_\si$ define a linear isomorphism $M\cong M^0$. 
In view of Lemma 2.3 the subspace $U_i$ of $M$ is mapped onto a similar 
subspace $U^0_i=\{v\in M^0\mid\oT_iv=0\}$ of $M^0$, for each $i=1,\ldots,n-1$. 
Hence the complex $K_\bullet\bigl(M;(U_i)\bigr)$ is isomorphic to 
$K_\bullet\bigl(M^0;(U^0_i)\bigr)$. But the latter complex is exact, as we 
have observed already.
\endproof

\proclaim
Corollary 2.5.
Under the same assumption about $M$ the complex $K_\bullet\bigl(M;(U_i)\bigr)$ 
is exact also when $U_i=(T_i-q)M$ for each $i=1,\ldots,n-1$.
\endproclaim

\Proof.
Let $\Mt$ be $M$ with the $\calH_n$-module structure twisted by the automorphism 
of $\calH_n$ sending $T_i$ to $q-1-T_i$ for each $i$. The $\calH_n$-module $\Mt$ 
has the same submodules as $M$, but with the twisted action of $\calH_n$. So 
it follows from Lemma 1.5 that all indecomposable direct summands of $\Mt$ 
have $1$-dimensional sources provided this holds for $M$. Since $T_i-q$ acts 
on $M$ as $-(T_i+1)$ acts on $\Mt$, we have $U_i=(T_i+1)\Mt$. Therefore 
Corollary 2.5 follows from Proposition 2.4 applied to $\Mt$.
\endproof

If $[n]_q\ne0$ then the conclusion of Corollary 2.5 holds without any 
restriction on $M$ in view of Lemma 2.1. In particular, this is true for 
$q=0$.

\section
3. Koszulness of the $R$-symmetric algebras

Let $A=\bigoplus_{n=0}^\infty A_n$ be a quadratic graded algebra generated by 
some vector space $V=A_1$. This means that $A\cong\bbT(V)/I$ where $I$ is the 
ideal of the tensor algebra $\bbT(V)$ generated by a vector subspace 
$U\sbs V^{\ot2}$. The books \cite{Man} and \cite{Pol-P} provide general 
reference on quadratic algebras. For each $n>1$ and $0<i<n$ put
$$
U_i^{(n)}=V^{\ot(i-1)}\ot U\ot V^{\ot(n-i-1)}\sbs V^{\ot n}=\bbT_n(V)
$$
and $\Up^{(n)}=\bigcap_{i=1}^{n-1}U_i^{(n)}$. Put also $\Up^{(0)}=\bbk$ and 
$\Up^{(1)}=V$. The right \emph{Koszul complex} $K_\bullet(A)$ is the complex 
of right $A$-modules
$$
\ldots\mapr{}\Up^{(i)}\ot A\mapr{\partial_i}\Up^{(i-1)}\ot A\mapr{}\ldots
\mapr{}V\ot A\mapr{}A\mapr{}0
$$
where $\partial_i\,$, for each $i>0$, is the restriction of the $A$-linear map
$$
V^{\ot i}\ot A\mapr{}V^{\ot(i-1)}\ot A,\qquad (t\ot v)\ot a\mapsto t\ot va
$$
for $t\in V^{\ot(i-1)}$, $v\in V$ and $a\in A$. The grading of $A$ gives rise 
to a decomposition of $K_\bullet(A)$ into a direct sum of subcomplexes
$$
K_\bullet^{(n)}(A):\quad
0\mapr{}\Up^{(n)}\mapr{}\Up^{(n-1)}\ot A_1\mapr{}\ldots
\mapr{}V\ot A_{n-1}\mapr{}A_n\mapr{}0\,.
$$
There is an isomorphism of complexes $K_\bullet^{(n)}(A)\cong 
K_\bullet\bigl(V^{\ot n};(U_i^{(n)})\bigr)$, the latter having been defined in 
section 2. Indeed, letting \ $\Up_i^{(n)}=\bigcap_{j<i}U_j^{(n)}$, 
\ $\Si_i^{(n)}=\sum_{j>i}U_j^{(n)}$, we have
$$
\Up_i^{(n)}=\Up^{(i)}\ot V^{\ot(n-i)},\quad 
\Si_i^{(n)}=V^{\ot i}\ot\Si_0^{(n-i)},\quad 
\Up_i^{(n)}\cap\Si_i^{(n)}=\Up^{(i)}\ot\Si_0^{(n-i)}.
$$
Since $\Si_0^{(k)}=\sum_{j=1}^{k-1}U_j^{(k)}$ coincides with the $k$th 
homogeneous component of the ideal $I$, it follows that 
$V^{\ot k}/\Si_0^{(k)}\cong A_k$ for each $k$. Hence
$$
\Up_i^{(n)}\!/(\Up_i^{(n)}\cap\Si_i^{(n)})\cong\Up^{(i)}\ot A_{n-i}\,.
$$
It is easy to see that these linear isomorphisms between the homogeneous 
components of the two complexes are compatible with the differentials.

Note that $H_0\bigl(K_\bullet(A)\bigr)\cong A/A_+\cong\bbk$ where 
$A_+=VA=\sum_{n>0}A_n$. The algebra $A$ is said to be \emph{Koszul}\/ if the 
complex $K_\bullet(A)$ is acyclic in all positive degrees, i.e. $K_\bullet(A)$ 
is a resolution of the trivial right $A$-module $\bbk$. There are several 
equivalent characterizations of this property (see \cite{Pol-P}). By a 
fundamental result of Backelin \cite{Back81} $A$ is Koszul if and only if 
$U_1^{(n)},\ldots,U_{n-1}^{(n)}$ generate a distributive lattice of subspaces 
of $V^{\ot n}$ for each $n>1$. From this it is easy to see that Koszulness of 
an algebra is a left-right symmetric property.

Assume further on that $\dim V<\infty$. Then $\dim A_n<\infty$ for all $n$. 
The \emph{Hilbert series} $h_A(t)$ of $A$ is the formal power series in one 
indeterminate $t$ whose coefficients are the dimensions of the homogeneous 
components $A_n$:
$$
h_A(t)=\sum\,(\dim A_n)\,t^n.
$$
For each $n$ identify $\bbT_n(V^*)$ with the dual of the vector space 
$\bbT_n(V)$ using the bilinear pairing
$$
\<f_1\ot\ldots\ot f_n,\,v_1\ot\ldots\ot v_n\>=\prod f_i(v_i).
$$
The \emph{quadratic dual} $A^!$ of $A$ is the factor algebra of the tensor 
algebra $\bbT(V^*)$ by the ideal generated by the subspace
$$
U^\perp=\{f\in\bbT_2(V^*)\mid\<f,U\>=0\}\sbs\bbT_2(V^*).
$$
The $n$th homogeneous component of this ideal is then the subspace
$$
\sum_{i=1}^{n-1}\,(U_i^{(n)})^\perp=(\Up^{(n)})^\perp\sbs\bbT_n(V^*).
$$
Hence $A^!_n=\bbT_n(V^*)/(\Up^{(n)})^\perp\cong(\Up^{(n)})^*$, and it follows 
that
$$
h_{A^!}(t)=\sum\,(\dim\Up^{(n)})\,t^n.
$$
If $A$ is Koszul, then $\sum_{i=0}^n(-1)^i(\dim A_i)(\dim\Up^{(n-i)})=0$ for 
each $n>0$, which entails the well-known relation between the Hilbert series 
of $A$ and $A^!\,$:
$$
h_A(-t)\,h_{A^!}(t)=1.
$$

Let $R$ be a Hecke symmetry on a vector space $V$, and let $0\ne q\in\bbk$ be 
the parameter of the Hecke relation satisfied by $R$. For each $n\ge0$ there is 
a representation of the Hecke algebra $\calH_n=\calH_n(q)$ in $\bbT_n(V)$ such 
that each generator $T_i\,$, $\,0<i<n$, acts by means of the linear operator
$$
R_i^{(n)}=\Id^{\ot(i-1)}\!\ot\,R\ot\Id^{\ot(n-i-1)}. 
$$
Recall that $\calH_0=\calH_1=\bbk$. We assume that $\dim V<\infty$.

The algebras $\bbS(V,R)$ and $\La(V,R)$ are defined as the factor algebras 
of $\bbT(V)$ by the ideals generated, respectively, by the subspaces
$$
\Im\,(R-q\cdot\Id)\quad{\rm and}\quad\Ker\,(R-q\cdot\Id)\quad{\rm of\ }\bbT_2(V).
$$

\proclaim
Theorem 3.1.
Suppose that $R$ satisfies the $1$-dimensional source condition. Then the 
$R$-symmetric algebra $\bbS(V,R)$ and the $R$-skewsymmetric algebra 
$\La(V,R)$ are Koszul. Their Hilbert series satisfy the relation 
\ $h_{\bbS(V,R)}(t)\,h_{\La(V,R)}(-t)=1$.
\endproclaim

\Proof.
In terms of the $\calH_n$-module structure arising from $R$ 
the previously defined subspaces of $\bbT_n(V)$ are
$$
U_i^{(n)}=\cases{
\Im\,(R_i^{(n)}-q\cdot\Id)=\Im\,(T_i-q)_{\bbT_n(V)} & 
when $A=\bbS(V,R)$,\cr
\noalign{\smallskip}
\Ker\,(R_i^{(n)}-q\cdot\Id)=\Ker\,(T_i-q)_{\bbT_n(V)}& 
when $A=\La(V,R)$.\cr
}
$$
The assumption about $R$ means that each indecomposable direct summand of the 
$\calH_n$-module $\bbT_n(V)$ has a $1$-dimensional source. Hence the complex 
$K_\bullet\bigl(V^{\ot n},(U_i^{(n)})\bigr)$ is exact, for each $n>0$, by 
Proposition 2.4 and Corollary 2.5. Hence so is the isomorphic complex 
$K_\bullet^{(n)}(A)$. Since $K_\bullet^{(0)}(A)$ is concentrated in 
degree 0, it follows that 
$H_i\bigl(K_\bullet(A)\bigr)=H_i\bigl(K_\bullet^{(0)}(A)\bigr)=0$ for all 
$i>0$.

Let $\Up^{(n)}=\bigcap_{i=1}^{n-1}U_i^{(n)}$ be defined with respect to 
$A=\La(V,R)$, i.e.
$$
\Up^{(n)}=\{u\in\bbT_n(V)\mid T_iu=qu\hbox{ for each }i=1,\ldots,n-1\}.
$$
This is the largest subspace of $\bbT_n(V)$ on which $\calH_n$ operates 
trivially. On the other hand, $\bbS_n(V,R)$ is the largest factor space of 
$\bbT_n(V)$ on which $\calH_n$ operates trivially. It follows that
$$
\eqalign{
\dim\La^!_n(V,R)&{}=\dim\Up^{(n)}=\dim\Hom_{\calH_n}(\bbk\triv\,,V^{\ot n}),\cr
\dim\,\bbS_n(V,R)&{}=\dim\Hom_{\calH_n}(V^{\ot n},\,\bbk\triv),\cr
}
$$
whence $\dim\La^!_n(V,R)=\dim\bbS_n(V,R)$ according to Lemma 1.2. The relation 
between the Hilbert series of $\bbS(V,R)$ and $\La(V,R)$ follows now from the 
relation between $h_A$ and $h_{A^!}$.
\endproof

There are several transformations of the Hecke symmetry $R$. Put 
$$
\Rt=(q-1)\Id-R=-qR^{-1}\qquad{\rm  and}\qquad R\op=\tau R\,\tau
$$
where $\tau$ is the flip operator $v_1\ot v_2\mapsto v_2\ot v_1$ on $\bbT_2(V)$. 
Let $R^*$ be the linear operator on $\bbT_2(V^*)\cong\bbT_2(V)^*$ adjoint to $R$.
Given an algebra $A$, we denote by $A\op$ the algebra with the same set of 
elements but with the opposite multiplication.

\proclaim
Lemma 3.2.
The operators $\Rt,$ $R\op,$ $R^*$ are Hecke symmetries with the same parameter 
$q$ as $R$. We have
$$
\tabskip=0pt
\openup1\jot
\vcenter{\halign{\hfil$#$&$#$\hfil&\qquad\hfil$#$&$#$\hfil\cr
\bbS(V,R\op) & {}\cong\bbS(V,R)\op, & \bbS(V^*\!,R^*) & {}=\La(V,R)^!, \cr
\La(V,R\op) & {}\cong\La(V,R)\op, & \La(V^*\!,R^*) & {}=\bbS(V,R)^!. \cr
}}
$$
If\/ $q\ne-1,$ then $\,\bbS(V,\Rt)=\La(V,R)\,$ and $\,\La(V,\Rt)=\bbS(V,R)\,$.
\endproclaim

All verifications are straightforward. Note also that these Hecke symmetries 
$\Rt$, $R\op$, $R^*$ satisfy the $1$-dimensional source condition provided so 
does $R$. If $R$ satisfies the trivial source condition, so do $R\op$ and $R^*$. 
The inverse operator $R^{-1}$ is a Hecke symmetry with parameter $q^{-1}$ 
giving rise to the same pair of quadratic graded algebras as the pair 
$\bbS(V,R)$, $\La(V,R)$ obtained from $R$.

\medbreak
We end this section with an example showing that $\bbS(V,R)$, $\La(V,R)$ are 
not always Koszul. Let $V$ be a 2-dimensional vector space with a basis $x,y$. 
Assume that $\chr\bbk\ne2$. We start with the $R$-matrix $R_{H0.2}$ in the 
notation of Hietarinta \cite{Hie93, p. 1732}:
$$
\pmatrix{1&0&0&1\cr 0&1&1&0\cr 0&1&-1&0\cr -1&0&0&1}
$$
In a slightly different form an equivalent matrix appeared under the label 
$R_1$ in the list of Hlavat\'y \cite{Hla87, p. 1663}. This matrix represents an 
operator satisfying the quantum Yang-Baxter equation. Composing with the flip 
of tensorands we obtain a diagonalizable operator satisfying the braid equation 
whose characteristic polynomial is $(t^2-2t+2)^2$. A final scaling yields the 
matrix
$$
{q-1\over2}\pmatrix{1&0&0&1\cr 0&1&-1&0\cr 0&1&1&0\cr -1&0&0&1\cr},\qquad 
q^2=-1,
$$
of a Hecke symmetry with eigenvalues $-1$, $q$ where $q$ is a primitive 4th 
root of 1. The matrix is written in the basis $x^2,xy,yx,y^2$ of $\bbT_2(V)$. 
One eigenspace of $R$ is spanned by $y^2-qx^2$, $xy-qyx$. It gives the defining 
relations $y^2=qx^2$, $xy=qyx$ of $\bbS(V,R)$. The first relation shows that 
$y^2$ is central in $\bbS(V,R)$. But $xy^2=-y^2x$ according to the second 
relation, whence $xy^2=0$. From this it is clear that $\bbS_3(V,R)=0$. 

Similarly, the algebra $\La(V,R)$ has the defining relations $x^2=qy^2$, 
$yx=qxy$. It is isomorphic to $\bbS(V,R)$. Hence
$$
h_{\bbS(V,R)}(t)=h_{\La(V,R)}(t)=1+2t+2t^2.
$$
Moreover, the quadratic dual algebras are isomorphic to the original ones. 
Thus the standard relation between their Hilbert series is not satisfied.

In this example $\calH_4$ is the first nonsemisimple algebra in the family of 
Hecke algebras. Since all its proper parabolic subalgebras are semisimple, the 
indecomposable $\calH_4$-modules with a 1-dimensional source are either 
1-dimensional or projective. It can be checked that the $\calH_4$-module 
$V^{\ot4}$ is a direct sum of simple 2-dimensional submodules. There are two 
nonisomorphic simple modules of dimension 2. One of them is projective, but 
the other is not. Thus $R$ does not satisfy the 1-dimensional source condition.

\section
4. Nondegeneracy of the multiplication maps

The aim of this section is to prove that the $R$-skewsymmetric algebra 
$\La=\La(V,R)$ is Frobenius under suitable assumptions. Recall that $\La=\bbT/I$ 
where $\bbT=\bbT(V)$ is the tensor algebra of $V$ and $I$ is its homogeneous 
ideal generated by the subspace $U=\Ker\,(R-q\cdot\Id)$ of $\bbT_2$. We will 
be omitting the sign $\ot$ when referring to the multiplication in $\bbT$.

The next lemma provides the main step in tackling the problem. 
Recall that we denote by $\frS_{1,k-1}$ the subgroup of $\frS_k$ generated 
by $\{\tau_i\mid1<i<k\}$.

\proclaim
Lemma 4.1.
Fix some $n>1$ and put
$$
\openup1\jot
\displaylines{
L_k=\{a\in\bbT_k\mid a\bbT_{n-k}\sbs I_n\}\qquad\hbox{\rm for }0\le k<n,\cr
y_k=\sum_{\si\in\calD(\frS_{1,k-1}\backslash\frS_k)}
(-1)^{\ell(\si)}q^{k-1-\ell(\si)}T_\si\in\calH_k\qquad\hbox{\rm for each }k>0.\cr
}
$$
If\/ $0<k<n,$ then $y_k$ maps $L_k$ into $VL_{k-1}$.
\endproclaim

\Proof.
Note that $VL_{k-1}=\{b\in\bbT_k\mid b\bbT_{n-k+1}\sbs VI_n\}$. Thus we have 
to show that $(y_kL_k)\bbT_{n-k+1}\sbs VI_n$. Let $a\in L_k$. Then
$a\bbT_{n-k}\sbs I_n$, whence
$$
a\bbT_{n-k+1}=a\bbT_{n-k}V\sbs I_nV\sbs I_{n+1}.
$$ 
We will work inside the $\calH_{n+1}$-module $\bbT_{n+1}$. In conformance with 
the notation of section 2 put $U_i=\Ker\,(T_i-q)_{\bbT_{n+1}}$ for each 
$i=1,\ldots,n$. Note that
$$
I_{n+1}=\sum_{i=1}^nU_i\,,\qquad
VI_n=\sum_{i=2}^nU_i\,,\qquad
I_nV=\sum_{i=1}^{n-1}U_i\,.
$$
By Lemma 2.2 applied with $\,\frS_\la=\frS_{1,n}$, $\,\frS_\mu=\frS_{n+1}\,$ 
and $n$ replaced by $n+1$ we get $\,y_{n+1}I_{n+1}\sbs VI_n\,$. Therefore
$$ 
y_{n+1}(a\bbT_{n-k+1})\sbs VI_n\,.
$$
Since $\calD(\frS_{1,k-1}\backslash\frS_k)=\{e,\tau_1,\tau_1\tau_2\,,\ldots,
\tau_1\tau_2\cdots\tau_{k-1}\}$, we have
$$
y_k=\sum_{i=0}^{k-1}(-1)^iq^{k-1-i}p_i\quad\hbox{where $p_0=1$, $p_1=T_1$, 
$p_i=T_1T_2\cdots T_i$ for $i>1$}.
$$
In particular, $y_{n+1}=\sum_{i=0}^n(-1)^iq^{n-i}p_i$.
As follows immediately from the braid relations between $T_1,\ldots,T_n$, the 
element $p_n$ has the property that $T_{i+1}p_n=p_nT_i$, and therefore
$$
(T_{i+1}-q)p_n=p_n(T_i-q),
$$
for each $i=1,\ldots,n-1$. Since $p_n$ is invertible in $\calH_{n+1}$, we deduce 
that $p_nU_i=U_{i+1}$ for those values of $i$. Hence $p_n$ maps $I_nV$, and in 
particular the subspace $a\bbT_{n-k+1}$, into $VI_n$.

If $i<n$, then $p_i=p_nT_n^{-1}\cdots T_{i+1}^{-1}$. For each $j=k+1,\ldots,n$ 
the element $T_j$ operates on the second space in the decomposition 
$\bbT_{n+1}=\bbT_k\ot\bbT_{n-k+1}$, which implies that 
$a\bbT_{n-k+1}=a\ot\bbT_{n-k+1}$ is stable under the action of $T_j$ and 
$T_j^{-1}$. It follows that
$$ 
p_i(a\bbT_{n-k+1})\sbs VI_n\quad\hbox{whenever $k\le i\le n$}.
$$
Now $\,y_k=q^{k-1-n}(y_{n+1}-\sum_{i=k}^n(-1)^iq^{n-i}p_i)$. The previous 
inclusions entail
$$ 
y_k(a\bbT_{n-k+1})\sbs VI_n\,.
$$
Since for each $j=1,\ldots,k-1$ the element $T_j$ operates on the 
first space in the decomposition $\bbT_{n+1}=\bbT_k\ot\bbT_{n-k+1}$, so too 
does $y_k$. Hence
$$
(y_ka)\bbT_{n-k+1}=y_k(a\bbT_{n-k+1})\sbs VI_n\,,
$$
yielding $\,y_ka\in VL_{k-1}$.
\endproof

The reader should note that the proof of Lemma 4.1 uses only the braid relations 
between $T_1,\ldots,T_n$. Therefore Lemma 4.1 holds more generally when $R$ is a 
linear operator on $V^{\ot2}$ satisfying the braid equation but not necessarily 
the quadratic Hecke equation, and $q$ is any eigenvalue of $R$ used in the 
definition of $\La$. Such an operator (called a Yang-Baxter operator in the 
literature) gives rise to representations of the Artin braid groups $B_k$.  
The elements $T_\si$ with $\si\in\frS_k$ make sense in $B_k$, and the element 
$y_k$ is defined for each $k$ in the group algebra of $B_k$. This observation 
will be essential later (see Lemma 6.5).

\proclaim
Lemma 4.2.
Suppose that $\La_n\ne0$ for some $n>1$. If $0\ne a\in\La_k$ where either 
$k=1$ or $k=2$, then $\,a\La_{n-k}\ne0\,$ and $\,\La_{n-k}\,a\ne0$.
\endproclaim

\Proof.
In the notation of Lemma 4.1 $L_0=0$ since $\bbT_0=\bbk$ and $\bbT_n\ne I_n$. 
By Lemma 4.1 $y_1L_1\sbs VL_0=0$ and $y_2L_2\sbs VL_1$. Since $y_1=1$, we get 
$L_1=0$. Hence $y_2L_2=0$. This means that $L_2\sbs\Ker\,(T_1-q)_{\bbT_2}=I_2$ 
since $y_2=q-T_1$. In fact $L_2=I_2$ since the opposite inclusion is obvious.
Thus $L_k=I_k$ for $k=1$ and for $k=2$. Passing to the factor algebra 
$\La=\bbT/I$, we deduce that $\La_{n-k}$ has zero left annihilator in $\La_k$. 
In view of Lemma 3.2 we can apply this conclusion also to $\La\op$. Hence 
$\La_{n-k}$ has zero right annihilator in $\La_k$.
\endproof

\proclaim
Corollary 4.3.
Suppose that\/ $\dim\La_n=1$ and $\La_{n+1}=0$. If $n=2$ or $n=3$, then $\La$ 
is a Frobenius algebra.
\endproclaim

\Proof.
By Lemma 4.2 the bilinear pairings $\La_k\times\La_{n-k}\to\La_n$ arising from 
the multiplication in $\La$ are nondegenerate for $k=1$ and $k=2$.
\endproof

\proclaim
Lemma 4.4.
Let $M$ be a finite dimensional $\calH_n$-module whose indecomposable direct 
summands all have trivial sources. Let $\frS_\la$ be a Young subgroup of 
$\frS_n$ and
$$
y=\sum_{\si\in\calD(\frS_\la\backslash\frS_n)}
(-1)^{\ell(\si)}q^{m-\ell(\si)}T_\si\quad
\hbox{with $\,m=\max\{\ell(\si)\mid\si\in\calD(\frS_\la\backslash\frS_n)\}$}.
$$
Put\/ $\Si(\la)=\sum_{\{i\,\mid\,\tau_i\in\frB_\la\}}U_i$ and\/ 
$\Si(n)=\sum_{i=1}^{n-1}U_i$ where $U_i=\Ker\,(T_i-q)_M$ for each $i$. Then
$$
\{a\in M\mid ya\in\Si(\la)\}=\Si(n).
$$
\endproclaim

\Proof.
The inclusion $y\,\Si(n)\sbs\Si(\la)$ is a special case of Lemma 2.2. So we have 
only to prove that $a\in\Si(n)$ for each $a\in M$ such that $ya\in\Si(\la)$. 
This assertion holds for any given $\calH_n$-module $M$ if and only if it holds 
for each indecomposable direct summand of $M$. Therefore we may assume that 
$M=\calH_n\ot_{\calH_\nu}\bbk\triv$ where $\calH_\nu$ is a parabolic subalgebra 
of $\calH_n$.

If $\calH_\nu\ne\bbk$, then $T_j\in\calH_\nu$ for some $j$. In this case 
$T_jc=qc$ where $c$ is the canonical generator of $M$, whence 
$c\in U_j\sbs\Si(n)$. Observing that $\,(T_i+1)M\sbs U_i\sbs\Si(n)\,$ for each 
$i=1,\ldots,n-1$, we see that $\Si(n)$ is stable under the action of all 
$T_1,\ldots,T_{n-1}$, i.e. $\Si(n)$ is an $\calH_n$-submodule of $M$. But then 
$\Si(n)=M$, and the desired conclusion is obviously true.

It remains to consider the case when $\calH_\nu=\bbk$, and therefore 
$M=\calH_n$ is a cyclic free $\calH_n$-module. Since $\calH_n$ is a free 
module over its subalgebra $\bbk+\bbk T_i$, we deduce that $U_i=(T_i+1)M$ 
for each $i=1,\ldots,n-1$. Hence $M/\Si(n)$ is the largest factor module of $M$ 
on which each $T_i$ operates as minus identity transformation. Clearly 
$$
\bigl(M/\Si(n)\bigr)^*\cong\Hom_{\calH_n}(M,\bbk\alt)\cong\bbk,
$$
which shows that $\Si(n)$ is a subspace of codimension 1 in $M$. 

The subspace $Y=\{a\in M\mid ya\in\Si(\la)\}$ contains $\Si(n)$ by Lemma 2.2. 
Suppose that $Y\ne\Si(n)$. Then we must have $Y=M$. In particular, $1\in Y$, 
which means that $y\in\Si(\la)$. However, 
$$
\Si(\la)=\sum_{\{i\mid\tau_i\in\frB_\la\}}(T_i+1)M=JM=J\calH_n
$$
where $J$ is the ideal of $\calH_\la$ generated by 
$\{T_i+1\mid\tau_i\in\frB_\la\}$, i.e. $J$ is the annihilator of the alternating 
representation of $\calH_\la$. Since $\calH_n$ is a free left $\calH_\la$-module 
with a basis $\{T_\si\mid\si\in\calD(\frS_\la\backslash\frS_n)\}$, the 
inclusion $y\in J\calH_n$ entails $1\in J$, which is impossible. This 
contradiction proves that $Y=\Si(n)$.
\endproof

\proclaim
Theorem 4.5.
Suppose that $R$ satisfies the trivial source condition. Suppose also that\/ 
$\dim\La_n(V,R)=1$ and $\La_{n+1}(V,R)=0$ for some $n>0$. Then\/ $\La(V,R)$ is 
a Frobenius algebra{\rm,} while\/ $\bbS(V,R)$ is a Gorenstein algebra of 
global dimension $n$.
\endproclaim

\Proof.
The left kernel of the bilinear pairing $\La_k\times\La_{n-k}\to\La_n$ arising 
from the multiplication in $\La$ is nothing else but the image of $L_k$ in 
$\La_k=\bbT_k/I_k$ where $L_k$ is the subspace of $\bbT_k$ introduced in Lemma 
4.1. To show that the left kernel vanishes we have to prove that $L_k=I_k$. 
But this can be done by induction on $k$. Indeed, if $L_{k-1}=I_{k-1}$, then 
$y_kL_k\sbs VI_{k-1}$ by Lemma 4.1. Note that 
$$
I_k=\sum_{i=1}^{k-1}U_i^{(k)},\qquad
VI_{k-1}=\sum_{i=2}^{k-1}U_i^{(k)}\quad{\rm where\ }
U_i^{(k)}=\Ker\,(T_i-q)_{\bbT_k}.
$$
The indecomposable direct summands of the $\calH_k$-module $M=\bbT_k$ all have 
trivial sources by the assumption about $R$, and we can apply Lemma 4.4 with 
$n$ replaced by $k$ and $\frS_\la=\frS_{1,k-1}$. In this case $y=y_k$, 
$\Si(\la)=VI_{k-1}$, $\Si(k)=I_k$, and the conclusion of Lemma 4.4 gives 
the desired inclusion $L_k\sbs I_k$.

Thus the pairings $\La_k\times\La_{n-k}\to\La_n$ have trivial left kernels for 
all $k=0,\ldots,n$. Then $\dim\La_k\le\dim\La_{n-k}$, and since this 
inequality holds also with $k$ replaced by $n-k$, we have in fact an equality. 
Thus all the above pairings are nondegenerate, which is a necessary and 
sufficient condition for $\La$ to be a Frobenius algebra.

By Theorem 3.1 $h_{\La}(t)=h_{\bbS}(-t)^{-1}=h_{\Sd}(t)$ where $\bbS=\bbS(V,R)$. 
Thus the two graded algebras $\La$ and $\Sd$ have homogeneous components of 
equal dimensions. In particular, $\dim\Sd_n=1$ and $\Sd_{n+1}=0$. But 
$\Sd=\La(V^*\!,R^*)$ by Lemma 3.2, and the Hecke symmetry $R^*$ 
satisfies the trivial source condition. We deduce that $\Sd$ is a Frobenius 
algebra by the already established part of Theorem 4.5. Since $\bbS$ is Koszul, 
its Gorensteinness follows then from \cite{Pol-P, Remark 2 on p. 25}.
\endproof

\section
5. Auxiliary results for the tensor product of two Hecke algebras

This section collects several results needed to deal with the algebras $A(R',R)$ 
and $E(R',R)$ in the next section. One of our goals is to investigate exactness 
of the complexes $K_\bullet\bigl(\calM;(U_i)\bigr)$ for certain collections of 
subspaces in a module $\calM$ over the tensor product $\calH'_n\ot\calH_n$ of 
two Hecke algebras $\calH_n=\calH_n(q)$ and $\calH'_n=\calH_n(q^{-1})$. Here 
$q\ne0$.

We identify $\calH_n$ and $\calH'_n$ with their canonical images in 
$\calH'_n\ot\calH_n$. Denote by $T_1,\ldots,T_{n-1}$ the standard generators 
of $\calH_n$ and by $T'_1,\ldots,T'_{n-1}$ those of $\calH'_n$. Put 
$\calT_i=T_iT'_i$,\/ i.e.
$$
\calT_i=T'_i\ot T_i\in\calH'_n\ot\calH_n\,,
$$
for $i=1,\ldots,n-1$. Since the elements of $\calH_n$ commute with those 
of $\calH'_n$, the elements $\calT_1,\ldots,\calT_{n-1}$ satisfy the braid 
relations. However, in general only the cubic relations 
$(\calT_i-1)(\calT_i+q)(\calT_i+q^{-1})=0$ hold rather than the quadratic ones.

\proclaim
Proposition 5.1.
Suppose that $\calM=M'\ot M$ where $M$ is an $\calH_n$-module and $M'$ is an 
$\calH'_n$-module such that all indecomposable direct summands of $M$ and 
$M'$ have $1$-dimensional sources. Let $U_i=(\calT_i-1)\calM$ for each 
$i=1,\ldots,n-1$. Then the complex $K_\bullet\bigl(\calM;(U_i)\bigr)$ is exact.
\endproclaim

Since direct sum decompositions of $M$ and $M'$ give rise to a direct sum 
decomposition of $K_\bullet\bigl(\calM;(U_i)\bigr)$, in proving the exactness 
of that complex we need only to consider the case when
$$
M=\calH_n\ot_{\calH_\la}\bbk(\chi),\qquad 
M'=\calH'_n\ot_{\calH'_\mu}\bbk(\chi')
$$
where $\,\chi:\calH_\la\to\bbk$ and $\,\chi':\calH'_\mu\to\bbk$ are 
1-dimensional representations of parabolic subalgebras $\calH_\la\,$, 
$\calH'_\mu$ of the respective Hecke algebras. By abuse of notation we will 
use the same letter $\chi$ for the 1-dimensional representation 
$\calH'_\mu\ot\calH_\la\to\bbk$ which restricts to the given representations 
of $\calH_\la$ and $\calH'_\mu$. In particular, $\chi(x)=\chi'(x)$ for 
$x\in\calH'_\mu$. The $\calH'_n\ot\calH_n$-module
$$
\calM=M'\ot M\cong(\calH'_n\ot\calH_n)\ot_{(\calH'_\mu\ot\calH_\la)}\bbk(\chi)
$$
has a generator $c$ such that $T_ic=\chi(T_i)c$ for all $i$ with 
$\tau_i\in\frB_\la$ and $T'_ic=\chi(T'_i)c$ for all $i$ with 
$\tau_i\in\frB_\mu$. Consider the standard bases $\{T_\si\mid\si\in\frS_n\}$, 
$\,\{T'_\si\mid\si\in\frS_n\}$ for $\calH_n$ and $\calH'_n$. Then $\calM$ has 
a vector space basis
$$
\{T_\pi T'_\si c\mid\pi\in\calD_\la,\ \si\in\calD_\mu\}\qquad
\bigl(\calD_\la=\calD(\frS_n/\frS_\la),\ 
\calD_\mu=\calD(\frS_n/\frS_\mu)\bigr).
$$
In $\calM$ we obtain a filtration of vector subspaces 
$0=F_{-1}\calM\sbs F_0\calM\sbs F_1\calM\sbs\ldots$ taking $F_p\calM$ to be 
the linear span of the elements
$$
\{T_\pi T'_\si c\mid\pi\in\calD_\la,\ \si\in\calD_\mu,
\ \ell(\pi)+\ell(\si)\le p\}.
$$
Our strategy is to relate the question we study for $\calM$ to a similar question 
for the associated graded vector space $\grFM$ (cf. \cite{Pol-P, Ch. 1, 
Cor. 7.3}):

\proclaim
Lemma 5.2.
Put $\Up_i=\bigcap_{j<i}U_j$ and $\Si_i=\sum_{j>i}U_j$. Assume that
$$
\gr^F\Up_i=\bigcap_{j<i}\gr^FU_j,\quad
\gr^F\Si_i=\sum_{j>i}\gr^FU_j,\quad
\gr^F(\Up_i\cap\Si_i)=\gr^F\Up_i\cap\gr^F\Si_i
$$
for all $i=0,\ldots,n$. If the complex $K_\bullet\bigl(\grFM;(\gr^FU_i)\bigr)$ 
is exact, then so too is the complex $K_\bullet\bigl(\calM;(U_i)\bigr)$.
\endproclaim

\Proof.
All subspaces of $\calM$ are endowed with the induced filtrations. 
Under the assumptions stated there is an exact sequence of complexes
$$
0\to K_\bullet\bigl(F_{p-1}\calM;(F_{p-1}U_i)\bigr)\to
K_\bullet\bigl(F_p\calM;(F_pU_i)\bigr)\to
K_\bullet\bigl(\gr^F_p\!\calM;(\gr^F_pU_i)\bigr)\to0
$$
where the last complex is exact, for each $p\ge0$. Induction on $p$ shows that 
the complex $K_\bullet\bigl(F_p\calM;(F_pU_i)\bigr)$ is exact. But 
$F_p\calM=\calM$, $F_pU_i=U_i$ for large $p$.
\endproof

With $M$ and $M'$ assumed to be fixed, all conditions needed for an 
application of Lemma 5.2 will be verified in Lemmas 5.3--5.8. This will 
accomplish a proof of Proposition 5.1.

In Lemma 5.3 the subspaces $\gr^FU_i$ of $\grFM$ will be determined explicitly. 
A module structure over the 0-Hecke algebra $\calH_n(0)$ will be constructed on 
$\grFM$ in Lemma 5.4. It will enable us to derive exactness of the complex 
$K_\bullet\bigl(\grFM;(\gr^FU_i)\bigr)$ from the results of section 2. 
Comparison of the subspaces $\Up_i$, $\Si_i$ in $\calM$ with their counterparts 
in $\grFM$ will be provided by Lemmas 5.7, 5.8.

Denote by $x_\calM$ the linear operator by which an element 
$x\in\calH'_n\ot\calH_n$ acts on $\calM$. The generators $T_i$, $T'_i$ act on 
$\calM$ by the formulas similar to those for the action of $T_i$ on $M$. Recall 
the subsets $A_i(\la)$, $B_i(\la)$ of $\calD_\la$ defined in section 1. There 
are similar subsets of $\calD_\mu$. We have
$$
T_iT_\pi T'_\si c=T_{\tau_i\pi}T'_\si c\quad\hbox{if $\pi\in A_i(\la)$},\qquad
T'_iT_\pi T'_\si c=T_\pi T'_{\tau_i\si}c\quad\hbox{if $\si\in A_i(\mu)$}.
$$
If $\pi\in B_i(\la)$, then $T_\pi T'_\si c$ is an eigenvector for the operator 
$(T_i)_\calM$ with the eigenvalue $\chi_\pi(T_i)=\chi(T_{\pi^{-1}(i)})$. 
If $\si\in B_i(\mu)$, then $T_\pi T'_\si c$ is an eigenvector for $(T'_i)_\calM$ 
with the eigenvalue $\chi_\si(T'_i)=\chi(T'_{\si^{-1}(i)})$.

\smallskip
Each homogeneous component $\gr^F_p\!\calM$ of $\grFM$ has a basis
$$
\{v_{\pi,\mskip1mu \si}\mid\pi\in\calD_\la,\ 
\si\in\calD_\mu,\ \ell(\pi)+\ell(\si)=p\}
$$
where $v_{\pi,\mskip1mu \si}=T_\pi T'_\si c+F_{p-1}\calM$. Hence 
$\{v_{\pi,\mskip1mu \si}\mid\pi\in\calD_\la,\ \si\in\calD_\mu\}$ is a basis 
for $\grFM$.

\proclaim
Lemma 5.3.
For each $i=1,\ldots,n-1$ the space $\gr^FU_i$ has a basis consisting of the 
following elements\/{\rm:}
$$
\eqalign{
v_{\tau_i\pi,\tau_i\si}\qquad&{\rm with}\ \pi\in A_i(\la),\ \si\in A_i(\mu),\cr
v_{\tau_i\pi,\mskip1mu \si}-qv_{\pi,\tau_i\si}\qquad&
{\rm with}\ \pi\in A_i(\la),\ \si\in A_i(\mu),\cr
v_{\tau_i\pi,\mskip1mu \si}\qquad&{\rm with}\ \pi\in A_i(\la),\ \si\in B_i(\mu),\cr
v_{\pi,\tau_i\si}\qquad&{\rm with}\ \pi\in B_i(\la),\ \si\in A_i(\mu),\cr
v_{\pi,\mskip1mu \si}\qquad&{\rm with}\ \pi\in B_i(\la),\ \si\in B_i(\mu),\ 
\chi_{\pi}(T_i)\chi_{\si}(T'_i)\ne1.\cr
}
$$
\endproclaim

\Proof.
Fixing $i$, let $\calM(\pi,\si)$ be the $\{T_i,T'_i\}$-invariant subspace 
of $\calM$ spanned by
$$
T_\pi T'_\si c\,,\qquad
T_{\tau_i\pi}T'_\si c\,,\qquad
T_\pi T'_{\tau_i\si}c\,,\quad
T_{\tau_i\pi}T'_{\tau_i\si}c\,.
$$
If $\pi\in A_i(\la)$ and $\si\in A_i(\mu)$, then these 4 elements are 
linearly independent. A basis for $\calM(\pi,\si)$ is formed by 2 elements 
$T_\pi T'_\si c$, $T_{\tau_i\pi}T'_\si c$ when $\pi\in A_i(\la)$, 
$\si\in B_i(\mu)$, and by 2 elements $T_\pi T'_\si c$, $T_\pi T'_{\tau_i\si}c$ 
when $\pi\in B_i(\la)$, $\si\in A_i(\mu)$. Computing the action of $\calT_i$, 
we deduce that $U_i\cap\calM(\pi,\si)$ is spanned in the first case by
$$
T_{\tau_i\pi}T'_{\tau_i\si}c-T_\pi T'_\si c\quad{\rm and}\quad
T_{\tau_i\pi}T'_\si c-qT_\pi T'_{\tau_i\si}c+(1-q)T_\pi T'_\si c\,,
$$
in the second by $\chi_\si(T'_i)T_{\tau_i\pi}T'_\si c-T_\pi T'_\si c$, in the 
third by $\chi_\pi(T_i)T_\pi T'_{\tau_i\si}c-T_\pi T'_\si c$. 
If $\pi\in B_i(\la)$ and $\si\in B_i(\mu)$, then $\calM(\pi,\si)$ is spanned 
by a single element $T_\pi T'_\si c$ which is an eigenvector for the operator
by which $\calT_i$ acts on $\calM$ with the eigenvalue 
$\chi_\pi(T_i)\chi_\si(T'_i)$. In this case $\calM(\pi,\si)\sbs U_i$ if and 
only if $\chi_ \pi(T_i)\chi_\si(T'_i)\ne1$.

In each case there is a basis for $\gr^F(U_i\cap\calM(\pi,\si))$ given by the 
respective elements in the statement of Lemma 5.3.

Note that $\calM$ is a direct sum of these subspaces $\calM(\pi,\si)$ with 
$\pi\in A_i(\la)\cup B_i(\la)$ and $\si\in A_i(\mu)\cup B_i(\mu)$. This direct 
sum decomposition is compatible with the filtration of $\calM$. Furthermore, 
we have $U_i=\bigoplus(U_i\cap\calM(\pi,\si))$ since each $\calM(\pi,\si)$ is 
stable under the action of $\calT_i$, whence 
$\gr^FU_i=\bigoplus\gr^F(U_i\cap\calM(\pi,\si))$.
\endproof

With the next goal to describe an $\calH_n(0)$-module structure it will be more 
convenient to index the basis elements of $\grFM$ by the pairs of cosets 
since this will allow us to exploit the natural actions of $\frS_n$ on 
$\frS_n/\frS_\la$ and on $\frS_n/\frS_\mu$.

For $x\in\frS_n/\frS_\la$ and $y\in\frS_n/\frS_\mu$ with their distinguished 
representatives $\pi\in\calD_\la$ and $\si\in\calD_\mu$ we put 
$\vxy=v_{\pi,\mskip1mu \si}$, and we will write $\chi_x(T_i)$, 
$\chi_y(T'_i)$, instead of $\chi_\pi(T_i)$, $\chi_\si(T'_i)$. Note that 
$\pi\in B_i(\la)$ if and only if $\tau_ix=x$. Hence $\chi_x(T_i)$ is defined 
if $\tau_ix=x$. Similarly, $\chi_y(T'_i)$ is defined if $\tau_iy=y$.

Consider the partial orders on $\frS_n/\frS_\la$ and $\frS_n/\frS_\mu$ 
transferred from the Bruhat orders on $\calD_\la$ and $\calD_\mu$. For $x$ and 
$\pi$ as in the preceding paragraph, we have $\tau_ix>x$ if and only if 
$\pi\in A_i(\la)$. Similarly, $\tau_iy>y$ if and only if $\si\in A_i(\mu)$. 
The next lemma applies to $N=\grFM$.

\proclaim
Lemma 5.4.
Let $N$ be a vector space with a basis $\{\vxy\mid 
x\in\frS_n/\frS_\la,\ y\in\frS_n/\frS_\mu\}$. Define linear operators 
$\oT_1,\ldots,\oT_{n-1}$ on $N$ by the rule
$$
\oT_i\vxy=\cases{
v_{\tau_ix,\mskip1mu\tau_iy} & \rm if $\tau_ix>x$, $\tau_iy\ge y$ or $\tau_ix=x$, 
$\tau_iy>y$,\cr
\noalign{\vskip3pt}
0 & \rm if $\tau_ix>x$, $\tau_iy<y$\cr
\noalign{\vskip3pt}
-\vxy+qv_{\tau_ix,\mskip1mu\tau_iy} & \rm if $\tau_ix<x$, $\tau_iy>y$,\cr
\noalign{\vskip3pt}
-\vxy & \rm if $\tau_ix<x$, $\tau_iy\le y$ or $\tau_ix=x$, 
$\tau_iy<y$,\cr
\noalign{\vskip3pt}
-\vxy & \rm if $\tau_ix=x$, $\tau_iy=y$, $\chi_x(T_i)\chi_y(T'_i)\ne1$,\cr
\noalign{\vskip3pt}
0 & \rm if $\tau_ix=x$, $\tau_iy=y$, $\chi_x(T_i)\chi_y(T'_i)=1$.\cr
}
$$
Then $\oT_1,\ldots,\oT_{n-1}$ satisfy the defining relations of the $0$-Hecke 
algebra $\calH_n(0)$.
\endproclaim

\Proof.
It is checked immediately that $\oT_i^2=-\oT_i$. Also, we have to show that
$$
\oT_i\oT_j\oT_i=\oT_j\oT_i\oT_j\quad\hbox{when }|i-j|=1,\qquad
\oT_i\oT_j=\oT_j\oT_i\quad\hbox{when }|i-j|>1.
$$
With the aim to do this express $\oT_i=\Phi_i+\Psi_i$ as the sum of two linear 
operators defined by the formulas
$$
\openup3\jot
\eqalign{ 
\Phi_i\vxy&{}=\cases{-\vxy & if $\tau_ix<x$ 
or $\tau_ix=x$, $\tau_iy<y$,\cr
\noalign{\vskip3pt}
-\vxy & if $\tau_ix=x$, $\tau_iy=y$, $\chi_x(T_i)\chi_y(T'_i)\ne1$,\cr
\noalign{\vskip3pt}
0 & in the remaining cases,\cr
}\cr
\Psi_i\vxy&{}=\cases{
v_{\tau_ix,\mskip1mu\tau_iy} & \rm if $\tau_ix>x$, $\tau_iy\ge y$ or $\tau_ix=x$, 
$\tau_iy>y$,\cr
\noalign{\vskip3pt}
qv_{\tau_ix,\mskip1mu\tau_iy} & \rm if $\tau_ix<x$, $\tau_iy>y$,\cr
\noalign{\vskip3pt}
0 & in the remaining cases.\cr
}
}
$$
In particular, each $\vxy$ is an eigenvector for $\Phi_i$. Note 
also that $\Psi_i\vxy$ is always a scalar multiple of 
$v_{\tau_ix,\mskip1mu\tau_iy}$, and $\Psi_i\vxy\ne0$ if and only 
if either $\tau_iy>y$ or $\tau_iy=y$ and $\tau_ix>x$.

Suppose first that $|i-j|=1$. We claim that
$$
\openup1\jot
\displaylines{
\Phi_i\Phi_j\Phi_i=\Phi_j\Phi_i\Phi_j\,,\qquad
\Psi_i\Phi_j\Psi_i=0\,,\qquad
\Psi_i\Psi_j\Psi_i=\Psi_j\Psi_i\Psi_j\,,\cr
\Psi_i\Psi_j\Phi_i=\Phi_j\Psi_i\Psi_j\,,\qquad
\Phi_i\Psi_j\Phi_i=\Psi_j\Phi_i\Phi_j+\Phi_j\Phi_i\Psi_j.
}
$$

Note that $\Phi_i\Phi_j\Phi_i\vxy\ne0$ if and only if both 
$\Phi_i\vxy\ne0$ and $\Phi_j\vxy\ne0$. If these inequalities hold, then 
$\Phi_i\Phi_j\Phi_i\vxy=-\vxy$. Since this description is symmetric in 
$i$ and $j$, we get \ $\Phi_i\Phi_j\Phi_i\vxy=\Phi_j\Phi_i\Phi_j\vxy$ 
\ in all cases.

Next, $\Phi_j\Psi_i\vxy$ is either 0 or equal to $-\Psi_i\vxy$. From the 
definition of $\Psi_i$ it is clear that $\Psi_i^2=0$. Hence 
$\Psi_i\Phi_j\Psi_i\vxy=0$.

\medskip
We have checked the first two identities. Before we proceed with the others 
let us make several remarks. The two transpositions $\tau_i$, $\tau_j$ generate 
a subgroup $\<\tau_i,\tau_j\>$ of $\frS_n$ isomorphic to $\frS_2$. Each 
$\<\tau_i,\tau_j\>$-orbit in $\frS_n/\frS_\la$ has a smallest element. In fact, 
a coset $x\in\frS_n/\frS_\la$ is minimal in its $\<\tau_i,\tau_j\>$-orbit if 
and only if $\tau_ix\ge x$ and $\tau_jx\ge x$, if and only if the distinguished 
representative $\pi_x$ of $x$ lies in 
$\calD(\<\tau_i,\tau_j\>\backslash\frS_n/\frS_\la)$.

If $x$ is minimal in its orbit, then the stabilizer 
$\<\tau_i,\tau_j\>\cap\pi_x\frS_\la\pi_x^{-1}$ of $x$ in $\<\tau_i,\tau_j\>$ 
is a parabolic subgroup by the general properties of Coxeter groups. Hence 
there are exactly 4 possibilities for this stabilizer: the trivial subgroup, 
the subgroups $\<\tau_i\>$, $\<\tau_j\>$ generated by one of the two 
transpositions, and the whole $\<\tau_i,\tau_j\>$. In the first case 
$\<\tau_i,\tau_j\>x$ is isomorphic as a poset to $\frS_2$ with the Bruhat 
order. In the second case the orbit contains 3 elements forming a chain 
$x<\tau_jx<\tau_i\tau_jx$. The third case has a similar description with $i$ 
and $j$ interchanged. In the last case $\<\tau_i,\tau_j\>x$ is the 
single element set $\{x\}$.

The action of the longest element $w=\tau_i\tau_j\tau_i=\tau_j\tau_i\tau_j$ of 
the group $\<\tau_i,\tau_j\>$ reverses order on each $\<\tau_i,\tau_j\>$-orbit 
in $\frS_n/\frS_\la$. From this it is clear that each 
$\<\tau_i,\tau_j\>$-orbit has a largest element, and $x\in\frS_n/\frS_\la$ is 
maximal in its orbit if and only if $\tau_ix\le x$ and $\tau_jx\le x$.

Since $w\tau_i=\tau_i\tau_j=\tau_jw$, we have $\tau_ix=x$ if and only if $wx$ 
is fixed by $\tau_j$. We claim that in this case 
$\chi_x(T_i)=\chi_{wx}(T_j)=\chi_{\tau_i\tau_jx}(T_j)$. 
If $\tau_jx\ne x$, then $x$ is either the smallest or the largest element in 
its $\<\tau_i,\tau_j\>$-orbit, and $wx=\tau_i\tau_jx$ is, respectively, the 
largest or the smallest element in this orbit with the distinguished 
representative $\pi_{wx}=\tau_i\tau_j\pi_x\in\calD_\la$ in both cases, so that
$$ 
\chi_{wx}(T_j)=\chi(T_{\pi_{wx}^{-1}\tau_j\pi_{wx}})
=\chi(T_{\pi_x^{-1}\tau_i\pi_x})=\chi_x(T_i).
$$
Suppose now that $\tau_jx=x$. Then both $\pi_x^{-1}\tau_i\pi_x$ and 
$\pi_x^{-1}\tau_j\pi_x$ are in $\frB_\la$. If these two transpositions are 
$\tau_k$ and $\tau_l$, then $\tau_k\tau_l\tau_k=\tau_l\tau_k\tau_l$ since 
$\tau_i\tau_j\tau_i=\tau_j\tau_i\tau_j$. This means that $|k-l|=1$, and so 
$\chi(T_k)=\chi(T_l)$. Hence
$$
\chi_{wx}(T_j)=\chi_x(T_j)=\chi(T_l)=\chi(T_k)=\chi_x(T_i).
$$

All the previous observations apply also to the $\<\tau_i,\tau_j\>$-orbits in 
$\frS_n/\frS_\mu$.

\medskip
Let us now turn to the third identity $\Psi_i\Psi_j\Psi_i=\Psi_j\Psi_i\Psi_j$. 
\ Consider
$$
\openup1\jot
\vcenter{\halign{\hfil$\displaystyle#$&$\displaystyle#$\hfil
&\quad\hfil$\displaystyle#$&$\displaystyle#$\hfil
&\quad\hfil$\displaystyle#$&$\displaystyle#$\hfil
&\quad\hfil$\displaystyle#$&$\displaystyle#$\hfil \cr
x_0&{}=x, & x_1&{}=\tau_ix, & x_2&{}=\tau_jx_1=\tau_j\tau_ix, 
& x_3&{}=\tau_ix_2=wx,\cr
y_0&{}=y, & y_1&{}=\tau_iy, & y_2&{}=\tau_jy_1=\tau_j\tau_iy, 
& y_3&{}=\tau_iy_2=wy,\cr
x^*_0&{}=x, & x^*_1&{}=\tau_jx, & x^*_2&{}=\tau_ix^*_1=\tau_i\tau_jx, 
& x^*_3&{}=\tau_jx^*_2=wx,\cr
y^*_0&{}=y, & y^*_1&{}=\tau_jy, & y^*_2&{}=\tau_iy^*_1=\tau_i\tau_jy, 
& y^*_3&{}=\tau_jy^*_2=wy.\cr
}}
$$
where $w=\tau_i\tau_j\tau_i$.
Note that $\Psi_i\Psi_j\Psi_i\vxy\ne0$ if and only if for each $k=1,2,3$ we 
have $y_k\ge y_{k-1}$, and if $y_k=y_{k-1}$ then $x_k>x_{k-1}$. Similarly, 
$\Psi_j\Psi_i\Psi_j\vxy\ne0$ if and only if for each $k=1,2,3$ we have 
$y_k^*\ge y^*_{k-1}$, and if $y^*_k=y^*_{k-1}$ then $x^*_k>x^*_{k-1}$.

Since $x^*_k=wx_{3-k}$, $y^*_k=wy_{3-k}$ and $w$ reverses order on 
$\<\tau_i,\tau_j\>$-orbits in $\frS_n/\frS_\la$ and $\frS_n/\frS_\mu$, we see 
that $\Psi_i\Psi_j\Psi_i\vxy\ne0$ if and only if 
$\Psi_j\Psi_i\Psi_j\vxy\ne0$. Furthermore, if these two elements are 
nonzero, then
$$
\Psi_i\Psi_j\Psi_i\vxy=q^av_{wx,\mskip1mu wy}\,,\qquad
\Psi_j\Psi_i\Psi_j\vxy=q^{a^*}v_{wx,\mskip1mu wy}
$$
where two numbers $a$, $a^*$ are the cardinalities of the sets of integers 
$1\le k\le 3$ such that $\,y_k>y_{k-1}$, $\,x_k<x_{k-1}\,$ in the case of $a$ 
and $\,y^*_k>y^*_{k-1}$, $\,x^*_k<x^*_{k-1}\,$ in the case of $a^*$. Since the 
assignment $k\mapsto4-k$ gives a bijection between these two sets, we have 
$a=a^*$, whence \ $\Psi_i\Psi_j\Psi_i\vxy=\Psi_j\Psi_i\Psi_j\vxy\,$.

\medskip
The fourth identity.
Recall that $\Psi_i\Psi_j\vxy$ is always a scalar multiple of 
$v_{x^*\!,\mskip1mu y^*}$ where $x^*=\tau_i\tau_jx$, $\,y^*=\tau_i\tau_jy$. For 
$\Psi_i\Psi_j\Phi_i\vxy\ne0$ to hold, it is necessary and sufficient that 
$\Psi_i\Psi_j\vxy\ne0$ and $\Phi_i\vxy\ne0$, while 
$\Phi_j\Psi_i\Psi_j\vxy\ne0$ holds if and only if $\Psi_i\Psi_j\vxy\ne0$ 
and $\Phi_jv_{x^*\!,y^*}\ne0$.

Since $\tau_i\tau_j\tau_i$ reverses the order on each $\<\tau_i,\tau_j\>$-orbit 
in $\frS_n/\frS_\la$, we have $\tau_ix<x$ if and only if $\tau_jx^*<x^*$. Also, 
$\tau_ix=x$ if and only if $\tau_jx^*=x^*$, and in this case 
$\chi_{x^*}(T_j)=\chi_x(T_i)$. Similarly, $\tau_iy<y$ if and only if 
$\tau_jy^*<y^*$, and $\tau_iy=y$ if and only if $\tau_jy^*=y^*$. If $\tau_iy=y$, 
then $\chi_{y^*}(T'_j)=\chi_y(T'_i)$.

Hence $\chi_{x^*}(T_j)\chi_{y^*}(T'_j)=\chi_x(T_i)\chi_y(T'_i)$ whenever 
$\tau_ix=x$ and $\tau_iy=y$ hold simultaneously. So it follows that 
$\Phi_i\vxy\ne0$ if and only if $\Phi_jv_{x^*\!,\mskip1mu y^*}\ne0$, in which case
$$
\Psi_i\Psi_j\Phi_i\vxy=-\Psi_i\Psi_j\vxy=\Phi_j\Psi_i\Psi_j\vxy\,.
$$

The fifth identity.
All three terms $\Phi_i\Psi_j\Phi_i\vxy$, $\Psi_j\Phi_i\Phi_j\vxy$, 
$\Phi_j\Phi_i\Psi_j\vxy$ equal 0 when $\Psi_j\vxy=0$. Suppose that 
$\Psi_j\vxy\ne0$. Then $\tau_jy\ge y$, and if $\tau_jy=y$ then $\tau_jx>x$. 
Hence $\tau_jy>y$ whenever $\tau_jx\le x$.

If $\tau_jx\ge x$, we must have $\Phi_j\vxy=0$ and 
$\Phi_jv_{\tau_jx,\tau_jy}\ne0$ by the definition of $\Phi_j$. In this case 
$\Psi_j\Phi_i\Phi_j\vxy=0$, while 
$\Phi_j\Phi_i\Psi_j\vxy=-\Phi_i\Psi_j\vxy$. Furthermore, the inequality 
$\Phi_i\Psi_j\vxy\ne0$ is only possible when either $\tau_i\tau_jx<\tau_jx$ 
or $\tau_i\tau_jx=\tau_jx$ and $\tau_i\tau_jy\le\tau_jy$. These conditions 
imply that $\tau_jx$ is the largest element of the $\<\tau_i,\tau_j\>$-orbit 
of $x$, and therefore $\tau_ix<x$ when $x<\tau_jx$. If $\tau_jx=x$, we have 
either $\tau_ix<x$ or $\tau_ix=x$ and $\tau_i\tau_jy\le\tau_jy$. In 
particular, $\tau_jy$ is the largest element of the $\<\tau_i,\tau_j\>$-orbit 
of $y$ in the case when $\tau_ix=\tau_jx=x$, but then $\tau_iy<y$ since 
$y<\tau_jy$. In each of these cases we deduce that 
$\Phi_i\vxy\ne0$, whence
$$
\Phi_i\Psi_j\Phi_i\vxy=-\Phi_i\Psi_j\vxy=\Phi_j\Phi_i\Psi_j\vxy
=\Psi_j\Phi_i\Phi_j\vxy+\Phi_j\Phi_i\Psi_j\vxy\,.
$$
The same equalities are trivially true when $\tau_jx\ge x$, but 
$\Phi_i\Psi_j\vxy=0$.

Consider the remaining case when $\tau_jx<x$. Here $\Phi_j\vxy\ne0$ and 
$\Phi_jv_{\tau_jx,\tau_jy}=\nobreak0$, so that $\Phi_j\Phi_i\Psi_j\vxy=0$, 
while $\Psi_j\Phi_i\Phi_j\vxy=-\Psi_j\Phi_i\vxy$. If now 
$\Phi_i\vxy=0$, then all the three terms we look at vanish. If 
$\Phi_i\vxy\ne0$, then $\tau_ix\le x$, which implies that $x$ is the 
largest element in its $\<\tau_i,\tau_j\>$-orbit. But this entails 
$\tau_i\tau_jx<\tau_jx$ since $\tau_jx<x$, and it follows that 
$\Phi_iv_{\tau_jx,\tau_jy}\ne0$. Hence
$$
\Phi_i\Psi_j\Phi_i\vxy=-\Psi_j\Phi_i\vxy=\Psi_j\Phi_i\Phi_j\vxy
=\Psi_j\Phi_i\Phi_j\vxy+\Phi_j\Phi_i\Psi_j\vxy\,.
$$

We have checked all the required relations between $\Phi_i,\Phi_j,\Psi_i,\Psi_j$ 
for any pair $i,j$ with $|i-j|=1$. By symmetry they hold also with $i$ and $j$ 
interchanged. The braid relation $\oT_i\oT_j\oT_i=\oT_j\oT_i\oT_j$ is now 
immediate.

Suppose that $|i-j|>1$. Considering the subgroup $\<\tau_i,\tau_j\>$ generated 
by $\tau_i,\tau_j$, it is still true that its longest element 
$\tau_i\tau_j=\tau_j\tau_i$ reverses order on each $\<\tau_i,\tau_j\>$-orbit in 
$\frS_n/\frS_\la$ and $\frS_n/\frS_\mu$. If $\tau_ix=x$ for some 
$x\in\frS_n/\frS_\la$, then $\chi_{\tau_jx}(T_i)=\chi_x(T_i)$, and if 
$\tau_iy=y$ for some $y\in\frS_n/\frS_\mu$, then 
$\chi_{\tau_jy}(T'_i)=\chi_y(T'_i)$. The arguments similar to those used in 
the case $|i-j|=1$, but this time much shorter, show that
$$
\Phi_i\Phi_j=\Phi_j\Phi_i\,,\qquad
\Psi_i\Psi_j=\Psi_j\Psi_i\,,\qquad
\Phi_i\Psi_j=\Psi_j\Phi_i.
$$
Hence \ $\oT_i\oT_j=(\Phi_i+\Psi_i)(\Phi_j+\Psi_j)
=(\Phi_j+\Psi_j)(\Phi_i+\Psi_i)=\oT_j\oT_i\,$.
\endproof

\proclaim
Lemma 5.5.
Put $N=\grFM$ and denote by $\oT_1,\ldots,\oT_{n-1}$ the canonical generators 
of the $0$-Hecke algebra $\calH_n(0)$. With the $\calH_n(0)$-module structure 
on $N$ defined by the formulas in the statement of Lemma\/ {\rm 5.4} we have 
$\gr^FU_i=\oT_iN$ for each $i,$ and there is an isomorphism of 
$\calH_n(0)$-modules
$$
N\cong\bigoplus_{\pi\in\,\muDlass}\calH_n(0)\ot_{\calH_{\nu(\pi)}(0)}\bbk(\xi_\pi)
$$
where $\nu(\pi)$ is the composition of $n$ such that 
$\frS_{\nu(\pi)}=\frS_\mu\cap\pi\frS_\la\pi^{-1},$ $\,\calH_{\nu(\pi)}(0)$ the 
corresponding parabolic subalgebra of $\calH_n(0),$ and 
$\xi_\pi:\calH_{\nu(\pi)}(0)\to\bbk$ the representation such that
$$
\xi_\pi(\oT_i)=\cases{-1 & \rm if $\tau_i\in\frB_{\nu(\pi)}$ and 
$\chi_\pi(T_i)\chi(T'_i)\ne 1$,\cr
\noalign{\smallskip}
0 & \rm if $\tau_i\in\frB_{\nu(\pi)}$ and $\chi_\pi(T_i)\chi(T'_i)=1$.\cr
}
$$
\endproclaim

\Proof.
In terms of the indexation of the basis elements of $N$ by the pairs 
$(x,y)\in\frS_n/\frS_\la\times\frS_n/\frS_\mu$ the image of the linear 
operator $\oT_i:N\to N$ is spanned by the elements $\vxy$ with
$$
\openup1\jot
\eqalign{
&\tau_ix<x,\ \tau_iy\le y,\quad{\rm or}\ \,\tau_ix=x,\ \tau_iy<y,\cr
&{\rm or}\ \,\tau_ix=x,\ \tau_iy=y,\ \chi_x(T_i)\chi_y(T'_i)\ne1,\cr
}
$$
and by the elements $\vxy-qv_{\tau_ix,\tau_iy}$ with $\tau_ix<x$, $\tau_iy>y$.

If $\pi\in\calD_\la$ is the distinguished representative of $x$, then 
$\tau_ix<x$ if and only if $\pi\in\tau_iA_i(\la)$, and $\tau_ix=x$ if and 
only if $\pi\in B_i(\la)$. Similarly, if $\si\in\calD_\mu$ is the 
distinguished representative of $y$, then $\tau_iy>y$, $\tau_iy<y$, or 
$\tau_iy=y$ depending on whether $\si$ is in $A_i(\mu)$, $\tau_iA_i(\mu)$, 
or $B_i(\mu)$, respectively. Comparison with the description given in Lemma 
5.3 yields the desired equality $\,\gr^FU_i=\oT_iN$.

We claim that the $\calH_n(0)$-module $N$ is generated by the set 
$$
\{v_{x,\mskip1mu e_\mu}\mid\hbox{$x\in\frS_n/\frS_\la$, $\tau_ix\ge x$ for all 
$\tau_i\in\frB_\mu$}\}.
$$
where $e_\mu=\frS_\mu$ is the coset of the identity element. If $\tau_iy>y$ 
for some $y\in\frS_n/\frS_\mu$ and $\tau_i\in\frB_n$, then $v_{x,\tau_iy}$ 
equals $\oT_iv_{\tau_ix,y}$ when $\tau_ix\le x$ and equals 
$q^{-1}(1+\oT_i)v_{\tau_ix,y}$ when $\tau_ix>x$. From this it follows by 
induction on $y$ that each basis element $\vxy$ lies in the submodule of $N$ 
generated by $\{v_{x'\!,\mskip1mu e_\mu}\mid x'\in\frS_n/\frS_\la\}$. On the 
other hand, $v_{\tau_ix'\!,\mskip1mu e_\mu}=\oT_iv_{x'\!,\mskip1mu e_\mu}$ if 
$\tau_ix'>x'$ for some $\tau_i\in\frB_\mu$. Induction on $x'$ shows that 
$v_{x'\!,\mskip1mu e_\mu}$ lies in the submodule of $N$ generated by 
$v_{x,\mskip1mu e_\mu}$ where $x$ is the smallest element of the 
$\frS_\mu$-orbit of $x'$. This proves the claim about the generating set of $N$.

Denote by $N(x)$ the submodule of $N$ generated by $v_{x,e_\mu}$. From the 
formulas in the statement of Lemma 5.4 it is clear that $N(x)$ is contained in 
the subspace of $N$ spanned by $\{v_{\si x,\mskip1mu\si e_\mu}\mid\si\in\frS_n\}$. 
If $\tau_ix\ge x$ for all $\tau_i\in\frB_\mu$, then the considerations in the 
preceding paragraph in fact show that $v_{\si x,\mskip1mu\si e_\mu}\in N(x)$ for 
all $\si\in\frS_n$. In this case $N(x)$ has a basis consisting of the elements 
$v_{x'\!,\mskip1mu y'}$ with $(x'\!,y')$ in the orbit of $(x,e_\mu)$ with 
respect to the diagonal action of $\frS_n$ on 
$\frS_n/\frS_\la\times\frS_n/\frS_\mu\,$, and therefore
$$
\dim N(x)=|\frS_n|/|\St(x,e_\mu)|=
\hbox{the index of $\St(x,e_\mu)$ in $\frS_n$}
$$
where $\St(x,e_\mu)$ stands for the stabilizer of $(x,e_\mu)$ in $\frS_n$ with 
respect to that action. The condition that $\tau_ix\ge x$ for all 
$\tau_i\in\frB_\mu$ means precisely that the distinguished representative of 
$x$ lies in $\muDla$. If $\pi$ is this representative, then 
$\St(x,e_\mu)=\frS_{\nu(\pi)}$.

Now we put $N(\pi)=N(\pi\frS_\la)$ for each $\pi\in\muDla$. In other words, 
$N(\pi)$ is the submodule of $N$ generated by $v_{\pi,\mskip1mu e}$. Then 
$\,\dim N(\pi)=(\frS_n:\frS_{\nu(\pi)})\,$ by the above. If 
$\tau_i\in\frB_{\nu(\pi)}$, then $\tau_i\pi\frS_\la=\pi\frS_\la$ and 
$\tau_i\frS_\mu=\frS_\mu$, whence
$$
\oT_iv_{\pi,e}=\xi_\pi(\oT_i)v_{\pi,e}
$$
by comparison of the definition of $\xi_\pi$ with the last two formulas in the 
statement of Lemma 5.4. Hence there is a surjective homomorphism of
$\calH_n(0)$-modules
$$
\calH_n(0)\ot_{\calH_{\nu(\pi)}(0)}\bbk(\xi_\pi)\mapr{}N(\pi).
$$
Comparing the dimensions, we deduce that this map is an isomorphism.

From the preceding discussion it is also clear that each basis element of $N$ 
lies in exactly one submodule $N(\pi)$. Thus 
$N=\bigoplus_{\pi\in\,\muDlass}N(\pi)$.
\endproof

\proclaim
Lemma 5.6.
The complex $K_\bullet\bigl(\grFM;(\gr^FU_i)\bigr)$ is exact.
\endproclaim

\Proof.
By Lemma 5.5 there is an $\calH_n(0)$-module structure on $N=\grFM$ with the 
property that $\gr^FU_i=\oT_iN$ for each $i$. Therefore Corollary 2.5 applies.
\endproof

\proclaim
Lemma 5.7.
Put $\,\Si=\sum\limits_{i=1}^{n-1}U_i,$ 
$\,\Si^{\gr}=\sum\limits_{i=1}^{n-1}\gr^FU_i,$ 
$\,\Up=\bigcap\limits_{i=1}^{n-1}U_i,$ 
$\Up^{\gr}=\bigcap\limits_{i=1}^{n-1}\gr^FU_i$.
We have $\,\gr^F\Si=\Si^{\gr}$ and $\,\gr^F\Up=\Up^{\gr}$.
\endproclaim

\Proof.
The inclusions $\,\Si^{\gr}\sbs\gr^F\Si\,$ and $\,\gr^F\Up\sbs\Up^{\gr}\,$ 
are always true, and so we need only to compare the dimensions.
By Lemma 5.5 $\,N=\grFM\,$ is a direct sum of $\calH_n(0)$-submodules 
$$
N(\pi)\cong\calH_n(0)\ot_{\calH_{\nu(\pi)}(0)}\bbk(\xi_\pi),\qquad\pi\in\muDla.
$$
Now $N/\Si^{\gr}=N/\sum\oT_iN$ is the largest factor module of $N$ annihilated 
by each $\oT_i$. Since
$$
\eqalign{
\Hom_{\calH_n(0)}(N,\bbk\triv)&{}\cong
\bigoplus_{\pi\in\,\muDlass}\Hom_{\calH_n(0)}\bigl(N(\pi),\bbk\triv\bigr)\cr
&{}\cong
\bigoplus_{\pi\in\,\muDlass}\Hom_{\calH_{\nu(\pi)}(0)}\bigl(\bbk(\xi_\pi),\bbk\triv\bigr),\cr
}
$$
$$
\eqalign{
\hbox{we get}\qquad
\dim N/\Si^{\gr}&{}=\dim\Hom_{\calH_n(0)}(N,\bbk\triv)\cr
&{}=\#\{\pi\in\muDla\mid\xi_\pi\hbox{ is the trivial 
representation}\}\cr
&{}=\#\{\pi\in\muDla\mid\chi_\pi(T_i)\chi(T'_i)=1
\hbox{ for all $i$ with }\tau_i\in\frS_{\nu(\pi)}\}.\cr
}
$$
As $(\oT_i+1)\oT_i=0$, we can also write $\gr^FU_i=\Ker\,(\oT_i+1)_N$. Hence 
$\Up^{\gr}$ is the largest submodule of $N$ on which each $\oT_i$ acts as $-\Id$. 
Since
$$
\Hom_{\calH_n(0)}(\bbk\alt,N)\cong
\bigoplus_{\pi\in\,\muDlass}\Hom_{\calH_n(0)}\bigl(\bbk\alt,N(\pi)\bigr),
$$
we deduce from Lemma 1.3 that
$$
\eqalign{
\dim\Up^{\gr}&{}=\dim\Hom_{\calH_n(0)}(\bbk\alt,N)\cr
&{}=\#\{\pi\in\muDla\mid\xi_\pi\hbox{ is the alternating representation}\}\cr
&{}=\#\{\pi\in\muDla\mid\chi_\pi(T_i)\chi(T'_i)\ne1
\hbox{ for all $i$ such that }\tau_i\in\frB_{\nu(\pi)}\}.\cr
}
$$

Next we are going to determine $\calM/\Si$ and $\Up$. For this we will need two 
different interpretations of the spaces $U_i$. The assignment $T_i\mapsto(T'_i)^{-1}$ 
extends to an algebra antiisomorphism $\calH_n\to\calH'_n$ under which 
$\calH_\mu$ is mapped onto $\calH'_\mu$. It allows us to view $M'$ as a right 
$\calH_n$-module and $\bbk(\chi')$ as a right $\calH_\mu$-module. Clearly,
$$
M'=\calH'_n\ot_{\calH'_\mu}\bbk(\chi')\cong\bbk(\chi')\ot_{\calH_\mu}\calH_n\,.
$$
The space $U_i$ is spanned by all elements 
$\calT_i(u\ot v)-u\ot v=u{T_i}^{-1}\ot T_iv-u\ot v$ with $u\in M'$ and $v\in M$. 
Replacing here $u$ with $uT_i$, we rewrite these elements as 
$u\ot T_iv-uT_i\ot v$. Since $T_1,\ldots,T_{n-1}$ generate $\calH_n$, the space 
$\Si$ is spanned by all elements $u\ot hv-uh\ot v$ with $u\in M'$, $v\in M$ and 
$h\in\calH_n$. It follows that
$$
\calM/\Si\cong M'\ot_{\calH_n}M\cong\bbk(\chi')\ot_{\calH_\mu}M\cong
\bigoplus_{\pi\in\,\muDlass}\bbk(\chi')\ot_{\calH_{\nu(\pi)}}\bbk(\chi_\pi)
$$
where the last isomorphism is a consequence of the Mackey decomposition formula 
since the restriction of $M$ to $\calH_\mu$ is a direct sum of modules 
$\calH_\mu\ot_{\calH_{\nu(\pi)}}\bbk(\chi_\pi)$ with $\pi\in\muDla$. 
Now $\bbk(\chi')\ot_{\calH_{\nu(\pi)}}\bbk(\chi_\pi)\ne0$ if and only if 
$\calH_{\nu(\pi)}$ operates in $\bbk(\chi')$ and in $\bbk(\chi_\pi)$ via 
the same homomorphism $\calH_{\nu(\pi)}\to\bbk$, i.e.  
$\chi(T'_i)^{-1}=\chi_\pi(T_i)$ for all $i$ such that 
$\tau_i\in\frB_{\nu(\pi)}$. We get
$$
\dim\calM/\Si=\#\{\pi\in\muDla\mid\chi_\pi(T_i)\chi(T'_i)=1
\hbox{ for all $i$ such that }\tau_i\in\frB_{\nu(\pi)}\}.
$$
Thus $\Si$ has the same codimension in $\calM$ as $\Si^{\gr}$ in $N$. Since the 
filtration on $\calM$ is exhaustive and separating, the dimension and the 
codimension of subspaces of $\calM$ are preserved under passage to the associated 
graded spaces. The first equality in the statement of Lemma 5.7 is now clear.

We may view $\calM=M'\ot M$ as an $\calH_n$-bimodule. By the previous description 
$U_i$ is spanned by all elements $T_i\psi-\psi T_i$ with $\psi\in\calM$. The 
dual vector space $\Md$ is a left $\calH_n$-module in a natural way. Identifying 
$\calM$ with $\Homk(\Md,M)$ by means of the canonical bijection, we see that 
$\Up$ consists precisely of those $\bbk$-linear maps $\Md\to M$ that satisfy 
condition (a) of Lemma 1.6.

Suppose first that $q\ne-1$. Then $\Up=\Hom_{\calH_n}(\Md,\Mt)$. 
By Lemmas 1.4, 1.5
$$
\Md\cong\calH_n\ot_{\calH_\mu}\bbk(\chi'),\qquad
\Mt\cong\calH_n\ot_{\calH_\la}\bbk(\chit).
$$
For each $i$ with $\tau_i\in\frB_\mu$ the generator $T_i$ of $\calH_\mu$ acts 
on $\bbk(\chi')$ as the multiplication by $\chi(T'_i)^{-1}$. 
Therefore Lemma 1.1 yields
$$
\dim\Up=\#\{\pi\in\muDla\mid\chi(T'_i)^{-1}=\chit_\pi(T_i)
\hbox{ for all $i$ such that }\tau_i\in\frB_{\nu(\pi)}\}.
$$
If $\tau_i\in\frB_{\nu(\pi)}$, then $\chit_\pi(T_i)=\chi_\pi(\Tt_i)$. Since 
$\chi_\pi(\Tt_i)+\chi_\pi(T_i)=q-1$ and $\chi_\pi(T_i)$ equals either $q$ or 
$-1$, we always have $\chi_\pi(\Tt_i)\ne\chi_\pi(T_i)$. Therefore 
$\chi(T'_i)^{-1}=\chi_\pi(\Tt_i)$ if and only if 
$\chi(T'_i)^{-1}\ne\chi_\pi(T_i)$. Thus
$$
\dim\Up=\#\{\pi\in\muDla\mid\chi_\pi(T_i)\chi(T'_i)\ne1
\hbox{ for all $i$ such that }\tau_i\in\frB_{\nu(\pi)}\}.
$$

If $q=-1$, then Lemma 1.7 with $M$ replaced by $\Mt$ shows that $\Up$ consists 
precisely of those $\calH_n$-module homomorphisms $\Md\to\Mt$ which factor 
through a free module. By the last assertion in that lemma
$$
\dim\Up=\#\{\pi\in\muDla\mid
\hbox{$\frS_{\nu(\pi)}$ is the trivial subgroup of $\frS_n$}\}.
$$
On the other hand, $\chi(T'_i)=\chi_\pi(T_i)=-1$, and therefore 
$\chi_\pi(T_i)\chi(T'_i)=1$, for each $i$ such that $\tau_i\in\frB_{\nu(\pi)}$. 
This means that in the case $q=-1$ the earlier formula for the dimension of 
$\Up^{\gr}$ counts only those permutations $\pi\in\muDla$ for which the 
set $\frB_{\nu(\pi)}$ is empty, i.e. $\frS_{\nu(\pi)}$ is the trivial group.

We conclude that $\dim\Up=\dim\Up^{\gr}$ both for $q\ne-1$ and for $q=-1$. 
This proves the second equality in the statement of Lemma 5.7. 
\endproof

\proclaim
Lemma 5.8.
Put\/ $\Up_i=\bigcap_{j<i}U_j$ and\/ $\Si_i=\sum_{j>i}U_j$. Then
$$
\gr^F\Up_i=\bigcap_{j<i}\gr^FU_j,\quad
\gr^F\Si_i=\sum_{j>i}\gr^FU_j,\quad
\gr^F(\Up_i\cap\Si_i)=\gr^F\Up_i\cap\gr^F\Si_i.
$$
\endproclaim

\Proof.
Let $\calH_{i,n-i}\sbs\calH_n$ and $\calH'_{i,n-i}\sbs\calH'_n$ be the parabolic 
subalgebras corresponding to the subgroup $\frS_{i,n-i}$ of $\frS_n$ generated 
by $\{\tau_j\in\frB_n\mid j\ne i\}$. Consider the Mackey decompositions
$$
M=\bigoplus_{\pi\in\calD(\frS_{i,n-i}\backslash\frS_n/\frS_\la)}M(\pi),\qquad
M'=\bigoplus_{\si\in\calD(\frS_{i,n-i}\backslash\frS_n/\frS_\mu)}M'(\si)
$$
where $M(\pi)$ is the $\calH_{i,n-i}$-submodule of $M$ generated by $T_\pi\ot1$ 
and $M'(\si)$ is the $\calH'_{i,n-i}$-submodule of $M'$ generated by 
$T'_\si\ot1$. They give rise to the decomposition of $\calM$ as a direct sum 
of $\,\calH'_{i,n-i}\ot\calH_{i,n-i}$-submodules
$$
\calM(\pi,\si)=M'(\si)\ot M(\pi).
$$
This decomposition is compatible with the filtration on $\calM$. Also, if 
$j\ne i$, then $U_j=\bigoplus\bigl(U_j\cap\calM(\pi,\si)\bigr)$ since each 
summand $\calM(\pi,\si)$ is stable under the action of $\calT_j$. It follows 
that $\Up_i=\bigoplus\bigl(\Up_i\cap\calM(\pi,\si)\bigr)$ and 
$\Si_i=\bigoplus\bigl(\Si_i\cap\calM(\pi,\si)\bigr)$. 

\smallskip
Next, $\calH_{i,n-i}\cong\calH_i\ot\calH^\tri_i$ where $\calH_i$ and 
$\calH^\tri_i$ are the subalgebras of $\calH_n$ generated, respectively, by 
$\{T_j\mid j<i\}$ and $\{T_j\mid j>i\}$. Since the $\calH_{i,n-i}$-module 
$M(\pi)$ is induced from a 1-dimensional module over a parabolic subalgebra 
of $\calH_{i,n-i}$, we have
$$
M(\pi)\cong M(\pi)_1\ot M(\pi)_2
$$
where $M(\pi)_1$ is an $\calH_i$-module induced from a 1-dimensional module 
over a parabolic subalgebra of $\calH_i$ and $M(\pi)_2$ is an 
$\calH^\tri_i$-module induced from a 1-dimensional module over a parabolic 
subalgebra of $\calH^\tri_i$. Let 
$\,\calH'_{i,n-i}\cong\calH'_i\ot\Hsptri\,$ and
$$
M'(\si)\cong M'(\si)_1\ot M'(\si)_2
$$
be similar decompositions. Then $\,\calM(\pi,\si)\cong 
\calM(\pi,\si)_1\ot\calM(\pi,\si)_2\,$ where
$$
\calM(\pi,\si)_1=M'(\si)_1\ot M(\pi)_1,\qquad
\calM(\pi,\si)_2=M'(\si)_2\ot M(\pi)_2.
$$
Note that $\calH'_i\ot\calH_i\cong\calH_i(q^{-1})\ot\calH_i(q)$ and 
$\Hsptri\ot\calH^\tri_i\cong\calH_{n-i}(q^{-1})\ot\calH_{n-i}(q)$.
The $\calH'_i\ot\calH_i$-module $\calM(\pi,\si)_1$ and the 
$\Hsptri\ot\calH^\tri_i$-module $\calM(\pi,\si)_2$ satisfy the same 
assumptions that we have imposed on the $\calH'_n\ot\calH_n$-module $\calM$.  
In particular, we obtain filtrations $F^1$, $F^2$ on these two modules by the 
construction we have done for $\calM$. Then the tensor product filtration 
$F^t$ on $\calM(\pi,\si)$ differs from the filtration induced from that on 
$\calM$ only by a shift of the filtration degrees.

Note that $\calT_j$ lies in $\calH'_i\ot\calH_i$ when $j<i$ and in 
$\Hsptri\ot\calH^\tri_i$ when $j>i$. Hence
$$
\openup1\jot
\eqalign{
U_j\cap\calM(\pi,\si)&{}=\cases{
U_j(\pi,\si)\ot\calM(\pi,\si)_2
& for $j<i$,\cr
\noalign{\smallskip}
\calM(\pi,\si)_1\ot U_j(\pi,\si)
& for $j>i$\cr
}\cr
{\rm where}\quad U_j(\pi,\si)&{}=\cases{
(\calT_j-1)\calM(\pi,\si)_1 & for $j<i$,\cr
\noalign{\smallskip}
(\calT_j-1)\calM(\pi,\si)_2 & for $j>i$.\cr
}\cr
} 
$$
With $\,\Up(\pi,\si)=\bigcap_{j<i}U_j(\pi,\si)\,$ and 
$\,\Si(\pi,\si)=\sum_{j>i}U_j(\pi,\si)\,$ we get
$$
\displaylines{
\Up_i\cap\calM(\pi,\si)=\Up(\pi,\si)\ot\calM(\pi,\si)_2,\qquad
\Si_i\cap\calM(\pi,\si)=\calM(\pi,\si)_1\ot\Si(\pi,\si)\cr
\noalign{\smallskip}
\Up_i\cap\Si_i\cap\calM(\pi,\si)=\Up(\pi,\si)\ot\Si(\pi,\si)\cr
}
$$
By Lemma 5.7 applied to $\calM(\pi,\si)_1$ and $\calM(\pi,\si)_2$ we have 
$$
\gr^{F^1}\!\Up(\pi,\si)=\bigcap_{j<i}\,\gr^{F^1}\!U_j(\pi,\si),\qquad
\gr^{F^2}\!\Si(\pi,\si)=\sum_{j>i}\,\gr^{F^2}\!U_j(\pi,\si),
$$
and it follows that
$$
\eqalign{
\gr^{F^t}\!\bigl(\Up_i\cap\calM(\pi,\si)\bigr)
&{}=\gr^{F^1}\!\Up(\pi,\si)\ot\gr^{F^2}\!\!\calM(\pi,\si)_2
=\bigcap_{j<i}\gr^{F^t}\!\bigl(U_j\cap\calM(\pi,\si)\bigr),\cr
\gr^{F^t}\!\bigl(\Si_i\cap\calM(\pi,\si)\bigr)
&{}=\gr^{F^1}\!\!\calM(\pi,\si)_1\ot\gr^{F^2}\!\Si(\pi,\si)
=\sum_{j>i}\gr^{F^t}\!\bigl(U_j\cap\calM(\pi,\si)\bigr),\cr
\gr^{F^t}\!\bigl(\Up_i\cap\Si_i\cap\calM(\pi,\si)\bigr)
&{}=\gr^{F^1}\!\Up(\pi,\si)\ot\gr^{F^2}\!\Si(\pi,\si)\cr
&{}=\gr^{F^t}\!\bigl(\Up_i\cap\calM(\pi,\si)\bigr)
\cap\gr^{F^t}\bigl(\Si_i\cap\calM(\pi,\si)\bigr).\cr
}
$$
The equalities of the left and right hand sides above hold then also with 
$F^t$ replaced by the original filtration $F$ on $\calM$ since this change 
results in the same associated graded spaces with shifted degrees of 
homogeneous components. Summing up over all pairs $(\pi,\si)\in 
\calD(\frS_{i,n-i}\backslash\frS_n/\frS_\la)\times 
\calD(\frS_{i,n-i}\backslash\frS_n/\frS_\mu)$, we arrive at the final 
conclusions.
\endproof

Now the proof of Proposition 5.1 is complete. There is also a version of this 
result for a different collection of subspaces in $\calM$:

\proclaim
Proposition 5.9.
With the same assumptions about $\calM$ as in Proposition\/ {\rm5.1} the 
complex $K_\bullet\bigl(\calM;(U_i)\bigr)$ is exact also in the case when 
$U_i=\Ker\,(\calT_i-1)_\calM$ for each $i=1,\ldots,n-1$.
\endproclaim

\Proof.
The dual space $\calM^*\cong\Md\ot M^*$ is a right $\calH'_n\ot\calH_n$-module 
in a natural way. We will view $M^*$ as a left $\calH_n$-module and $\Md$ as 
a left $\calH'_n$-module by means of the antiautomorphisms of $\calH_n$ and 
$\calH'_n$ such that $T_i\mapsto T_{n-i}$, $T'_i\mapsto T'_{n-i}$ for each 
$i=1,\ldots,n-1$. Since these antiautomorphisms map parabolic subalgebras 
onto parabolic subalgebras, the class of modules induced from 1-dimensional 
representations of parabolic subalgebras is preserved under passing to the duals 
in this way. Hence all indecomposable direct summands of $M^*$ and $\Md$ have 
$1$-dimensional sources, and so $\calM^*$ satisfies the assumptions of 
Proposition 5.1.

For each subspace $S\sbs\calM$ put $S^\perp=\{f\in\calM^*\mid f(S)=0\}$. Then
$\Si_i^\perp=\bigcap_{j>i}U_j^\perp$, $\Up_i^\perp=\sum_{j<i}U_j^\perp$, and
$$
\bigl(\Up_i/(\Up_i\cap\Si_i)\bigr)^*\cong(\Up_i^\perp+\Si_i^\perp)/\Up_i^\perp
\cong\Si_i^\perp/(\Up_i^\perp\cap\Si_i^\perp),
$$
which is the component of the complex 
$K_\bullet\bigl(\calM^*;(U^\perp_{n-1},\ldots,U^\perp_1)\bigr)$ in degree 
$n-i$. Thus $K_\bullet\bigl(\calM;(U_i)\bigr)\vphantom{)}^*$ is isomorphic to 
the complex $K_\bullet\bigl(\calM^*;(U^\perp_{n-1},\ldots,U^\perp_1)\bigr)$ 
with the degrees shifted by $n$. Note also that 
$U_{n-i}^\perp=\calM^*(\calT_{n-i}-1)=(\calT_i-1)\calM^*$ for each $i$. 
Therefore the latter complex is exact by Proposition 5.1.
\endproof

The next lemma provides a key ingredient in the proof of Theorem 6.6.

\proclaim
Lemma 5.10.
Suppose that $\calM=M'\ot M$ where $M$ is an $\calH_n$-module and $M'$ is an 
$\calH'_n$-module such that all indecomposable direct summands of $M$ and $M'$ 
have trivial sources. Put $\,\Si=\sum_{i=1}^{n-1}U_i,$ 
$\,\Si_1=\sum_{i>1}U_i\,$ where $U_i=\Ker\,(\calT_i-1)_\calM$ for each $i,$ and
$$
y=\sum_{\si\in\calD(\frS_{1,n-1}\backslash\frS_n)}
(-1)^{\ell(\si)}\calT_\si\in\calH'_n\ot\calH_n\qquad
\hbox{where $\calT_\si=T_\si T'_\si$}.
$$
Then $\,y\Si\sbs\Si_1\,$ and{\rm,} moreover{\rm,} $\,y_\calM^{-1}(\Si_1)=\Si$.
\endproclaim

\Proof.
The inclusion $y\Si\sbs\Si_1$ is proved exactly as in Lemma 2.2. For this one 
needs only the braid relations satisfied by $\calT_1,\ldots,\calT_{n-1}$. The 
equality $y_\calM^{-1}(\Si_1)=\Si$ will follow from injectivity of the map 
$\ph:\calM/\Si\to\calM/\Si_1$ induced by the linear operator $y_\calM$. Since 
all verifications can be done on direct summands of $\calM$ it suffices to 
consider the case when
$$
M=\calH_n\ot_{\calH_\la}\bbk\triv\,,\quad M'=\calH'_n\ot_{\calH'_\mu}\bbk\triv
\,,\quad\calM\cong(\calH'_n\ot\calH_n)\ot_{(\calH'_\mu\ot\calH_\la)}\bbk\triv
$$
for some parabolic subalgebras $\calH_\la\,$, $\calH'_\mu$. Let $c$ be the 
canonical generator of $\calM$. By Lemma 5.11 stated below a basis for 
$\calM/\Si$ is formed by the cosets of elements $T_\pi c$ with $\pi\in\muDla^e$ 
where
$$
\muDla^e=\{\pi\in\calD(\frS_\mu\backslash\frS_n/\frS_\la)\mid
\frS_\mu\cap\pi\frS_\la\pi^{-1}=e\}
$$
is the set of distinguished representatives of the double cosets with the 
trivial intersection property. There is a similar basis for $\calM/\Si_1$ 
obtained as follows. Consider the Mackey decompositions 
$$ 
M=\bigoplus_{\al\in\calD(\frS_{1,n-1}\backslash\frS_n/\frS_\la)}M(\al),\qquad
M'=\bigoplus_{\be\in\calD(\frS_{1,n-1}\backslash\frS_n/\frS_\mu)}M'(\be)
$$
with respect to the parabolic subalgebras $\calH_{1,n-1}$ and $\calH'_{1,n-1}$. 
We have
$$
M(\al)\cong\calH_{1,n-1}\ot_{\calH_{\nu(\al)}}\bbk\triv,\qquad
M'(\be)\cong\calH'_{1,n-1}\ot_{\calH'_{\nu'(\be)}}\bbk\triv
$$
where $\nu(\al)$ and $\nu'(\be)$ are the compositions of $n$ such that
$$
\frS_{\nu(\al)}=\frS_{1,n-1}\cap\al\frS_\la\al^{-1},\qquad
\frS_{\nu'(\be)}=\frS_{1,n-1}\cap\be\frS_\mu\be^{-1}.
$$
The $\calH'_{1,n-1}\ot\calH_{1,n-1}$-submodule $\calM(\al,\be)$ of $\calM$ 
generated by $T_\al T'_\be c$ is isomorphic to $M'(\be)\ot M(\al)$.
Since $\calM$ is a direct sum of these submodules for different $\al$ and $\be$, 
we get
$$
\calM/\Si_1\cong\bigoplus_{\al,\,\be}\,\calM(\al,\be)/\Si(\al,\be)\quad
\hbox{where $\Si(\al,\be)=\Si_1\cap\calM(\al,\be)$}.
$$
Note that $\calH_{1,n-1}\cong\calH_{n-1}(q)$, 
$\calH'_{1,n-1}\cong\calH_{n-1}(q^{-1})$. Thus Lemma 5.11 applies to each 
$\calM(\al,\be)$ viewed as an $\calH_{n-1}(q^{-1})\ot\calH_{n-1}(q)$-module.
We will need only those summands in the decomposition of $\calM/\Si_1$ which 
are indexed by the pairs $(\al,\be)$ with $\be=e$. Put $\mu^1=\nu'(e)$, so that 
$\frS_{\mu^1}=\frS_{1,n-1}\cap\frS_\mu$, and put
$$
\calD^e(\frS_{\mu^1}\backslash\frS_{1,n-1}/\frS_{\nu(\al)})
=\{\si\in\calD(\frS_{\mu^1}\backslash\frS_{1,n-1}/\frS_{\nu(\al)})\mid
\frS_{\mu^1}\cap\si\frS_{\nu(\al)}\si^{-1}=e\}.
$$
By Lemma 5.11 $\calM(\al,e)/\Si(\al,e)$ has a basis formed by the cosets of elements 
$$
T_\si T_\al c=T_{\si\al} c\quad\hbox{with 
$\si\in\calD^e(\frS_{\mu^1}\backslash\frS_{1,n-1}/\frS_{\nu(\al)})$}.
$$
The equality $T_\si T_\al=T_{\si\al}$ here is explained by the fact that 
$\si\in\frS_{1,n-1}$, while $\al$ is the shortest element in the coset 
$\frS_{1,n-1}\al$, so that $\ell(\si\al)=\ell(\si)+\ell(\al)$.

We claim that the assignment $(\si,\al)\mapsto\si\al$ gives a bijection
$$
\openup1\jot
\eqalign{
\{(\si,\al)\mid\al\in\calD(\frS_{1,n-1}\backslash\frS_n/\frS_\la),
\ \si\in\calD^e(\frS_{\mu^1}\backslash\frS_{1,n-1}/\frS_{\nu(\al)})\}
\hskip2cm & \cr
\mapr{}\muspDla=\{\pi\in\calD(\frS_{\mu^1}\backslash\frS_n/\frS_\la)\mid
\frS_{\mu^1}\cap\pi\frS_\la\pi^{-1}=e\}. \cr
}
$$
By the Mackey decomposition of coset representatives (see \cite{Geck-P, Lemma 
2.1.9}) $\calD_\la$ consists precisely of those elements $\pi\in\frS_n$ which 
can be written as $\pi=\si\al$ for some 
$\al\in\calD(\frS_{1,n-1}\backslash\frS_n/\frS_\la)$ and 
$\si\in\calD(\frS_{1,n-1}/\frS_{\nu(\al)})$. The pair $(\si,\al)$ is uniquely 
determined by $\pi$ since $\al$ is the shortest element in the double coset 
$\frS_{1,n-1}\pi\frS_\la$. Furthermore, for each $\tau_i\in\frB_{\mu^1}$ we 
have $\ell(\tau_i\pi)=\ell(\tau_i\si)+\ell(\al)$ since 
$\tau_i\si\in\frS_{1,n-1}$, whence $\tau_i\pi>\pi$ if and only if 
$\tau_i\si>\si$. This shows that $\pi\in\muspD$ if and only if $\si\in\muspD$. 

Lastly, for each $\rho\in\frS_{\mu^1}$ we have $\rho\si\in\frS_{1,n-1}$, and 
by the Mackey decomposition the equality $\rho\pi=(\rho\si)\al$ implies that 
$\rho\pi\in\calD_\la$ if and only if 
$\rho\si\in\calD(\frS_{1,n-1}/\frS_{\nu(\al)})$. 
The double coset $\frS_{\mu^1}\pi\frS_\la$ has the trivial intersection 
property if and only if $\rho\pi\in\calD_\la$ for all $\rho\in\frS_{\mu^1}$, 
while $\frS_{\mu^1}\si\frS_{\nu(\al)}$ has the trivial intersection 
property if and only if $\rho\si\in\calD(\frS_{1,n-1}/\frS_{\nu(\al)})$ for all 
$\rho\in\frS_{\mu^1}$. We see that these properties are equivalent. Thus 
$\pi\in\muspDla$ if and only if 
$\si\in\calD^e(\frS_{\mu^1}\backslash\frS_{1,n-1}/\frS_{\nu(\al)})$, and 
bijectivity of the map considered above has been established.

It follows from the preceding discussion that the vector space
$$
Q=\bigoplus_{\al\in\calD(\frS_{1,n-1}\backslash\frS_n/\frS_\la)}
\calM(\al,e)/\Si(\al,e)
$$
has a basis formed by the cosets of all elements $T_\pi c$ with $\pi\in\muspDla$. 
Consider the map $\psi:\calM/\Si\to Q$ obtained as the composite of $\ph$ with 
the projection $p_Q$ of $\calM/\Si_1$ onto $Q$. We will check that $\psi$ is 
injective. Once this has been done, the injectivity of $\ph$ will be clear, 
and the proof of Lemma 5.10 will be complete.

We have $\calD(\frS_{1,n-1}\backslash\frS_n)
=\{\si_j\mid j=0,\ldots,n-1\}$ where $\si_0=e$, $\si_1=\tau_1$, and generally 
$\si_j=\tau_1\tau_2\cdots\tau_j$ for $j>1$. Therefore 
$y=\sum_{j=0}^{n-1}(-1)^j\calT_{\si_j}$.

Let $\pi\in\muDla^e$, i.e. $\pi$ is the distinguished representative of a 
$\frS_\mu\,$-$\,\frS_\la$ double coset with the trivial intersection property. 
Then the coset $T_\pi c+\Si$ is a basis element of $\calM/\Si$ which is sent by 
$\psi$ to the element $p_Q(yT_\pi c+\Si_1)\in Q$. We have
$$
yT_\pi c=\sum_{j=0}^{n-1}(-1)^j\calT_{\si_j}T_\pi c
=\sum_{j=0}^{n-1}(-1)^jT_{\si_j}T_\pi T'_{\si_j}c\,.
$$
Note that $\si_j\tau_i>\si_j$ for all $i=1,\ldots,n-1$ except for $i=j$. 
Hence $\si_j\in\calD_\mu$ for $j>0$ if and only if $\tau_j\notin\frS_\mu$. If 
$\tau_j\in\frS_\mu$, then 
$T'_{\si_j}c=T'_{\si_{j-1}}T'_jc=q^{-1}T'_{\si_{j-1}}c$. For each $j$ it 
follows by induction that
$$
T'_{\si_j}c=q^{k(j)-j}\,T'_{k(j)}c\in\calM(e,\si_{k(j)})
$$
where $k(j)$ is the largest integer $k$ such that $0\le k\le j$ and 
$\si_k\in\calD_\mu$. This entails
$$
T_{\si_j}T_\pi T'_{\si_j}c\in\calH_n T'_{k(j)}c\sbs
\bigoplus_{\al\in\calD(\frS_{1,n-1}\backslash\frS_n/\frS_\la)}\calM(\al,\si_{k(j)}).
$$
Denote by $m$ the largest integer such that $0\le m<n$ and $k(m)=0$. Then 
$k(j)>0$, and therefore $\si_{k(j)}\ne e$ for all $j>m$. In particular, 
$p_Q(T_{\si_j}T_\pi T'_{\si_j}c+\Si_1)=0$ for $j>m$.

If $j\le m$, then $\si_j\in\frS_\mu$, whence $T_{\si_j}T_\pi=T_{\si_j\pi}$ 
with $\si_j\pi\in\calD_\la$ by the conditions on $\pi$. Moreover, 
$\si_j\pi\in\muspDla$. Indeed, for each $\rho\in\frS_{\mu^1}$ we have 
$\rho\si_j\pi\in\calD_\la$ since $\rho\si_j\in\frS_\mu$. This means that the 
double coset $\frS_{\mu^1}\si_j\pi\frS_\la$ has the trivial intersection 
property. But $\ell(\rho\si_j)=\ell(\rho)+\ell(\si_j)$ since 
$\rho\in\frS_{1,n-1}$, and therefore
$$
\ell(\rho\si_j\pi)=\ell(\rho\si_j)+\ell(\pi)=\ell(\rho)+\ell(\si_j)+\ell(\pi)
=\ell(\rho)+\ell(\si_j\pi). 
$$
This shows that $\si_j\pi\in\muspD$.

We conclude that $\psi(T_\pi c+\Si)$ equals 
$\sum_{j=0}^m(-1)^jq^{-j}(T_{\si_j\pi}c+\Si_1)$, which is a linear combination 
of distinct basis elements of $Q$ with nonzero coefficients. Note that all 
elements $\si_j\pi$ with $j=0,\ldots,m$ belong to the same coset $\frS_\mu\pi$
having $\pi$ as its shortest representative. Therefore the expressions for the 
images under $\psi$ of two different basis elements of $\calM/\Si$ involve 
disjoint sets of basis elements of $Q$. Injectivity of $\psi$ and $\ph$ is 
now clear.
\endproof

\proclaim
Lemma 5.11.
Let $M=\calH_n\ot_{\calH_\la}\bbk\triv$ and $M'=\calH'_n\ot_{\calH'_\mu}\bbk\triv$. 
With the notation as in Lemma\/ {\rm5.10} the set $\{T_\pi c+\Si\mid\pi\in\muDla^e\}$ 
is a basis for $\calM/\Si$.
\endproclaim

\Proof.
Consider $M'$ as a right $\calH_n$-module by means of the algebra antiisomorphism 
$\calH_n\to\calH'_n$ sending $T_i$ to $(T'_i)^{-1}$. Then 
$M'\cong\bbk\triv\ot_{\calH_\mu}\calH_n$ and $\calM=M'\ot M$ is an 
$\calH_n$-bimodule with respect to the $\calH_n$-module structures on $M$ 
and $M'$. We have
$$
U_i=\{u\in\calM\mid T_iT'_iu=u\}=\{u\in\calM\mid T_iu=uT_i\}
$$
for each $i=1,\ldots,n-1$. If $u=T_im-m\Tt_i$ for some $m\in\calM$ where 
$\Tt_i=q-1-T_i$, then $u\in U_i$ since
$$
T_iu-uT_i=T^2_im-T_im(T_i+\Tt_i)+mT_i\Tt_i=(T^2_i-(q-1)T_i-q)m=0.
$$
This shows that $U_i\sps\Ut_i$ where $\Ut_i=\{T_im-m\Tt_i\mid m\in\calM\}$. If 
$q\ne-1$, then $U_i=\Ut_i$ since each element $u\in U_i$ can be written as 
$$
u=(q+1)^{-2}(T_i-\Tt_i)^2u=T_im-m\Tt_i\quad
\hbox{with $m=(q+1)^{-2}(T_i-\Tt_i)u$}.
$$
Note that $\Ut_i$ is spanned by all elements $w\ot T_iv-w\Tt_i\ot v$ with 
$v\in M$ and $w\in M'$. Setting $\,\Sit=\sum_{i=1}^{n-1}\Ut_i\,$, we deduce that
$$
\calM/\Sit\cong\Mtt\ot_{\calH_n}M
$$
where $\Mtt$ is the right $\calH_n$-module obtained from $M'$ by composing the 
original action of $\calH_n$ with the automorphism of $\calH_n$ sending $T_i$ 
to $\Tt_i$ for each $i$. By the right hand version of Lemma 1.5 we have 
$\Mtt\cong\bbk\alt\ot_{\calH_\mu}\calH_n$. Hence 
$$ 
\Mtt\ot_{\calH_n}M\cong\bbk\alt\ot_{\calH_\mu}M
\cong\bigoplus_{\pi\in\,\muDlass}\bbk\alt\ot_{\calH_{\nu(\pi)}}\bbk\triv
$$
since $M\cong\bigoplus\,(\calH_\mu\ot_{\calH_{\nu(\pi)}}\bbk\triv)$ by the 
Mackey formula. Here $\nu(\pi)$ is the composition of $n$ such that 
$\frS_{\nu(\pi)}=\frS_\mu\cap\pi\frS_\la\pi^{-1}$. For each $\pi$ the 
respective summand in the above decomposition of $\Mtt\ot_{\calH_n}M$ is 
spanned by the image of $T_\pi c$.

In view of the previous isomorphisms $\calM/\Sit$ has a basis consisting of the 
cosets of elements $T_\pi c$ with $\pi\in\muDla$ such that the alternating 
representation of $\calH_{\nu(\pi)}$ coincides with the trivial representation.  
If $q\ne-1$, this condition on $\pi$ means precisely that $\frS_{\nu(\pi)}=e$, 
i.e. $\pi\in\muDla^e$. In this case we also have $\Si=\Sit$, and the 
conclusion of Lemma 5.11 follows.

Suppose further that $q=-1$. In this case $\calM/\Si$ is a quotient of 
$\calM/\Sit$ since $\Sit\sbs\nobreak\Si$. Hence $\calM/\Si$ is spanned by the 
cosets of elements $T_\pi c$ with $\pi\in\muDla$. If 
$\tau_i\in\frS_{\nu(\pi)}$ for some $i$, then $T_\pi c\in U_i\sbs\Si$ since 
$T_iT_\pi c=-T_\pi c$ and $T'_iT_\pi c=-T_\pi c$. It follows that $\calM/\Si$ 
is spanned by the cosets of elements $T_\pi c$ with $\pi\in\muDla^e$ only. It 
remains to prove that such cosets are linearly independent. But the dual space 
$(\calM/\Si)^*$ is identified with the subspace of the $\calH_n$-bimodule
$$
\Homk\bigl(M,\Md\bigr)\cong\calM^*
$$
consisting of all $\bbk$-linear maps $f:M\to\Md$ with the property that for 
each $i$, $0<i<n$, there exists a $\bbk$-linear map $f_i:M\to\Md$ 
such that $f=T_if_i-f_iT_i$, i.e. $f=T_if_i+f_i\Tt_i+2f_i$. Here $\Md$ is a left 
$\calH_n$-module with respect to the action of $\calH_n$ arising naturally from 
the right action on $M'$. By Lemmas 1.4, 1.5
$$
\Md\cong\calH_n\ot_{\calH_\mu}\bbk\,,\qquad
\Mt\cong\calH_n\ot_{\calH_\la}\bbk.
$$
By Lemma 1.7 the $\bbk$-linear maps $f$ considered above are precisely those 
$\calH_n$-module homomorphisms $\Mt\to\Md$ that factor through a free 
module, and the space of such homomorphisms has a basis indexed by the set 
$\laDmu^e$. Since this set is in a bijection with $\muDla^e$ by the map 
$\si\mapsto\si^{-1}$, the space $\calM/\Si$ has dimension equal to the 
cardinality of $\muDla^e$, and we are done.
\endproof

\section
6. Intertwining algebras for a pair of Hecke symmetries

Let $V$ and $V'$ be two finite dimensional vector spaces over the field $\bbk$.
Let $R$ be a Hecke symmetry on $V$ and $R'$ a Hecke symmetry on $V'$ 
satisfying the Hecke relation with the same parameter $q$. For each 
$n\ge0$ we will view $\bbT_n(V)$ and $\bbT_n(V')$ as left modules over the 
Hecke algebra $\calH_n=\calH_n(q)$ with respect to the representations arising 
from $R$ and $R'$.

The tensor algebra $\bbT(\Vd\ot V)$ embeds canonically into 
$\bbT(\Vd)\ot\bbT(V)$. Under this embedding $\bbT_n(\Vd\ot V)$ is mapped 
onto $\bbT_n(\Vd)\ot\bbT_n(V)$. Identify $\bbT_n(V^*)$ and $\bbT_n(\Vd)$ with 
the duals of the vector spaces $\bbT_n(V)$  and $\bbT_n(V')$ as in section 3. 

Let $\Rd$ be the Hecke symmetry on $\Vd$ adjoint to $R'$ (see Lemma 3.2). The 
inverse operator $(\Rd)^{-1}$ is a Hecke symmetry with parameter $q^{-1}$. 
Denote by $\calR$ the linear operator on $\bbT_2(\Vd\ot V)$ which corresponds 
to the operator $(\Rd)^{-1}\ot R$ acting on $\bbT_2(\Vd)\ot\bbT_2(V)$.

Define $A(R',R)$ and $E(R',R)$ as the factor algebras of $\bbT(\Vd\ot V)$ by 
the ideals generated, respectively, by $\Im\,(\calR-\Id)$ and 
$\Ker\,(\calR-\Id)$. Under the isomorphism 
$$
\bbT_2(\Vd\ot V)\cong\bbT_2(\Vd)\ot\bbT_2(V)
$$
these two subspaces of $\bbT_2(\Vd\ot V)$ are mapped, respectively, onto
$$
\Im\,(\Id\ot R-\Rd\ot\Id)\quad{\rm and}\quad\Ker\,(\Id\ot R-\Rd\ot\Id).
$$
In the case when $R'=R$ the first subspace gives a well-known presentation of 
the FRT bialgebra $A(R)=A(R,R)$ (see \cite{Hay92, section 4}).

\proclaim
Lemma 6.1.
The quadratic dual algebras $A(R',R)^!$ and $E(R',R)^!$ are isomorphic{\rm,} 
respectively{\rm,} to $E(R,R')$ and $A(R,R')$. If $q\ne-1,$ then 
$A(R',R)=E(R',\Rt)$ and $E(R',R)=A(R',\Rt)$ where $\Rt=(q-1)\Id-R$.  
\endproclaim

\Proof.
Identifying the dual space of $\bbT_2(\Vd)\ot\bbT_2(V)$ with 
$\bbT_2(V^*)\ot\bbT_2(V')$, we have
$$
\displaylines{
\bigl(\Im\,(\Id\ot R-\Rd\ot\Id)\bigr)^\perp=\Ker\,(\Id\ot R'-R^*\ot\Id),\cr
\bigl(\Ker\,(\Id\ot R-\Rd\ot\Id)\bigr)^\perp=\Im\,(\Id\ot R'-R^*\ot\Id).\cr
}
$$
The subspaces of $\bbT_2(V^*)\ot\bbT_2(V')$ in the left hand sides of these 
equalities define the algebras $A(R',R)^!$ and $E(R',R)^!$.

If $q\ne-1$, then $\Id\ot R$ and $\Rd\ot\Id$ are commuting diagonalizable 
operators with two eigenvalues $-1$ and $q$. Since $\bbT_2(\Vd)\ot\bbT_2(V)$ 
is a sum of common eigenspaces of these two operators, it is clear that
$$
\displaylines{
\openup1\jot
\Im\,(\Id\ot\Rt-\Rd\ot\Id)=\Ker\,(\Id\ot R-\Rd\ot\Id),\cr
\Ker\,(\Id\ot\Rt-\Rd\ot\Id)=\Im\,(\Id\ot R-\Rd\ot\Id).\cr
}
$$
We thus obtain the second conclusion.
\endproof

Put $\calH'_n=\calH_n(q^{-1})$, as in section 5. Consider $\bbT_n(\Vd)$ as a 
left $\calH'_n$-module with respect to the representation arising from the 
Hecke symmetry $(\Rd)^{-1}$. We thus obtain a left $\calH'_n\ot\calH_n$-module 
structure on $\bbT_n(\Vd\ot V)\cong\bbT_n(\Vd)\ot\bbT_n(V)$.

Recall that $\calT_i=T'_i\ot T_i\in\calH'_n\ot\calH_n$ where 
$T'_1,\ldots,T'_{n-1}$ are the standard generators of $\calH'_n$. The action 
of $\calT_i$ on $\bbT_n(\Vd\ot V)$ is given by the operator 
$$
\calR_i^{(n)}=\Id^{\ot(i-1)}\ot\calR\ot\Id^{\ot(n-i-1)}.
$$

\proclaim
Theorem 6.2.
Suppose that both $R$ and $R'$ satisfy the $1$-dimensional source condition. 
Then the graded algebras $A(R',R)$ and $E(R',R)$ are Koszul. Their Hilbert 
series satisfy the relation \ $h_{A(R',R)}(t)\,h_{E(R',R)}(-t)=1$.
\endproclaim

\Proof.
Put $U_i^{(n)}=\bbT_{i-1}(\Vd\ot V)\ot U\ot\bbT_{n-i-1}(\Vd\ot V)\sbs
\bbT_n(\Vd\ot V)$ for a subspace $U\sbs\bbT_2(\Vd\ot V)$ and $i=1,\ldots,n-1$. 
We have
$$
U_i^{(n)}=\cases{
\Im\,(\calR_i^{(n)}-\Id)=\Im\,(\calT_i-1)_{\bbT_n(V'^*\ot V)} &
if $\,U=\Im\,(\calR-\Id)$,\cr
\noalign{\smallskip}
\Ker\,(\calR_i^{(n)}-\Id)=\Ker\,(\calT_i-1)_{\bbT_n(V'^*\ot V)} &
if $\,U=\Ker\,(\calR-\Id)$.\cr
}
$$
In both cases the complex $K_\bullet\bigl(\bbT_n(\Vd\ot V),(U_i^{(n)})\bigr)$ 
is exact for each $n>0$ by Propositions 5.1 and 5.9. But these complexes are 
precisely the direct summands in the decomposition of the right Koszul 
complex $K_\bullet(A)$ (see section 3) where $A$ is $A(R',R)$ in the first 
case and $E(R',R)$ in the second. Hence $K_\bullet(A)$ is acyclic in all 
positive degrees.

Thus $A(R',R)$ and $E(R',R)$ are Koszul. The Hilbert series of the algebra 
$E(R',R)^!$ is $\sum\,(\dim\Up^{(n)})t^n$ where the spaces $\Up^{(n)}$ are 
determined in Lemma 6.4 below. Making use also of Lemmas 6.3 and 1.2, we get
$$
\openup1\jot
\eqalign{
\dim A_n(R',R)&{}=\dim\Hom_{\calH_n}\bigl(\bbT_n(V),\bbT_n(V')\bigr)\cr 
&{}=\dim\Hom_{\calH_n}\bigl(\bbT_n(V'),\bbT_n(V)\bigr)=\dim\Up^{(n)}.\cr
}
$$
This shows that $h_{A(R',R)}(t)=h_{E(R',R)^!}(t)$, and the final conclusion in 
the statement of Theorem 6.2 reduces to the standard relation between the 
Hilbert series of the Koszul algebra $E(R',R)$ and its quadratic dual.
\endproof

\proclaim
Lemma 6.3.
$\,\,A_n(R',R)\cong\Hom_{\calH_n}\bigl(\bbT_n(V),\bbT_n(V')\bigr)^*$.
\endproclaim

\Proof.
The ideal $I$ of the algebra $\bbT(\Vd\ot V)$ defining its 
factor algebra $A(R',R)$ has homogeneous components $I_n=0$ for $n\le1$ and 
$I_n=\sum_{i=1}^{n-1}U_i^{(n)}$ for $n>1$ where $U_i^{(n)}$ are the subspaces 
of $\bbT_n(\Vd\ot V)$ defined in the proof of Theorem 6.2 with 
$U=\Im\,(\calR-\Id)$.

The right $\calH_n$-module structure on $\bbT_n(V')^*$ obtained in a natural 
way from the left module structure on $\bbT_n(V')$ allows us to view 
$\bbT_n(\Vd\ot V)\cong\bbT_n(V')^*\ot\bbT_n(V)$ as an $\calH_n$-bimodule. 
The left action of $T'_i$ on $\bbT_n(\Vd)\cong\bbT_n(V')^*$ is the same as the 
right action of $T_i^{-1}$. Then $T'_ia=aT_i^{-1}$, and so 
$\calT_ia=T_iT'_ia=T_iaT_i^{-1}$, for all $a\in\bbT_n(\Vd\ot V)$. 
Hence $U_i^{(n)}$ is spanned by all elements $T_iaT_i^{-1}-a$ or, 
equivalently, by all elements $T_ia-aT_i$ with $a\in\bbT_n(\Vd\ot V)$.
Under the canonical isomorphisms of $\calH_n$-bimodules
$$
\bbT_n(\Vd\ot V)^*\cong\bbT_n(V')\ot\bbT_n(V)^*\cong
\Homk\bigl(\bbT_n(V),\bbT_n(V')\bigr)
$$
the orthogonal of $U_i^{(n)}$ in that vector space is
$$
(U_i^{(n)})^\perp=\{f\in\Homk\bigl(\bbT_n(V),\bbT_n(V')\bigr)\mid T_if=fT_i\}.
$$
Hence
$\,A_n(R',R)^*\cong I_n^\perp=\bigcap_{i=1}^{n-1}(U_i^{(n)})^\perp
=\Hom_{\calH_n}\bigl(\bbT_n(V),\bbT_n(V')\bigr)$.
\endproof

\proclaim
Lemma 6.4.
The right Koszul complex for the algebra $E(R',R)$ has components 
$$
\Up^{(n)}\ot E(R',R)\quad\hbox{with 
$\,\Up^{(n)}\cong\Hom_{\calH_n}\bigl(\bbT_n(V'),\bbT_n(V)\bigr)$}.
$$
\endproclaim

\Proof.
Recall from section 3 that $\Up^{(n)}=\bigcap_{i=1}^{n-1}U_i^{(n)}$ where 
$U_i^{(n)}$ are as defined in the proof of Theorem 6.2 with $U=\Ker\,(\calR-\Id)$. 
In terms of the $\calH_n$-bimodule structure on $\bbT_n(\Vd\ot V)$ we have 
$$
U_i^{(n)}=\{a\in\bbT_n(\Vd\ot V)\mid T_ia=aT_i\}
$$
for each $i$. It follows that $\Up^{(n)}=\{a\in\bbT_n(\Vd\ot V)\mid xa=ax 
\hbox{ for all }x\in\calH_n\}$. The canonical isomorphisms of 
$\calH_n$-bimodules
$$
\bbT_n(\Vd\ot V)\cong\bbT_n(V)\ot\bbT_n(V')^*\cong 
\Homk\bigl(\bbT_n(V'),\bbT_n(V)\bigr)
$$
map $\Up^{(n)}$ onto $\Hom_{\calH_n}\bigl(\bbT_n(V'),\bbT_n(V)\bigr)$.
\endproof

\proclaim
Lemma 6.5.
Let $I$ be the ideal of the algebra $\bbT=\bbT(\Vd\ot V)$ defining its 
factor algebra $E(R',R)$. Fix some $n>1$ and put
$$
\openup1\jot
\displaylines{
L_k=\{a\in\bbT_k\mid a\bbT_{n-k}\sbs I_n\}\qquad\hbox{\rm for }0\le k<n,\cr
y_k=\sum_{\si\in\calD(\frS_{1,k-1}\backslash\frS_k)}
(-1)^{\ell(\si)}\calT_\si\in\calH'_k\ot\calH_k\qquad\hbox{\rm for each }k>0\cr
}
$$
where $\calT_\si=T_\si T'_\si\,$. If\/ $0<k<n,$ then $\,y_kL_k\sbs\bbT_1L_{k-1}$.
\endproclaim

This lemma is proved by exactly the same arguments as those used for Lemma 4.1 
(see the remarks following the proof of that lemma).

\proclaim
Theorem 6.6.
Suppose that both $R$ and $R'$ satisfy the trivial source condition. 
If $\dim E_n(R',R)=1$ and $E_{n+1}(R',R)=0$ for some integer 
$n>0,$ then $E(R',R)$ is a Frobenius algebra{\rm,} while $A(R',R)$ 
is a Gorenstein algebra of global dimension $n$.
\endproclaim

\Proof.
We proceed as in the proof of Theorem 4.5. Let $\bbT$ and $I$ be as in Lemma 6.5. 
By induction on $k$ we can show that $L_k=I_k$ for each $k=0,\ldots,n-1$. 
Indeed, if the equality $L_{k-1}=I_{k-1}$ holds for some $k$, then 
$y_kL_k\sbs\bbT_1I_{k-1}$ by Lemma 6.5. Now we apply Lemma 5.10 with $n$ 
replaced by $k$ and $\calM=\bbT_k$. In the notation of that lemma we then have 
$y=y_k$, $\Si=I_k$ and $\Si_1=\bbT_1I_{k-1}$. The inclusion $yL_k\sbs\Si_1$ 
entails $L_k\sbs y_\calM^{-1}(\Si_1)=\Si$, i.e. $L_k=I_k$. Thus the 
multiplication pairing
$$
E_k(R',R)\times E_{n-k}(R',R)\to E_n(R',R)
$$
has zero left kernel. Since this holds also with $k$ replaced by $n-k$, 
comparison of dimensions shows that the pairing is nondegenerate. This means 
that the algebra $E(R',R)$ is Frobenius. By Lemma 6.1 $A(R',R)^!\cong E(R,R')$. 
Since
$$
h_{E(R,R')}(t)=h_{A(R',R)}(-t)^{-1}=h_{E(R',R)}(t),
$$
we have $\dim E_n(R,R')=1$ and $E_{n+1}(R,R')=0$. Hence $E(R,R')$ is also 
Frobenius, and $A(R',R)$ is Gorenstein by \cite{Pol-P, Remark 2 on p. 25}. 
\endproof

\section
7. Monoidal equivalences

Let $V,V',V''$ be three finite dimensional vector spaces over $\bbk$ and 
$R,R',R''$ Hecke symmetries on the respective spaces with the same parameter 
$q$. For each $n\ge0$ we equip $\,\bbT_n=\bbT_n(V)$, $\,\bbT'_n=\bbT_n(V')\,$ 
and $\,\bbT''_n=\bbT_n(V'')\,$ with the $\calH_n$-module structures arising 
from $R,R',R''$. There is a $\bbk$-linear map
$$
\De_n:\ A_n(R',R'')\mapr{}A_n(R',R)\ot A_n(R,R'')
$$
whose dual map $\,\De_n^*:A_n(R',R)^*\ot A_n(R,R'')^*\mapr{}A_n(R',R'')^*\,$ 
is, in terms of the identifications of Lemma 6.3, the map
$$
\Hom_{\calH_n}(\bbT_n,\bbT'_n)\ot
\Hom_{\calH_n}(\bbT''_n,\bbT_n)\mapr{}
\Hom_{\calH_n}(\bbT''_n,\bbT'_n)
$$
given by the composition of homomorphisms. In particular, $A_n(R)$ is endowed 
with a comultiplication dual to the multiplication in the algebra 
$\,\End_{\calH_n}\!\bbT_n$. Thus $A_n(R)$ is a coalgebra. Also, $A_n(R',R)$ 
has an $A_n(R'),A_n(R)$-bicomodule structure dual to the 
$\,\End_{\calH_n}\!\bbT'_n\,,\End_{\calH_n}\!\bbT_n\,$-bimodule structure on 
$\Hom_{\calH_n}(\bbT_n,\bbT'_n)$.

Let $\calH_{m,n}$ be the parabolic subalgebra of $\calH_{m+n}$ generated by 
$\{T_i\mid i\ne m\}$. Then
$$
\eqalign{
A_m(R',R)^*\ot A_n(R',R)^*
& \cong\Hom_{\calH_m}(\bbT_m,\bbT'_m)\ot
\Hom_{\calH_n}(\bbT_n,\bbT'_n)\cr
& \cong
\Hom_{\calH_m\ot\calH_n}(\bbT_m\ot\bbT_n,\,\bbT'_m\ot\bbT'_n)\cr
& \cong\Hom_{\calH_{m,n}}(\bbT_{m+n},\bbT'_{m+n}).\cr
}
$$
Noting that the multiplication maps 
$\,A_m(R',R)\ot A_n(R',R)\mapr{}A_{m+n}(R',R)\,$ are dual to the inclusion maps
$$
\Hom_{\calH_{m+n}}(\bbT_{m+n},\bbT'_{m+n})\hrar
\Hom_{\calH_{m,n}}(\bbT_{m+n},\bbT'_{m+n}),
$$
it is easy to see that the map
$$
\De_{R'\!,R,R''}:\ A(R',R'')\mapr{}A(R',R)\ot A(R,R'')
$$
obtained from the previously defined maps $\De_n$ on the homogeneous 
components, is an algebra homomorphism. In particular, $A(R)$ is a bialgebra, 
while $A(R',R)$ is an $A(R'),A(R)$-bicomodule algebra.

If $C$ is a coalgebra, $\rho:X\to X\ot C$ and $\la:Y\to C\ot Y$ are right and 
left comodule structures on two vector spaces, then the \emph{cotensor 
product}\/ $X\square_CY$ is defined as the kernel of the $\bbk$-linear map
$$
X\ot Y\lmapr9{\rho\ot\id-\id\ot\la}X\ot C\ot Y.
$$
Suppose that $X$, $Y$ and $C$ are finite dimensional. Passing to the dual 
spaces, we get then an exact sequence
$$
X^*\ot C^*\ot Y^*\lmapr{10}{\rho^*\ot\id-\id\ot\la^*}X^*\ot Y^*
\mapr{}(X\square_CY)^*\mapr{}0
$$
which shows that $(X\square_CY)^*\cong X^*\ot_{C^*}Y^*$. Since $\De_n^*$ 
factors through
$$
\Hom_{\calH_n}(\bbT_n,\bbT'_n)\ot_{\End_{\calH_n}\!\!\bbT_n}
\Hom_{\calH_n}(\bbT''_n,\bbT_n),
$$
it follows that $\,\Im\De_n\sbs A_n(R',R)\square_{A_n(R)}A_n(R,R'')$ for each 
$n$. Hence $\De_{R'\!,R,R''}$ is in fact an algebra homomorphism
$\,A(R',R'')\mapr{}A(R',R)\square_{A(R)}A(R,R'')$.

\proclaim
Lemma 7.1.
Suppose that for each $n>1$ each indecomposable direct summand of the 
$\calH_n$-module $\bbT_n(V')$ is isomorphic to a direct summand of the 
$\calH_n$-module $\bbT_n(V)$. Then $\De_{R'\!,R,R''}$ gives an isomorphism of 
algebras
$$
A(R',R'')\cong A(R',R)\square_{A(R)}A(R,R'').
$$
\endproclaim

\Proof.
For any $\calH_n$-modules $X,Y,Z$ there is a canonical map
$$
\Hom_{\calH_n}(X,Y)\ot_{\End_{\calH_n}\!\!X}\Hom_{\calH_n}(Z,X)\mapr{}
\Hom_{\calH_n}(Z,Y)
$$
which is obviously bijective when $Y=X$. Since the collection of these maps 
with varying $Y$ give a natural transformation of additive functors of $Y$, 
such a map is bijective also when $Y$ is a direct sum of modules isomorphic to 
direct summands of $X$. In particular, we may take $X=\bbT_n(V)$, 
$Y=\bbT_n(V')$, $Z=\bbT_n(V'')$. Passing to the dual spaces, we deduce that 
$\De_n$ maps $A_n(R',R'')$ bijectively onto 
$A_n(R',R)\square_{A_n(R)}A_n(R,R'')$.
\endproof

\proclaim
Theorem 7.2.
Suppose that for each $n>1$ the indecomposable $\calH_n$-modules isomorphic to 
direct summands of\/ $\bbT_n(V')$ are the same as those isomorphic to direct 
summands of\/ $\bbT_n(V)$. Then the functors
$$
A(R',R)\square_{A(R)}\hbox{}
?\qquad{\rm and}\qquad?\,\hbox{}\square_{A(R)}A(R,R')
$$
are braided monoidal equivalences $\,\ARM\mapr{}\ARpM\,$ and $\,\MAR\mapr{}\MARp$.
\endproclaim

\Proof.
The functor $F=A(R',R)\square_{A(R)}\hbox{}?$ has a quasiinverse 
$F'=A(R,R')\square_{A(R')}\hbox{}?$ since
$$
F'F(X)\cong
\bigl(A(R,R')\square_{A(R')}A(R',R)\bigr)\square_{A(R)}X\cong
A(R)\square_{A(R)}X\cong X
$$
for left $A(R)$-comodules $X$ by the associativity of cotensor products and by 
Lemma 7.1, and, similarly, $FF'\cong\Id$. Thus $F$ is a category equivalence.

There are homomorphisms of $A(R')$-comodules $\xi_{XY}:F(X)\ot F(Y)\to F(X\ot Y)$, 
natural in $X,Y\in\ARM$, obtained as restrictions of the maps
$$
\displaylines{
\bigl(A(R',R)\ot X\bigr)\ot\bigl(A(R',R)\ot y\bigr)\to A(R',R)\ot(X\ot Y),\cr
(a\ot x)\ot(b\ot y)\mapsto ab\ot(x\ot y).\cr
}
$$
We will show that $\xi_{XY}$ is an isomorphism. Since $A(R)$ is a direct sum 
of its subcoalgebras $A_n(R)$, each left $A(R)$-comodule $X$ can be written as 
$X=\oplus_{n=0}^\infty X_n$ where $X_n$ is a left $A_n(R)$-comodule for each 
$n$. Therefore it suffices to prove bijectivity of $\xi_{XY}$ assuming $X$ to 
be an $A_m(R)$-comodule, $Y$ an $A_n(R)$-comodule for some $m,n$. Since every 
comodule is a sum of finite dimensional subcomodules, we may also assume that 
$\dim X<\infty$ and $\dim Y<\infty$. Then
$\,F(X)=A_m(R',R)\square_{A_m(R)}X$, and
$$
F(X)^*\cong A_m(R',R)^*\ot_{A_m(R)^*}X^*\cong
\Hom_{\calH_m}(\bbT_m,\bbT'_m)\ot_{\End_{\calH_m}\!\!\bbT_m}X^*.
$$
where $\,\bbT_m=\bbT_m(V)$, $\,\bbT'_m=\bbT_m(V')$. \ Similarly,
$$
F(Y)^*\cong
\Hom_{\calH_n}(\bbT_n,\bbT'_n)\ot_{\End_{\calH_n}\!\!\bbT_n}Y^*.
$$
Identifying $\calH_m\ot\calH_n$ with the subalgebra $\calH_{m,n}$ of 
$\calH_{m+n}$, we get
$$
\eqalign{
&F(X)^*\ot F(Y)^*\cr
&\qquad\qquad\cong
\Hom_{\calH_m\ot\calH_n}(\bbT_m\ot\bbT_n,\bbT'_m\ot\bbT'_n)
\ot_{\End_{\calH_m\ot\calH_n}(\bbT_m\ot\bbT_n)}(X^*\ot Y^*)\cr
&\qquad\qquad\cong
\Hom_{\calH_{m,n}}(\bbT_{m+n},\bbT'_{m+n})\ot_{\End_{\calH_{m,n}}\!\!\bbT_{m+n}}
(X^*\ot Y^*).\cr
}
$$
The dual map $F(X\ot Y)^*\to F(X)^*\ot F(Y)^*$ is identified with the 
canonical map
$$
\eqalign{
&\Hom_{\calH_{m+n}}(\bbT_{m+n},\bbT'_{m+n})\ot_{\End_{\calH_{m+n}}\!\!\bbT_{m+n}}
(X^*\ot Y^*)\cr
&\qquad\qquad\qquad\qquad
\mapr{}\Hom_{\calH_{m,n}}(\bbT_{m+n},\bbT'_{m+n})\ot_{\End_{\calH_{m,n}}\!\!\bbT_{m+n}}
(X^*\ot Y^*)\cr
}
$$
arising from the inclusions 
$\Hom_{\calH_{m+n}}(\bbT_{m+n},\bbT'_{m+n})\sbs\Hom_{\calH_{m,n}}(\bbT_{m+n},\bbT'_{m+n})$ 
and $\End_{\calH_{m+n}}\!\!\bbT_{m+n}\sbs\End_{\calH_{m,n}}\!\!\bbT_{m+n}$. 
This map is bijective by Lemma 7.3 below.

We have shown that $F(X)\ot F(Y)\cong F(X\ot Y)$, naturally in $X$ and $Y$. 
Coherence of these isomorphisms is clear from the construction. The trivial 
$A(R)$-comodule $A_0(R)=\bbk$ is sent by $F$ to the trivial $A(R')$-comodule 
$A_0(R',R)=\bbk$. Thus $F$ is a monoidal equivalence.

Let $b$ and $b'$ be the braidings in the categories $\ARM$ and $\ARpM$, 
respectively, such that $b_{\hbox{$\scriptstyle V^*V^*$}}=R^*$ and 
$b_{\hbox{$\scriptstyle\Vd\Vd$}}=\Rd\!$. We have to verify commutativity of 
the diagrams
$$ 
\openup1\jot
\diagram{
F(X)\ot F(Y) & \lmapr5{\xi_{XY}} & F(X\ot Y)\hphantom{.} \cr
\lmapd{20pt}{b'_{F(X)F(Y)}}{} && \lmapd{20pt}{}{F(b_{XY})} \cr
F(Y)\ot F(X) & \lmapr5{\xi_{YX}} & F(Y\ot X). \cr
}
$$
In fact, it suffices to do this only for $X=Y=V^*$. By the general properties 
of the braidings the diagram will then be commutative for $X=\bbT_m(V^*)$, 
$Y=\bbT_n(V^*)$, and therefore also when $X$ and $Y$ are subfactors of direct sums 
of left $A(R)$-comodules isomorphic to tensor powers of $V^*$. But every 
finite dimensional left $A(R)$-comodule is realized in this way since 
$\bbT_n(V^*)\cong\bbT_n^*$ is a faithful right $\End_{\calH_n}\!\bbT_n$-module 
for each $n$. For infinite dimensional comodules commutativity of the diagram 
will follow from the fact that $F$ commutes with inductive direct limits. Now
$$ 
F(\bbT_n^*)^*\cong
\Hom_{\calH_n}(\bbT_n,\bbT'_n)\ot_{\End_{\calH_n}\!\!\bbT_n}\bbT_n\cong\bbT'_n\,.
$$
The last isomorphism here is explained by the fact that the evaluation map 
$$
\Hom_{\calH_n}(\bbT_n,M)\ot_{\End_{\calH_n}\!\!\bbT_n}\bbT_n\to M
$$
is obviously bijective when $M=\bbT_n$, and therefore it is bijective whenever 
$M$ is a direct sum of $\calH_n$-modules isomorphic to direct summands of 
$\bbT_n$. This can be applied with $M=\bbT'_n$. Note that the bijection 
$F(\bbT_n^*)^*\cong\bbT'_n$ obtained in this way is an isomorphism of 
$\calH_n$-modules.

In particular, we have $F(V^*)\cong\Vd$ and $F(\bbT_2^*)\cong{\bbT'_2}^{\!*}$. 
The generator $T_1$ of $\calH_2$ acts via $R$ on $\bbT_2$ and via $R'$ on 
$\bbT'_2$. This entails the commutativity of the diagram
$$
\openup1\jot
\diagram{
{\bbT'_2}^{\!*} & {}\cong{} & F(V^*)\ot F(V^*) & \lmapr7{\xi_{V^*V^*}} & 
F(\bbT_2^*)\hphantom{.} \cr \lmapd{20pt}{\Rd}{} && 
\lmapd{20pt}{}{b'_{F(V^*)F(V^*)}} && \lmapd{20pt}{}{F(R^*)} \cr 
{\bbT'_2}^{\!*} & {}\cong{} & F(V^*)\ot F(V^*) & \lmapr7{\xi_{V^*V^*}} & 
F(\bbT_2^*). \cr
}
$$
Thus we have verified all the required properties of the functor $F$. Consideration 
of the other case in the statement of Theorem 7.2 is quite similar.
\endproof

\proclaim
Lemma 7.3.
Let $A$ be a ring{\rm,} $B$ its subring{\rm,} $X$ a left $A$-module{\rm,} and 
let $M$ be a left $\End_BX$-module. The canonical map
$$ 
\Hom_A(X,Y)\ot_{\End_A\!X}M\mapr{}\Hom_B(X,Y)\ot_{\End_B\!
X}M
$$ 
is bijective whenever $Y$ is a finite direct sum of left $A$-modules 
isomorphic to direct summands of $X$.
\endproclaim

\Proof.
The conclusion is obvious when $Y=X$. Since we deal here with a natural 
transformation of two additive functors of $Y$, the conclusion then holds 
in full generality.
\endproof

\references
\nextref Back81
\auth{J.,Backelin}
\paper{A distributiveness property of augmented algebras and some related homological results}
PhD thesis, Stockholm, 1981.

\nextref Bich03
\auth{J.,Bichon}
\paper{The representation category of the quantum group of a non-degenerate bilinear form}
\journal{Comm. Algebra}
\Vol{310}
\Year{2003}
\Pages{4831-4851}

\nextref Bj-B
\auth{A.,Bj\"orner;F.,Brenti}
\book{Combinatorics of Coxeter Groups}
\publisher{Springer}
\Year{2005}

\nextref Dip-J86
\auth{R.,Dipper;G.,James}
\paper{Representations of Hecke algebras of general linear groups}
\journal{Proc. London Math. Soc.}
\Vol{52}
\Year{1986}
\Pages{20-52}

\nextref Du92
\auth{J.,Du}
\paper{The Green correspondence for the representations of Hecke algebras of type $A\sb{r-1}$}
\journal{Trans. Amer. Math. Soc.}
\Vol{329}
\Year{1992}
\Pages{273-287}

\nextref Geck-P
\auth{M.,Geck;G.,Pfeiffer}
\book{Characters of Finite Coxeter Groups and Iwahori-Hecke Algebras}
\publisher{Clarendon Press}
\Year{2000}

\nextref Gur90
\auth{D.I.,Gurevich}
\paper{Algebraic aspects of the quantum Yang-Baxter equation\inRus}
\journal{Algebra i Analiz}
\Vol{2:4}
\Year{1990}
\Pages{119-148};
\etransl{Leningrad Math. J.}
\Vol{2}
\Year{1991}
\Pages{801-828}

\nextref Hai97
\auth{P.H.,Hai}
\paper{Koszul property and Poincar\'e series of matrix bialgebras of type $A\sb n$}
\journal{J.~Algebra}
\Vol{192}
\Year{1997}
\Pages{734-748}

\nextref Hai02
\auth{P.H.,Hai}
\paper{Realizations of quantum hom-spaces, invariant theory, and quantum determinantal ideals}
\journal{J.~Algebra}
\Vol{248}
\Year{2002}
\Pages{50-84}

\nextref Hai05
\auth{P.H.,Hai}
\paper{On the representation categories of matrix quantum groups of type $A$}
\journal{Vietnam J. Math.}
\Vol{33}
\Year{2005}
\Pages{357-367}

\nextref Hay92
\auth{T.,Hayashi}
\paper{Quantum groups and quantum determinants}
\journal{J.~Algebra}
\Vol{152}
\Year{1992}
\Pages{146-165}

\nextref Hie93
\auth{J.,Hietarinta}
\paper{Solving the two-dimensional constant quantum Yang-Baxter equation}
\journal{J.~Math. Phys.}
\Vol{34}
\Year{1993}
\Pages{1725-1756}

\nextref Hla87
\auth{L.,Hlavat\'y}
\paper{Unusual solutions to the Yang-Baxter equation}
\journal{J.~Phys.~A}
\Vol{20}
\Year{1987}
\Pages{1661-1667}

\nextref Lar-T91
\auth{R.G.,Larson;J.,Towber}
\paper{Two dual classes of bialgebras related to the concepts of ``quantum group'' and ``quantum Lie algebra''}
\journal{Comm. Algebra}
\Vol{19}
\Year{1991}
\Pages{3295-3345}

\nextref Lyu87
\auth{V.V.,Lyubashenko}
\paper{Superanalysis and solutions to the triangles equation}
PhD thesis, Kiev, 1987.

\nextref Man
\auth{Yu.,Manin}
\book{Quantum Groups and Non-Commutative Geometry}
\publisher{CRM, Univ. Montr\'eal}
\Year{1988}

\nextref Mro14
\auth{C.,Mrozinski}
\paper{Quantum groups of $GL(2)$ representation type}
\journal{J.~Noncommut. Geom.}
\Vol{8}
\Year{2014}
\Pages{107-140}

\nextref Pol-P
\auth{A.,Polishchuk;L.,Positselski}
\book{Quadratic Algebras}
\publisher{Amer. Math. Soc.}
\Year{2005}

\nextref Pr70
\auth{S.B.,Priddy}
\paper{Koszul resolutions}
\journal{Trans. Amer. Math. Soc.}
\Vol{152}
\Year{1970}
\Pages{39-60}

\nextref Resh-TF89
\auth{N.Yu.,Reshetikhin;L.A.,Takhtajan;L.D.,Faddeev}
\paper{Quantization of Lie groups and Lie algebras\inRus}
\journal{Algebra i Analiz}
\Vol{1:1}
\Year{1989}
\Pages{178-206};
\etransl{Leningrad Math. J.}
\Vol{1}
\Year{1990}
\Pages{193-225}

\nextref Sch96
\auth{P.,Schauenburg}
\paper{Hopf bigalois extensions}
\journal{Comm. Algebra}
\Vol{24}
\Year{1996}
\Pages{3797-3825}

\endreferences
\bye